

\input amstex
\documentstyle{amsppt}

\magnification = \magstep 1

\vsize 7.5
in
\hsize 5.3
in

\def \ff {${\Cal F}$}
\def \cc {\Bbb C}
\def \zz {\Bbb Z}
\def \hak {^{\vee}}

\def \Odst {\ \ \ \ \ \ \ \ \ \ }

\def \sigc {\sigma\!\!\!^{^{\circ}}}

\let\nkr=\overline

\let\le=\leqslant
\let\ge=\geqslant

\let\tk=\triangledown

\bigskip\bigskip\medskip

\centerline {\bf FORMULAS FOR LAGRANGIAN AND}
\medskip
\centerline {\bf ORTHOGONAL DEGENERACY LOCI;}
\medskip
\centerline {\bf The $\widetilde Q$-Polynomials Approach
\footnote {\eightrm 1991 Mathematics Subject Classification
14M15, 14C17, 05E05}}
\medskip

\bigskip\bigskip\medskip
\centerline {\bf Piotr Pragacz\footnote{\eightrm Research carried out during
the author's stay at the Max-Planck-Institut f\"ur Mathematik as a fellow
of the Alexander von Humboldt Stiftung. While preparing this work the
author was partially supported by KBN grant No 2 P301002 05.
While revising the paper he was supported by the Max-Planck Institut
f\"ur Mathematik.}}
\centerline {\eightrm {Max-Planck Institut f\"ur Mathematik,
Gottfried-Claren Strasse 26, D-53225 Bonn, Germany.}}
\bigskip\medskip
\centerline {\bf Jan Ratajski\footnote{\eightrm Supported by KBN grant
No 2 P301002 05 and while revising the paper -- by the R. Bosch Foundation.}}
\centerline {\eightrm {Institute of Mathematics, Polish Academy
of Sciences, \'Sniadeckich 8, PL-00950 Warsaw, Poland.}}

\bigskip\medskip\smallskip

{\bf Contents}
\medskip

\ Introduction
\smallskip

1. Schubert subschemes and their desingularizations.
\smallskip

2. Isotropic Schubert calculus and the class of the diagonal.
\smallskip

3. Subbundles intersecting an $n$-subbundle in dim $\geqslant k$.
\smallskip

4. $\widetilde Q$-polynomials and their properties.
\smallskip

5. Divided differences and isotropic Gysin maps;

\noindent
\ \ \ \ \ \ \ \ orthogonality of $\widetilde Q$-polynomials.
\smallskip

6. Single Schubert condition.
\smallskip

7. Two Schubert conditions.
\smallskip

8. An operator proof of Proposition 3.1.
\smallskip

9. Main results in the generic case.
\smallskip

 \  Appendix A: Quaternionic Schubert calculus.
\smallskip
 \  Appendix B: Introduction to Schubert polynomials
\`a\ la polonaise.
\smallskip

 \  References

\bigskip\medskip

{\bf Introduction}
\medskip

In this paper we give formulas for the fundamental classes of
Schubert subschemes in Lagrangian and orthogonal Grassmannians
of maximal rank subbundles as well as some globalizations of them.
Our motivation to deal with this subject came essentially from 3 examples
where such degeneracy loci appear in algebraic geometry: \
$1^o$ The Brill-Noether loci for Prym varieties, as defined by
Welters [W],
$2^o$ \ The loci of curves with sufficiently many
theta characteristics, as considered by Harris [Har],
$3^o$ \ Some "higher" Brill-Noether loci in the moduli spaces
of higher rank vector bundles over curves, considered by Bertram
and Feinberg [B-F] and, independently, by Mukai [Mu].
\medskip

The common denominator of these 3 situations is a simple and beautiful
construction of Mumford [M].
With a vector bundle over a curve equipped with a nondegenerate quadratic form
with values in the sheaf of 1-differentials, Mumford associates an even
dimensional vector space endowed with a nondegenerate quadratic form
and 2 maximal isotropic subspaces such that the space of global sections
of the initial bundle is the intersection of the two isotropic subspaces.
A globalization of this construction allows one to present in a similar way the
varieties in $1^o$ and $2^o$ above as loci where two isotropic rank $n$
subbundles
of a certain rank $2n$ bundle equipped with a quadratic nondegenerate
form, intersect in dimension exceeding a given number. On the other hand,
the locus in $3^o$ admits locally this kind of presentation using
an appropriate symplectic form.
\bigskip

These varieties are particular cases of Schubert subschemes in
Lagrangian and
orthogonal Grassmannian bundles and their globalizations. The formulas
for such loci are the main theme of this paper. More specifically,
given a vector bundle $V$ on a variety $X$, endowed with a nondegenerate
symplectic or orthogonal form, we pick $E$ and $F_1\subset F_2\subset ...
\subset F_n=F$ - isotropic subbundles of $V$ ($rank \ E=n, rank \ F_i=i$)
and for a given sequence $1\le a_1<...<a_k\le n$ we look at the locus:
$$
D(a.):=\bigl\{\ x\in X | \ dim\bigl(E\cap F_{a_p}\bigr)_x\geqslant p,
 \ p=1,...,k \bigr\}.
$$
\medskip
We distinguish three cases:

1. Lagrangian: $rank \ V=2n$, the form is symplectic;

2. odd orthogonal: $rank \ V=2n+1$, the form is orthogonal;

3. even orthogonal: $rank \ V=2n$, the form is orthogonal.
\smallskip

(In the latter case the definition of $D(a.)$ must be slightly modified - see
Section 9.)
\medskip

Let us remark that the loci $D(a.)$ (for the Lagrangian case)
admit an important specialization to the loci introduced by
Ekedahl and Oort in the moduli space of abelian varieties with fixed
dimension and polarization, in characteristic $p$ (see, e.g. [O],
the references therein \  and [E-vG]). This comes from certain filtrations
on the de Rham cohomology defined with the help of the Frobenius-
and "Verschiebung"-maps.
The formulas of the present paper
are well suited to computations of the fundamental
classes of such loci in the Chow groups of the moduli spaces
- for details see a forthcoming
paper by T. Ekedahl and G. van der Geer [E-vG].
\bigskip

The goal of this paper is to give an algorithm for computing the fundamental
classes of $D(a.)$ as polynomials in the Chern classes of $E$ and $F_i$.
Formulas given here can be thought of as Lagrangian and orthogonal
analogs of the formulas due independently to Kempf-Laksov \cite {K-L} and
Lascoux \cite {L1} (notice, however, that the formulas given in \cite {K-L}
are proved under a weaker assumption of \ "expected" dimension).
\bigskip

The method for computing the fundamental class of a subscheme of a given
(smooth) scheme which we use here stems from a paper by the first
author \cite {P3, Sect.5}.
It depends on a desingularization of the subscheme in question and
the knowledge of the class of the diagonal of the ambient space.
It appears that the diagonals in the fibre products of Lagrangian
or orthogonal Grassmannian-  and flag bundles are not given as
the subschemes of zeros of sections of bundles over the corresponding
products. This makes an additional difficulty (e.g. in comparison with
[K-L]) which is overcomed here using again a result from [P3, Sect.5]
allowing to compute the class od the diagonal with the help of
an appropriate "orthogonality" property of Gysin maps.
\bigskip

To establish formulas for the classes of these diagonals, we use essentially
two tools.
The first one is Theorem 6.17  of \cite {P2} interpreting (cohomology dual to)
the classes of Schubert subvarieties in Lagrangian and orthogonal Grassmannians
as Schur's $Q$- and $P$-polynomials.
The importance of these polynomials to algebraic geometry
was illuminated by the first author in \cite {P1} and then developed
in \cite {P2}.
In fact in [P2, Sect.6], a variant of these polynomials was used
to give a full description of Schubert Calculus on Grassmannians of
maximal dimensional isotropic subspaces associated with a nondegenerate
symplectic and orthogonal form. These familes of symmetric polynomials
are called $\widetilde Q$- and $\widetilde P$-polynomials in the present
paper.
Perhaps the "orthogonality" proved in Theorem 5.23 is their
central property. This is, in fact, the second tool in our computation of the
classes of the diagonals in isotropic Grassmannian bundles which
allows us to apply the technique of [P3, Sect.5].
The results of \cite {P2, Sect.6}, recalled in Theorem 2.1 below,
are a natural source of the ubiquity of $\widetilde Q$- and $\widetilde
P$-polynomials in
various formulas of this paper. As a general rule, these are $\widetilde
Q$-polynomials
that appear in the Lagrangian case and $\widetilde P$-polynomials that appear
in the orthogonal cases.
\bigskip

In general, our approach gives an efficient algorithm for finding
formulas for Lagrangian and orthogonal Schubert subschemes. In several
cases, however, we are able to give compact expressions. At first, these
are the cases of one (i.e. $k=1$) and two Schubert conditions (the case of one
Schubert condition is usually referred to as a {\it special} Schubert
subscheme). The corresponding formulas are given in Section 6 and 7.
\medskip

The derivation of those formulas uses a formula for the
push-forward of $\widetilde Q$-polynomials
(Theorems 5.10, 5.14, 5.20) from
isotropic Grassmannian bundles. For instance, in the Lagrangian case,
$\pi : LG_nV \to X$ with the tautological subbundle $R$, the element
$\widetilde Q_IR\hak$ has a nonzero
image under $\pi_*$ only if each number $p$,
$1\le p \le n$, appears as a part of $I$ with an odd multiplicity $m_p$.
If the latter condition holds then
$$
\pi_* \widetilde Q_IR\hak
= \prod_{p=1}^n \bigl((-1)^p c_{2p}V\bigr)^{(m_p-1)/2}.
$$
\medskip

We also give formulas for
the push-forward of $S$-polynomials (Theorems 5.13, 5.15, 5.21)
from isotropic Grassmannian bundles. For example, in the Lagrangian case,
the element $s_IR\hak$ has a nonzero image under $\pi_*$ only if
the partition $I$ is of the form $2J+\rho_n$ for some partition $J$
(here, $\rho_n=(n,n-1,...,1)$).
If $I=2J+\rho_n$ then
$$
\pi_*s_IR\hak = s_J^{^{[2]}}V \quad,
$$
where the right hand side is defined as follows:
if $s_J= P(e.)$ is a unique presentation of $s_J$ as a polynomial
in the elementary symmetric functions $e_i$, $E -$ a vector bundle,
then $s_J^{[2]}(E) := P$ with $e_i$ replaced by $(-1)^ic_{2i}E, \
 i=1,2,\ldots$ \ .
\bigskip

Another case (corresponding to the Schubert condition $a.=(n-k+1,...,n)$)
that leads to compact formulas is the variety of maximal rank isotropic
subbundles which
intersect a fixed maximal rank isotropic subbundle
in dimension exceeding a given number (Proposition 3.2
and its analogs).
Thanks to the Cohen-Macaulayness of Schubert subschemes in isotropic
Grassmannians proved in [DC-L],
one gets globalizations of those formulas (as well as the other ones) to
more general loci.
For instance, the latter case $a.=(n-k+1,...,n)$
globalizes to the Mumford type locus discussed above where two maximal
rank isotropic subbundles $E$ and $F$ intersect in dimension exceeding
$k$, say.\footnote{It is mentioned in [F1,2] that the problem
of finding formulas for the classes in this case was posed originally
by Professor J. Harris several years ago.}
\bigskip

Our formulas (see Theorems 9.1, 9.5 and 9.6) are quadratic expressions in
$\widetilde Q$- and $\widetilde P$-polynomials of the subbundles.
More explicitly in the corresponding cases we have
{\parindent=25pt
\medskip
\item{1.} Lagrangian: \ \ \ \ \ \ \ \ \ \ $\sum \widetilde
Q_IE\hak\cdot \widetilde Q_{(k,k-1,\ldots,1)\smallsetminus I} F\hak$;
\medskip
\item{2.} odd orthogonal:   \ \ \ \ \ \ \ \ $\sum \widetilde
P_IE\hak\cdot \widetilde P_{(k,k-1,\ldots,1)\smallsetminus I} F\hak$;
\medskip
\item{3.} even orthogonal:  \ \ \ \ \ \ \ \ $\sum \widetilde
P_IE\hak\cdot \widetilde P_{(k-1,k-2,\ldots,1)\smallsetminus I} F\hak$;
\par}
\medskip
\noindent
where in 1. and 2. the sum is over all subsequences $I$ in
$(k,k-1,\ldots,1)$, in 3. the sum is over all subsequences $I$ in
$(k-1,k-2,\ldots,1)$ and $(k,k-1,\ldots,1)\smallsetminus I$ denotes
the strict partition whose parts complements the ones of $I$ in
$\{k,k-1,\ldots,1\}$.
\bigskip

Formula 3. has been recently used by C. De Concini and the first
named author in \cite {DC-P} to compute the fundamental classes of the
Brill-Noether loci $V^r$ for the Prym varieties (see \cite W), thus
solving a problem of Welters, left open since 1985. The formula of \
\cite {DC-P} asserts that if either $V^r$ is empty or of pure codimension
$r(r+1)/2$ in the Prym variety then its fundamental class in the
numerical equivalence ring, or its cohomology class is equal to
$$
2^{r(r-1)/2}\prod_{i=1}^r \bigl ((i-1)!/(2i-1)!\bigr) \ [\Xi]^{r(r+1)/2},
$$
where $\Xi$ is the theta divisor on the Prym variety.
\bigskip

The paper is organized as follows.

\smallskip

Section 1 contains definitions and properties of Schubert varieties
in Lagrangian and orthogonal Grassmannian bundles. Also, some
desingularizations of these varieties, used in later sections, are described.

\smallskip

Section 2 contains some recollections of Schubert calculus for Lagrangian
and orthogonal Grassmannians from [P2, Sect.6] and computation of the classes
of the diagonals in the Chow rings of Lagrangian and orthogonal
Grassmannian bundles. This computation relies on the Gysin maps technique
from [P3, Sect.5] and on the orthogonality theorem 5.23 which is proved
independently later.

\smallskip

Section 3 contains an explicit computation of Gysin maps needed to
determine the formulas for the fundamental
classes of Schubert varieties $\Omega(n-k+1,\ldots,n)$
parametrizing subbundles intersecting
an $n$-subbundle in dimension exceeding $k$.
This is done using an elementary Schubert Calculus--type technique
based on linear algebra.

\smallskip

In Section 4 we introduce a family of symmetric polynomials
called $\widetilde Q$-polynomials which is modelled on Schur's
$Q$-polynomials (but is different from the latter family).
These polynomials are the basic algebraic tools of the present paper.
We prove several elementary but useful properties of
$\widetilde Q$-polynomials and give some examples.

\smallskip

In Section 5 we establish some new algebraic properties of
of $\widetilde Q$-polynomials and $S$-polynomials; these are either
certain determinantal identities like Propositions 5.2 and 5.11, or the
computation of the values of these polynomials under some
divided differences operators. These algebraic results are then
interpreted using Gysin maps for Lagrangian and orthogonal
Grassmannian bundles. Perhaps the most important result of
this section is the "orthogonality" Theorem 5.23. This theorem, interpreted
geometrically (using a result of [P3, Sect.5]), gives us the
classes of the diagonals of Lagrangian and orthogonal Grassmannian
bundles which are crucial for our computations.

\smallskip

Sections 6., 7. and 8. have a supplementary character.
They contain some examples and a certain alternative (to the
content of the previous sections) way of computing.
Section 6 contains formulas for Schubert varieties defined by one Schubert
condition in Lagrangian and orthogonal cases. Section 7 contains
similar computations for two Schubert conditions
in the Lagrangian and odd orthogonal cases. Section 8 contains
another (purely algebraic) proof, using divided differences, of
Proposition 3.1 that describes the Gysin maps for some flag bundles.

\smallskip
In Section 9 we formulate previous results in the general setup
of degeneracy loci and give some examples. A special emphasis
is put on formulas answering J. Harris' problem concerning
the Mumford-type degeneracy loci described above.

\smallskip
In Appendix A  we collect a number of useful results about
Quaternionic Grassmannians. We use them to reprove some results
proved earlier using different methods and to show how some problems
concerning Grassmannians of nonmaximal Lagrangian subspaces can be
reduced to those of maximal Lagrangian subspaces; this sort of applications
we plan to develop elsewhere.
\smallskip

Finally, in Appendix B, we give an introduction to a theory of symplectic
Schubert polynomials which has grown up from the present work. This theory
(see [L-P-R]) seems to be well suited to the needs of algebraic geometry
because it generalizes in a natural way $\widetilde Q$-polynomials which
govern the Schubert calculus on Lagrangian Grassmannians.
\medskip

In Sections 2, 3, 5, 6, 7, 8 and 9 we work in the Chow rings; all
results therein, however, are equally valid in the cohomology rings.

\medskip

Some of the results of this paper were announced in \cite {P-R0}.
\medskip

The paper is a revised version of the Max-Planck Institut f\"ur Mathematik
Preprint MPI / 94 -132.
\bigskip\medskip

\underbar {Acknowgledgements}
\smallskip

We gratefully thank:
\smallskip

 - A. Collino and L. Tu for convincing
us several years ago about
the importance of orthogonal degeneracy loci in geometry by informing us
about the examples $1^o$ and $2^o$ mentioned above;
\smallskip

 - C. De Concini and G. van der Geer for their interest in this paper
and pointing out some defects in its previous version;
\smallskip

 - A. Lascoux and J.-Y. Thibon for valuable comments about $Q'$- and
$\widetilde Q$-polynomials;
\smallskip

 -  Professor W. Fulton for encouraging us to write the present paper;
\smallskip
\noindent
and

 - Professor F. Hirzebruch for some helpful comments
concerning quaternionic manifolds and drawing our attention
to the paper [Sl].
\medskip

We also thank the referee for pointing out several defects in the previous
version.
\bigskip

The preparation of this paper has been facilitated by the use of the
program system SYMMETRICA [K-K-L].
\medskip

The authors greatly benefited from the hospitality and stimulating
atmosphere of the Max-Planck Institut f\"ur Mathematik during the
preparation of this paper.

\bigskip\medskip

\underbar {Background}

Several results of this paper: e.g. Propositions 3.2, 3.4 and 3.6 as well
as their globalizations in Theorems 9.1,
9.5 and 9.6 were obtained already in Spring 1993 when
we tried to deduce formulas for the loci $D(a.)$ by combining the ideas
of the paper of Kempf and Laksov \cite {K-L} with the $Q$-polynomials
technique developed in \cite {P1,2}.
These results were announced together with outlines of their proofs in
\cite {P-R0}.
\smallskip

In summer '93, we received an e-mail message from Professor W. Fulton
informing us about his (independent) work on the same subject and
announcing another expressions for the loci considered in Proposition 3.2, 3.4
and 3.6 of the present paper. Responding, we informed Professor Fulton
about our results of [P-R0] mentioned above.
In February '94 we
obtained from Professor W. Fulton his preprints
[F1,2]
containing details of his e-mail announcement. Both the form of the formulas
obtained as well as the approach used in [F1,2] are totally different from
the content of our work and just a simple comparison of the results of [F1,2]
with ours leads to very nontrivial new identities which are interesting in
themselves.
It would be desirable to develop, in a systematic way, the comparison
of formulas given in \cite{F1,2} from
one side with those
in the present paper and \cite {P-R0} - from the other one.

\bigskip\bigskip

\underbar {Conventions}
\bigskip

Partitions are weakly decreasing sequences of positive integers (as in
[Mcd1] and
are denoted by capital Roman letters (as in [L-S1]).
We identify partitions with their
Ferrers' diagrams visualized as in [L-S1].
The relation "$\subset$" for partitions is induced from
that for diagrams.
\medskip

For a given partition $I=(i_1,i_2,\ldots)$ we denote by $|I|$ (the {\it weight}
of $I$) the
partitioned number (i.e. the sum of all parts of $I$) and by $l(I)$
(the {\it length} of $I$)
the number of nonzero parts of $I$.
Moreover, $I^{\sim}$ denotes the dual partition of $I$, i.e. $I^{\sim}=
(j_1,j_2,\dots)$ where $j_p=card\{h |i_h \ge p\}$, and $(i)^k$ - the partition
$(i,\dots,i)$ \ ($k$-times).
\medskip

Given sequences \ $I=(i_1,i_2,\ldots)$ and $J=(j_1,j_2,\ldots)$ we denote
by $I\pm J$ the sequence $(i_1\pm j_1,i_2\pm j_2,\ldots)$.
\medskip

By {\it strict} partitions we mean
those whose (positive) parts are all different.
\medskip

In this paper, we denote by $s_i(E)$ the
complete symmetric polynomial of degree $i$ with variables
specialized to the Chern
roots of a vector bundle $E$.
\medskip

The reader should be careful with our notion of $\widetilde Q$-polynomials
here. Namely,
since we are mainly interested in the polynomials in the Chern
classes of vector bundles, we introduce $\widetilde Q$-polynomials given
by the Pfaffian of an antisymmetric matrix whose entries
are quadratic expressions in the elementary symmetric polynomials
rather than in the "one row" Schur's $Q$-polynomials. Therefore these
polynomials are different from the original Schur's $Q$-polynomials.
Note that nonzero $\widetilde Q$-polynomials $\widetilde Q_I(x_1,\ldots,x_n)$
are indexed by "usual" partitions $I$
but the parts of these partitions cannot exceed the number of variables;
on the contrary, nonzero Schur's $Q$-polynomials $Q_I(x_1,\ldots,x_n)$
are indexed by strict partitions $I$
only but the parts of these partitions can be bigger than the number of
variables.
\medskip

Also, the specialization of $\widetilde Q_I(x_1,\dots,x_n)$
with $(x_i)$ equal to the sequence of the Chern roots of
a rank $n$ vector bundle $E$, denoted here
- accordingly - by $\widetilde Q_IE$, is a different cohomology
class
than the one associated with $E$ in [P1] and [P2, Sect.3 and 5],
and denoted by $Q_IE$ therein. (Notice, however, that the
$\widetilde Q$-polynomials appeared already in an
implicit way in [P2, Sect.6].)
The reader should make a proper distinction between Schur's $Q$-polynomials
and $\widetilde Q$-polynomials that are mainly used in the present paper.
\medskip

For a vector bundle $V$, by $G_nV$ we denote the usual Grassmannian bundle
parametrizing rank $n$ subbundles of $V$.
Moreover, $\Bbb P(V)=G_1V$.
We follow mostly [F] for the terminology in algebraic geometry.
In many situations when the notation starts to be too cumbersome, we
omit some pullback-indices of the induced vector bundles.
\medskip

A good reference for "changes of alphabets" in the $\lambda$-ring sense
is [L-S1].

\bigskip\bigskip

\centerline {\bf 1. Schubert subschemes and their desingularizations}
\bigskip

We start with the Lagrangian case. Let $K$ be an arbitrary ground field.

Assume that $V$ is a rank $2n$ vector bundle over a smooth scheme $X$ over $K$
equipped
with a nondegenerate symplectic form.
Moreover, assume that a flag $V.$ :
$V_1 \subset V_2 \subset \ldots \subset V_n$
of Lagrangian (i.e. isotropic) subbundles w.r.t. this form is fixed,
with $rank \ V_i = i$.
Let $\pi: {L}G_n(V)\to X$ denote
the Grassmannian bundle parametrizing Lagrangian rank $n$
subbundles of $V$.
$G=LG_n(V)$ is endowed with the tautological Lagrangian
bundle $R \subset V_G$.
Given a sequence \ $a_.=(1 \leqslant a_1 < \ldots < a_k \leqslant n)$ \ we
consider in $G$ a closed subset:
$$
\Omega (a_.) =\Omega (a_.;V.) =
\{g\in G | \ dim\bigl(R\cap V_{a_i}\bigr)_g\geqslant i,\ i=1,\ldots,k\}.
$$
The locus $\Omega (a_.)$, called a {\it Schubert subscheme}
is endowed with a reduced scheme structure induced from
the reduced one of the corresponding Schubert subscheme
in the Grassmannian $G_nV$ -
this is discussed in detail , e.g., in [L-Se].

\bigskip
The following desingularization of \ $\Omega = \Omega(a_.)$ \
should be thought of as a Lagrangian
analogue of the construction used in [K-L].
Let \ ${\Cal F} = {\Cal F}(a_.) = {\Cal F}(V_{a_1} \subset \ldots \subset
V_{a_k})$ \ be the scheme parametrizing flags \ $A_1 \subset A_2 \subset
\ldots \subset A_k \subset A_{k+1}$ \  such that \ $rank \ A_i = i$ \ and \
$A_i \subset V_{a_i}$ \ for \ $i = 1,\ldots,k$ ; \ $rank \ A_{k+1} = n$ \ and
\ $A_{k+1}$ \ is Lagrangian. ${\Cal F}$ is endowed with the tautological flag
\ $D_1 \subset D_2 \ldots \subset D_k \subset D_{k+1}$, \ where \
$rank \ D_i = i , \ i=1,\ldots,k$ \ and \ $rank \ D_{k+1} = n$~.
We will write $D$ instead of \ $D_{k+1}$.
\smallskip

We have a fibre square:
$$
\CD
G \times _{_X} {\Cal F}  @> p_2 >>  {\Cal F} \\
@V p_1 VV                         @V p VV \\
G @> >>X
\endCD
$$
\smallskip

\noindent
Let $\alpha:{\Cal F}\to G$ be the map defined by:
$(A_1\subset A_2\subset\ldots\subset A_{k+1})\mapsto A_{k+1}$,
in other words $\alpha$ is a "classifying map" such that $\alpha^*R=D$.
It is easily verified that $\alpha$ maps ${\Cal F}$ onto $\Omega$ and
$\alpha$ is an isomorphism over the open subset of $\Omega$ parametrizing
rank $n$ Lagrangian subbundles $A$ of $V$ such that $rank(A\cap V_{a_i})=i$,
$i=1,\ldots,k$.
Moreover, $\alpha$ induces a section $s$ of $p_2$.
Set $Z:=s({\Cal F})\subset G\times_{_X}{\Cal F}$.
Alternatively, we can describe $Z$ as $(1\times\alpha)^{-1}(\Delta)$
where $\Delta$ is the diagonal in $G\times_{_X}G$.
The map $p_1$ restricted to $Z$ is a desingularization of $\Omega$.
Therefore $[\Omega]=(p_1)_*([Z])$.
On the other hand, $[Z]=(1\times\alpha)^*([\Delta])$ (see [K-L, Lemma 9]).
Note that \ff \ is obtained as a composition
of the following flag- and Grassmannian bundles.
Let $Fl= Fl(a.)= Fl(V_{a_1} \subset \ldots \subset V_{a_k})$ be the "usual"
flag bundle parametrizing flags $A_1 \subset \ldots \subset A_k$
where $rank\ A_i = i$ and $A_i \subset V_{a_i}, \ i=1,\ldots,k$.
Let $C_1 \subset \ldots \subset C_k$ be the tautological flag on
$Fl$.
We will write $C$ instead of $C_k$. Then $\Cal F$ is the Lagrangian
Grassmannian bundle $LG_{n-k}(C^{\perp}/C)$ over $Fl$, where
$C^{\perp}$ is the subbundle of $V_{Fl}$ consisting of all $v$
that are orthogonal to $C$ w.r.t. the given symplectic form. Note that
$C\subset C^{\perp}$ because $C$ is Lagrangian, $rank(C^{\perp}/C)
=2(n-k)$ and the vector bundle $C^{\perp}/C$ is endowed with a
nondegenerate symplectic form induced from the one on $V$.
Of course the tautological Lagrangian rank $n-k$ subbundle on
$LG_{n-k}(C^{\perp}/C)$ is identified with $D/C_{\Cal F}$.
In other words, \ff \ is a composition of a flag bundle (with the fiber being
$Fl(K^{a_1} \subset \ldots \subset K^{a_k})$) and
a Lagrangian Grassmanian bundle (with the fiber being
${L}G_{n-k}(K^{2(n-k)}$).
In particular, $$dim\ \Omega=dim\ {\Cal F}=dim\ Z=\sum\limits_{i=1}^k(a_i-i)+
(n-k)(n-k+1)/2 + dim\ X.$$
\medskip

The following particular cases will be treated in a detailed way in this paper:
$a_. = (n-k+1,\ n-k+2,\ldots, n)$ (then $\Omega(a_.)$
parametrizes Lagrangian rank $n$ subbundles $L$ of $V$ such that
\ $rank(L\cap V_n) \geqslant k$~); $a.=(n+1-i)$ , i.e. $k=1$ ;
and $a.=(n+1-i,n+1-j)$, i.e. $k=2$.
\bigskip

Now consider the odd orthogonal case.
Let $K$ be a ground field of characteristic different from $2$.
Assume, that $V$ is a rank $2n+1$ vector bundle over a smooth scheme $X$
over $K$
equipped with a nondegenerate orthogonal form.
We assume throughout this paper that the form restricts to a hyperbolic
form on each fiber (i.e. each fiber has an $n$-dimensional
isotropic subspace; if $K$ is algebraically closed, this is automatically
satisfied.)
Let $OG_nV$ be the Grassmannian bundle parametrizing
rank $n$ isotropic subbundles of $V$. Whenever, in this paper, we speak
about $OG_nV$, we assume that there exists a rank $n$ isotropic
subbundle in $V$.
All definitions, notions and notation concerning Schubert subschemes
and their desingularizations
are used mutatis mutandis (just instead of
"symplectic" use "orthogonal" and instead of
"Lagrangian" use "isotropic"). The formula for the dimension of
$\Omega(a.)$ in the odd orthogonal case is the same as in the Lagrangian case.
Of course, \ff \ is now a composition of the same flag bundle $Fl$ and the odd
orthogonal Grassmannian bundle ${O}G_{n-k}(C^{\perp}/C)$, where $C$ is the
rank $k$ tautological subbundle on $Fl$.

\bigskip

Assume now that $V$ is a rank $2n$ vector bundle over a smooth
connected scheme $X$
over a field $K$ of characteristic different from $2$ equipped
with a nondegenerate orthogonal form.
We assume throughout this paper that there exists an isotropic rank $n$
subbundle of $V$.
The scheme parametrizing
isotropic rank $n$ subbundles of $V$ breaks up into two
connected components denoted $OG_n'V$ and $OG_n''V$.
Let $V_n$ be a rank $n$ isotropic subbundle of $V$ fixed once and for all.
Then $OG_n'V$ (resp. $OG_n''V$) parametrizes rank $n$ isotropic
subbundles $E\subset V$ such that $dim(E\cap V_n)_x\equiv n (mod \ 2)$
 \ (resp. $dim(E\cap V_n)_x\equiv n+1 (mod \ 2)$ ) for every $x\in X$.
 Write $G':= OG_n'V$ and $G'':= OG_n''V$.
Two isotropic rank $n$ subbundles are in the same component iff they
intersect fiberwise in dimension congruent to $n$ modulo $2$.

\bigskip

Let $V.: V_1\subset V_2\subset\ldots\subset V_n$ be a flag of
isotropic subbundles of $V$ with $rank\ V_i=i$.
Given a sequence $a_.=(1\leqslant a_1<\ldots<a_k\leqslant n)$ such that
$k\equiv n\ (mod\ 2)$, we consider in $G'$ a Schubert subvariety:
$$
\Omega(a_.)=\Omega(a_.;V_.)=\Bigl\{g\in G' | \ dim(R\cap V_{a_i})_{_g}
\geqslant i,\ i=1,\ldots,k\Bigr\}
$$
($R\subset V_{G'}$ is here the tautological bundle).
Similarly, given a sequence $a_.=(1\leqslant a_1<\ldots<a_k\leqslant n)$
such that $k\equiv n+1\ (mod\ 2)$, we consider in $G''$ a Schubert subvariety
$$
\Omega(a_.)=\Omega(a_.;V_.)=\Bigl\{g\in G'' | \ dim(R\cap V_{a_i})_{_g}
\geqslant i,\ i=1,\ldots,k\Bigr\}.
$$
(Over a point, say, the interiors of the $\Omega(a.)$'s form a
cellular decomposition of $G'$ and respectively $G''$.)
Here, the definition of the scheme structure is more delicate than in the
previous two cases (roughly speaking, instead of minors one should use
the Pfaffians of the "coordinate" antisymmetric matrix of $G'$ and
$G''$). We refer the reader
for details to \cite {L-Se} and references therein.

\medskip

The Schubert subvarieties
$\Omega(a_.)$ in $G'$ (resp. $\Omega(a_.)$ in $G''$) are desingularized
using the same construction as above but instead of the scheme \ff\ one must
now use the following scheme \ff$'$ (resp. \ff$''$).
Let ${\Cal F}'={\Cal F}'(a_.)={\Cal F}'(V_{a_1}\subset\ldots\subset V_{a_k})$
be a scheme parametrizing flags
$A_1\subset A_2\subset\ldots\subset A_k\subset A_{k+1}$ such that $rank\ A_i=i$
and $A_i\subset V_{a_i}$ for $i=1,\ldots,k$; \ $rank\ A_{k+1}=n$, $A_{k+1}$
is isotropic and $rank(A_{k+1}\cap V_n)_x \equiv n\ (mod\ 2)$ for any
$x\in X$.
The definition of ${\Cal F}''={\Cal F}''(a_.)$ is the same with exception
of the last condition now replaced by:
$rank(A_{k+1}\cap V_n)_x \equiv n+1\ (mod\ 2)$ for any $x\in X$.
Let $p':\Cal F \to X$
(resp. $p'':\Cal F \to X$) denote the projection maps.
Of course, \ff$'$ (resp. \ff$''$) now is a composition of the same flag
bundle $Fl$ and the even orthogonal Grassmannian bundle ${O}G_{n-k}'
(C^{\perp}/C)$
(resp. ${O}G_{n-k}''(C^{\perp}/C)$), where $C$ is the rank $k$ tautological
subbundle on $Fl$ and $V_n/C$ is the rank $n-k$ isotropic bundle used
to define ${O}G_{n-k}'(C^{\perp}/C)$
and ${O}G_{n-k}''(C^{\perp}/C)$.
\medskip

The formula for dimension now is different:
$$
dim\ {\Cal F}'=dim\ {\Cal F}''=\sum\limits^k_{i=1}(a_i-i)+(n-k)(n-k-1)/2 +
dim\ X.
$$
\medskip

We finish this section with the following lemma which will be of
constant use in this paper.
\medskip

\proclaim{\bf Lemma 1.1}
Consider cases 1., 2., 3. of a vector bundle endowed with a
nondegenerate form $\Phi$ that are specified in the Introduction.
Let $C\subset V$ be an isotropic subbundle and $C^{\perp}$ be the
subbundle of $V$ consisting of all $v\in V$ such that $\Phi(v,c)=0$
for any $c\in C$.

1. Then one has an exact sequence
$$
\CD
0\longrightarrow C^{\perp}\longrightarrow V @>\phi >> C\hak
\longrightarrow 0
\endCD
$$
where the map $\phi$ is defined by \ $v\mapsto \Phi (v,-)$.
In particular, in the Grothendieck group, $[V]=[C^{\perp}]+[C\hak]$,
 $[C^{\perp}/C]=[V]-[C]-[C\hak]$ and the Chern classes of $C^{\perp}/C$
are the same as the ones of $[V]-[C\oplus C\hak]$.
\smallskip

2. Assume now that $C$ is a maximal isotropic subbundle of $V$. Then
in cases 1. and 3. we have $C=C^{\perp}$ and $c_.(V)=c_.(C\oplus C\hak)$;
in case 2. one has $rank(C^{\perp}/C)=1$ and
$2c_.(V)=2c_.(C\oplus C\hak)$.

\endproclaim
\medskip

The latter equality of assertion 2 in case 2.
follows from the fact that the form
$\Phi$ induces an isomorphism $(C^{\perp}/C)^{\otimes 2}\cong \Cal O_X$.
This assertion will be used in the proof of Theorem 5.14 and 5.15
and is well suited for this purpose because of the appearance of the factor
"$2^n$" on the right hand side of the formulas of the theorems.

\bigskip\medskip

\centerline {\bf 2. Isotropic Schubert Calculus and the class of the diagonal}
\bigskip

Let us first recall the following result on Lagrangian and orthogonal Schubert
Calculus from [P1,2].
\smallskip

We need two families of polynomials in the Chern classes of a vector
bundle $E$ over a smooth variety $X$.
Their construction is inspired by I. Schur's paper
\cite {S}.
The both families will be indexed by partitions (i.e. by sequences
$I=(i_1\ge \ldots \ge i_k\ge 0)$ of integers).
Set, in the Chow ring $A^*(X)$ of $X$, for $i\ge j\ge 0$:
$$
\widetilde Q_{i,j}E:=c_iE\cdot c_jE+2\sum\limits^j_{p=1}(-1)^pc_{i+p}E\cdot
c_{j-p}E,
$$
so, in particular $\widetilde Q_iE:=\widetilde Q_{i,0}E=c_iE$  \ for $i\ge 0$.
In general, for a partition $I=(i_1,\ldots,i_k)$, $k-$even
(by putting $i_k=0$ if necessary), we set in $A^*(G)$:
$$
\widetilde Q_IE:=Pf\Bigl( \widetilde Q_{i_p,i_q}E \Bigr)_{_{1\leqslant
p<q\leqslant k}},
$$
where $Pf$ means the Pfaffian of
the given antisymmetric matrix.
For the definition and basic properties of Pfaffians we send interested
readers to [A], [Bou] and [B-E, Chap.2]. Also,
we refer the reader to the beginning
of Section 4 for an alternative recurrent definition of $\widetilde Q_IE$:
just replace the polynomial $\widetilde Q_I(X_n)$ from Section 4 by
the element $\widetilde Q_IE$.
\smallskip

The member of the second family, associated with a partition $I$,
is defined by
$$
\widetilde P_IE:=2^{-l(I)}\widetilde Q_IE.
$$
Observe that in particular $\widetilde P_iE=c_iE/2$ (so here
we must assume that $c_iE$ is divisible by $2$),
and
$$
\widetilde P_{i,j}E=\widetilde P_iE\cdot \widetilde
P_jE+2\sum\limits^{j-1}_{p=1}(-1)^p\widetilde P_{i+p}E\cdot \widetilde P_{j-p}E
+(-1)^j\widetilde P_{i+j}E.
$$

It should be emphasized that $\widetilde Q$- and $\widetilde P$-polynomials
are especially
important and useful for isotropic (sub)bundles.
\smallskip

The following result from \cite{P1, (8.7)} and \cite{P2, Sect.6},
gives a basic geometric interpretation of $\widetilde Q$-
and $\widetilde P$-polynomials.
(It can be interpreted as a Giambelli-type formula for isotropic Grassmannians;
recall that a Pieri-type formula for these Grassmannians was given in [H-B] -
consult also [P-R1] for a simple proof of the latter result.)

\proclaim{\bf Theorem 2.1} {\rm [P2, Sect.6]}
{\parindent=25pt
  \item{(i) \ \ } Let $V$ be a $2n-$dimensional vector space over
a field $K$ endowed with a nondegenerate symplectic form.
Then, one has in $A^*({L}G_nV)$,
$$
\bigl[\Omega(a_.)\bigr]=\widetilde Q_IR\hak\ ,
$$
where $R$ is the tautological subbundle on ${L}G_nV$ and
$i_p=n+1-a_p$, $p=1,\ldots,k$.
  \item{(ii) \ } Let $V$ be a $(2n+1)-$dimensional vector space over a field
$K$
of $char.\ne 2$ endowed with a nondegenerate orthogonal form.
Then, one has in $A^*({O}G_nV)$,
$$
\bigl[\Omega(a_.)\bigr]= \widetilde P_IR\hak\ ,
$$
where $R$ is the tautological subbundle on ${O}G_nV$ and
$i_p=n+1-a_p$, $p=1,\ldots,k$.
  \item{(iii) } Let $V$ be a $2n-$dimensional vector space over a field $K$
of $char.\ne 2$ endowed with a nondegenerate orthogonal form.
Then one has in $A^*({O}G_n'V)$ (resp. $A^*({O}G_n''V)$),
$$
\bigl[\Omega(a_.)\bigr]=\widetilde P_IR\hak\ ,
$$
where $R$ is the tautological subbundle on ${O}G_n'V$
(resp. ${O}G_n''V$) and $i_p=n-a_p$, $p=1,...,k$.
( Notice that the indexing family of $I$'s runs here over all
strict partitions contained in $\rho_{n-1}$. )
  \item{}
Observe that by Lemma 1.1, $R\hak$ is the tautological quotient
bundle on $LG_nV$, $OG_n'V$ and $OG_n''V$. Moreover, the Chern classes
of the tautological quotient bundle on $OG_nV$ and $R\hak$ are equal.

\par}
\endproclaim
Note that this result has been reproved recently by S. Billey and M. Haiman
in [B-H].

\bigskip

Assume now that $V$ is a vector bundle over a smooth variety $X$ and
$V.$ is a flag of isotropic bundles on $X$. Then, using Noetherian
induction, one shows that $\{\widetilde Q_IR\hak \}_{I\subset \rho_n}$,
$\{\widetilde P_IR\hak \}_{I\subset \rho_n}$ and $\{\widetilde P_IR\hak
\}_{I\subset \rho_{n-1}}$
are $A^*(X)$-bases respectively of $A^*({L}G_nV)$, $A^*({O}G_nV)$ and
$A^*(OG'_nV)$ (resp. $A^*(OG''_nV)$). Moreover, there is an
expression for $\Omega(a.;V.)$ as a polynomial in the Chern classes
of $R\hak$ and $V_i$. (This follows, e.g., from the existence of
desingularizations given in Section 1 and formulas for Gysin
push forwards -- for "usual" flag bundles they are obtained by iterating
a well known projective bundle case;
for isotropic Grassmannian bundles, they are given for the first
time in Section 5 of the present paper). Then the maximal degree term
in $c.(R\hak)$ of this expression, in respective
cases (i), (ii), (iii),  coincides with that in Theorem 2.1. We will
call it {\it the dominant term} (w.r.t. $R$).

\medskip
Let \ $G_1,G_2$ \ be two copies of the Lagrangian Grassmannian
bundle $LG_nV$ over a smooth variety $X$,
equipped with the tautological subbundles
$R_1$ and $R_2$. Write \ $GG := G_1 \times_{_X} G_2$~.
Consider the following diagonal
$$
\Delta = \Bigl\{ (g_1,g_2) \in GG | \
\bigl( (R_1)_{_{GG}} \bigr)_{_{(g_1,g_2)}}
= \bigl( (R_2)_{_{GG}} \bigr)_{_{(g_1,g_2)}} \Bigr\}.
$$

Our goal is to write down a formula for the class of this
diagonal. We first record:

\proclaim{\bf Lemma 2.2}
Let $G$ be a smooth complete variety such that the "$\times$-map"
{\rm(cf. [F, end of Sect.1])}
gives an isomorphism \ $A^*(G\times G) \cong A^*(G) \otimes A^*(G)$~.
Assume that there exists a family \ $\{b_{\alpha}\}$~, \
$b_{\alpha} \in A^{n_{\alpha}}(G)$~, such that \
$A^*(G) = \oplus \zz b_{\alpha}$~, and for every $\alpha$ there is a unique
$\alpha'$ such that \ $n_{\alpha}+n_{\alpha'}= dim \ G$ \ and \
$\int_X b_{\alpha}\cdot b_{\alpha'} \ne 0$~.
Suppose \ $\int_X b_{\alpha}\cdot b_{\alpha'} = 1$~.
Then the class $[\Delta]$ in \ $A^*(G\times G)$ \ is given by \
$\sum\limits_{\alpha} b_{\alpha}\times b_{\alpha'}$~.
\endproclaim

\demo{Proof}
It follows from the assumptions that in \ $A^*(G\times G)~, \
[\Delta] = \sum m_{\alpha \beta} b_{\alpha}\times b_{\beta}$~,
for some integers \ $m_{\alpha \beta}$ \ and \ $n_{\alpha}+n_{\beta}=dim \ G$
\ for all pairs \ $(\alpha,\beta)$ \ indexing the sum.
We have by a standard property of intersection theory for $g, h \in A^*(G)$
$$
\int_{X\times X} [\Delta] \cdot (g \times h) =
\int_X g \cdot h.
$$
Hence the coefficients \ $m_{\alpha\beta}$ \ satisfy:
$$
m_{\alpha\beta} = \int_{X\times X}[\Delta] \cdot (b_{\alpha'}
\times b_{\beta'}) =
\int_X  b_{\alpha'}\cdot b_{\beta'}.
$$
The latter expression, according to our assumption is not zero only if
\ $\alpha' = (\beta')'$ \ i.e. \ $\beta = \alpha'$, when it equals 1.
This proves the lemma.
\qed
\enddemo
\medskip

For a given positive integer $k$, put \ $\rho_k = (k,k-1,\ldots,2,1)$~.
For a strict partition \ $I\subset \rho_k$ \ (i.e.
$i_1\leqslant k,\  i_2\leqslant k-1, \ldots$) we denote by
$\rho_k\smallsetminus I$  the strict partition whose parts complement the parts
of $I$ in the set
$\{k,k-1,\ldots,2,1\}$.
\smallskip

The Lagrangian Grassmannian (over a point, say) satisfies the assumptions
of the lemma with \ $\bigl\{ \widetilde Q_IR\hak\bigr\}_{strict \ I \subset
\rho_n}$
playing the role of $\{b_{\alpha}\}$~ and for $\alpha=I$ we have $\alpha'
=\rho_n\smallsetminus I$.
This is a direct consequence the existence of a well known cellular
decomposition of such a Grassmannian into Schubert cells and the results of
[P2] recalled in Theorem 2.1(i) together with a description of Poincar\'e
duality in $A^*(LG_nV)$ from loc.cit.
Thus in this situation we get by the lemma:

\proclaim{\bf Lemma 2.3}
The class of the diagonal $\Delta$ of the Lagrangian Grassmannian equals
$$
[\Delta] = \sum \widetilde Q_I(R_1\hak) \times \widetilde
Q_{\rho_n\smallsetminus I}(R_2\hak) ,
$$
the sum over all strict  \ $I \subset \rho_n$~.
\endproclaim

We now want to show that the same formula holds true for an arbitrary
smooth base space $X$ of a vector bundle $V$. Our argument is based on the
following result expressing the class of the (relative) diagonal in
terms of Gysin maps. This result was due to the first
author in [P3, Sect.5]
and is accompanied here by its proof for the reader's convenience.

\proclaim{\bf Lemma 2.4} \ {\rm [P3]}
 \ Let $\pi:G\to X$ be a proper morphism of smooth varieties such that
$\pi^*$ makes $A^*(G)$ a free $A^*(X)$-module,
$A^*(G)=\oplus_{\alpha \in \Lambda} A^*(X)\cdot b_{\alpha}$,
where $b_{\alpha}\in A^{n_{\alpha}}(G)$ and
for any $\alpha$ there is a unique $\alpha'$ such that
$n_{\alpha}+n_{\alpha'}=dim\ G - dim\ X$ and
$\pi_*(b_{\alpha}\cdot b_{\alpha'})\ne 0$;
suppose $\pi_*(b_{\alpha}\cdot b_{\alpha'})=1$.
Moreover, denoting by $p_i:G\times_{_X}G\to G$ ($i=1,2$) the projections,
assume that, for a smooth $G\times_{_X}G$, the homomorphism
$A^*(G)\otimes_{_{A^*(X)}}A^*(G)\to
A^*(G\times_{_X}G)$, defined by $g\otimes h\mapsto p_1^*g\cdot p_2^*h$, is an
isomorphism. Then
\smallskip

(i) The class of the diagonal $\Delta$ in $G\times_{_X}G$ equals
 \ $\big[\Delta\big] = \sum_{\alpha,\beta} m_{\alpha\beta}
b_{\alpha}\otimes b_{\beta},$ \ where,

for any $\alpha,\beta$, \ $m_{\alpha\beta}=P_{\alpha\beta}
(\{\pi_*(b_{\kappa} \cdot
b_{\lambda})\})$
for some polynomial $P_{\alpha\beta}\in \Bbb Z[\{x_{\kappa\lambda}\}]$.
\smallskip

(ii) If \  $\pi_* (b_\alpha \cdot b_\beta) \neq 0$ iff $\beta=\alpha'$,
then the class of the diagonal $\Delta \subset G \times_{_X} G$ equals
 \ $\big[\Delta\big] = \sum_{\alpha} b_{\alpha}\otimes b_{\alpha'}$.
\endproclaim
\smallskip

\demo{Proof}
Denote by $\delta:G\to G\times_{_X} G$, $\delta':G\to G\times_{_K}G$ (the
Cartesian
product)
the diagonal embeddings and by
$\gamma$ the morphism $\pi\times_{_X}\pi:G\times_{_X} G\to X$.
For $g,h\in A^*(G)$ we have
$$
\pi_*(g\cdot h) = \pi_*\Bigl((\delta')^*(g\times h)\Bigr) =
 \pi_*\Bigl(\delta^*(g\otimes h)\Bigr) =
\gamma_*\delta_*\Bigl(\delta^*(g\otimes h)\Bigr)
= \gamma_*\Bigl([\Delta]\cdot(g\otimes h)\Bigr),
$$
using $\pi=\gamma\circ\delta$ and standard properties of intersection theory
([F]). Hence, writing
$[\Delta] = \sum m_{\mu\nu}b_{\mu}\otimes b_{\nu}$,
we get
$$
\aligned
   \pi_*(b_{\alpha}\cdot b_{\beta})
&= \gamma_*\bigl([\Delta]\cdot(b_{\alpha}\otimes b_{\beta}\bigr)=
(\pi_*\otimes\pi_*)\Bigl(\bigl(\sum m_{\mu\nu}b_{\mu}\otimes
b_{\nu}\bigr)\cdot(b_{\alpha}\otimes b_{\beta})\Bigr) \\
&= \sum_{\mu, \nu}
 m_{\mu\nu}\pi_*(b_{\mu}\cdot b_{\alpha})
\cdot \pi_*(b_{\nu}\cdot b_{\beta}).
\endaligned
\tag *
$$

\medskip
\noindent
(i) By the assumption and (*) we get
$$
m_{\alpha\beta}=\pi_*(b_{\alpha'} \cdot b_{\beta'})-
\sum_{\mu \ne \alpha, \nu \ne \beta}
 m_{\mu\nu}\pi_*(b_{\mu}\cdot b_{\alpha'})
\cdot \pi_*(b_{\nu}\cdot b_{\beta'}).
\tag**
$$
where the degree of $m_{\mu\nu} \in A^*(X)$ such that
$\mu \ne \alpha$ or $\nu \ne \beta$ and
$\pi_*(b_{\mu}\cdot b_{\alpha'})
\cdot \pi_*(b_{\nu}\cdot b_{\beta'})\ne 0$,
is smaller than the degree
of $m_{\alpha\beta}$. The assertion now follows by induction
on the degree of $m_{\alpha\beta}$.
\smallskip

\noindent
(ii) By virtue of the assumption, Equation (**) now reads \
$\pi_*(b_{\alpha'} \cdot b_{\beta'})= m_{\alpha\beta}$ \  and immediately
implies the assertion.
\qed
\enddemo
\smallskip

Let us remark that in [DC-P, Proposition 2], where a weaker variant of
Theorem 9.6 of the present paper is used, already assertion (i)
(plus results of Section 5)
are sufficient to conclude the proof.

\bigskip

Let now $G=LG_nV\to X$ denote a Lagrangian Grassmannian bundle.
\smallskip
\proclaim{\bf Proposition 2.5}
The class of the diagonal of the Lagrangian Grassmannian bundle in
$A^*(G\times_{_X}G)$
equals
$$
[\Delta] = \sum \widetilde Q_I(R_1\hak)_{GG} \cdot
\widetilde Q_{\rho_n\smallsetminus I}(R_2\hak)_{GG} ,
$$
the sum over all strict $I\subset \rho_n$, $GG=G\times_{_X}G$
and $R_i$, $i=1,2$, are the tautological (sub)bundles on the corresponding
factors.
\endproclaim
\demo{Proof}
The assertion follows from Lemma 2.4(ii) applied to
$b_I=\widetilde Q_IR\hak$ ($I$ strict $\subset \rho_n$)
and Theorem 5.23(i) which
will be proved (independently) later.
\qed
\enddemo
\medskip

\proclaim{\bf Corollary 2.6}
With the notation of Section 1 and $G\Cal F:=G\times_{_X} \Cal F$, the class
of $Z$ in  $A^*(G{\Cal F})$
(i.e. the image of the class of the diagonal of \ $G\times_{_X}G$ \ via \
$(1\times\alpha)^*$) equals
$$
\sum\limits_{strict \ I\subset\rho_n} \widetilde Q_I(D_{G\Cal F}\hak) \cdot
\widetilde Q_{\rho_n\smallsetminus I}(R_{G{\Cal F}}\hak).
$$
\endproclaim
\smallskip

Thus the problem of computing the classes of the $\Omega(a.)$'s
is essentially that
of calculation  $p_*(\widetilde Q_ID\hak)$ where $p:{\Cal F}\to X$ is the
projection map; then
we use a base change.
\medskip

Consider now the case of the orthogonal Grassmannian parametrizing rank $n$
subbundles of $V$,
where $rank\ V=2n+1$. The results of Lemma 2.3, Proposition 2.5 and
Corollary 2.6 translate mutatis mutandis to this case with $\widetilde
Q$-polynomials
replaced by $\widetilde P$-polynomials (using essentially Theorem 5.23(ii) ).
Thus the problem of computing the classes of the $\Omega(a.)$'s
is essentially that
of calculation $p_*(\widetilde P_ID\hak)$ where $p:{\Cal F}\to X$ is the
projection.
\medskip

Finally, consider the even orthogonal case.
Supppose that $V$ is a vector bundle of $rank\ V=2n$
endowed with a nondegenerate
orthogonal form. Let $G=OG_n'V$ or $G=OG_n''V$ following the notation of
Section 1.
The even orthogonal analog of Proposition 2.5 and Corollary 2.6
is obtained using Theorem 5.23(iii) and reads as follows:
\smallskip
\proclaim{\bf Proposition 2.7}
The class of the diagonal of the Grassmannian bundle in
$A^*(G\times_{_X}G)$
equals
$$
[\Delta] = \sum \widetilde P_I(R_1\hak)_{GG} \cdot
\widetilde P_{\rho_{n-1}\smallsetminus I}(R_2\hak)_{GG} ,
$$
the sum over all strict $I\subset \rho_{n-1}$, $GG=G\times_{_X}G$
and $R_i$, $i=1,2$, are the tautological (sub)bundles on the corresponding
factors.
With $G\Cal F:=G\times_{_X} \Cal F'$
(resp. $G\Cal F:=G\times_{_X} \Cal F''$),
the class of $Z$ in  $A^*(G{\Cal F})$
(i.e. the image of the class of the diagonal of \ $G\times_{_X}G$ \ via \
$(1\times\alpha)^*$) equals
$$
\sum\limits_{strict \ I\subset\rho_{n-1}} \widetilde P_I(D_{G\Cal F}\hak) \cdot
\widetilde P_{\rho_{n-1}\smallsetminus I}(R_{G{\Cal F}}\hak).
$$
\endproclaim
\smallskip

Thus the problem of computing the classes of the $\Omega(a.)$'s
is essentially that of
calculation $p_*'(\widetilde P_ID\hak)$ and  $p_*''(\widetilde P_ID\hak)$
where $p':{\Cal F}'\to X$ and
$p'':{\Cal F}'' \to X$ are the projection maps.

\bigskip\medskip

\centerline{\bf 3. Subbundles intersecting an $n$-subbundle in dim
$\geqslant k$}
\bigskip

We will now show an explicit computation in case $a_.=(n-k+1,\ldots,n)$~.
This computation relies on a simple linear algebra argument.
Another proof of Proposition 3.1,
using the algebra of divided differences, will be given
in Section 8.

We start with the Lagrangian case and follow the notation from Section 1.

\proclaim{\bf Proposition 3.1}
Assume \ $a_.=(n-k+1,\ldots,n)$~.
Let \ $I\subset \rho_n$ \ be a strict partition. If
\break
$(n,n-1,\ldots,k+1)\not\subset I$, then \ $p_*\widetilde Q_ID\hak=0$~.
In the opposite case, write \ $I=(n,n-1,\ldots,
\break
k+1,j_1,\ldots,j_l)$, where $j_l>0$ and $l\leqslant k$~.
Then \ $p_*\widetilde Q_ID\hak= \widetilde Q_{j_1,\ldots,j_l}V_n\hak$~.
\endproclaim

\demo{Proof}
It suffices to prove the formula for a vector bundle $V\to B$ endowed
with a symplectic form, $X$ equal to ${L}G_nV$
and $V_n$ equal to the tautological subbundle
on $LG_nV$. (Recall that $\Omega(n-k+1,\ldots,n;V.)$ depends only on $V_n$
; more precisely, it parametrizes Lagrangian rank $n$ subbundles $L$ of
$V$ such that
\ $rank(L\cap V_n) \geqslant k$~.)
The variety \ff \ in this case parametrizes triples $(L,M,N)$ of vector
bundles over $B$ such that $L$ and $N$ are Lagrangian rank $n$
subbundles of $V$
and $M$ is a rank $k$ subbundle of $L\cap N$.
Let $W.$:
$W_1\subset W_2\subset\ldots\subset W_n$ be a flag of Lagrangian subbundles
of $V$ with
$rank \ W_i=i$.
For a partition $J=(j_1>\ldots>j_l>0) \subset \rho_k$ \ ,
$$
\alpha_J = \Omega(n+1-j_1,\ldots,n+1-j_l;W.) =
\bigl\{L\in X | \ rank (L\cap W_{n+1-j_h})\ge h, \ h=1,\dots,l \bigr\}
$$
defines a Schubert cycle whose class has the dominant term (w.r.t. $V_n$)
equal to $\widetilde Q_JV_n\hak\in A^*(X)$. It is well known
that $\alpha_J$ is an irreducible subvariety of $X$ provided $B$
is irreducible.

\medskip

Similarly, for a partition $I=(i_1>\ldots>i_l>0) \subset \rho_n$, \
$q:\Cal F\to LG_nV$ the projection on the third factor,
\medskip

\noindent
$A_I = q^*\Omega(n+1-i_1,\ldots,n+1-i_l;W.)=$
\smallskip

\noindent
$ \ \ \ \ \ \ \ \ \ \ \ \ \ \ \ \ \ \ \ \ \ \ \ \ \ \ \ \ \ \ \ \ \ \ \ \ \ \ \
\ \ \ \ \ \ \
= \bigl\{(L,M,N)\in{\Cal F} | \ rank (N\cap W_{n+1-i_h})\ge h, \ h=1,\dots,l
\bigr\}$
\smallskip

\noindent
defines a cycle whose class has the dominant term (w.r.t. $D$)
equal to $\widetilde Q_ID\hak\in
A^*({\Cal F})$. Also, $A_I$ is an irreducible subvariety of $\Cal F$ provided
$B$ is
irreducible.
\bigskip

We will show (the push-forward is taken on the cycles level) that:
{\parindent=25pt
\smallskip
\item{\bf 1)} If $I\not\supset(n,n-1,\ldots,k+1)$ then $p_*A_I=0$.
          Passing to the rational equivalence classes, this implies
	  $p_*\widetilde Q_ID\hak= 0$.
\medskip
\item{\bf 2)} If $I\supset(n,n-1,\ldots,k+1)$ i.e. $I=(n,n-1,\ldots,
          k+1,j_1,\ldots,j_l)$, where $j_l>0$ and $l\le k$, then
	  $p_*A_I=\alpha_J$ where $J=(j_1,...,j_l)$.
	  Then, passing to the rational equivalence classes
	  (and using the projection formula), we get the following
	  equality involving the dominant terms: $p_*\widetilde Q_ID\hak
	  = \widetilde Q_JV_n\hak$.
\par}

\bigskip

Observe that {\bf {1)}} holds if $l(I)\le n-k$  because we then have
$codim_{_\Cal F}A_I
= |I| < n+(n-1)+\ldots+(k+1)$, which is the dimension
of the fiber of $p$.
We will need the following:
\smallskip
\noindent
\underbar{Claim} \ Let $I \subset \rho_n$ be a strict
partition. Let $l = card \{h \ | \ i_{n-k+h} \ne 0\}$.
Assume that $l>0$.
Then one has
$$p(A_I) \subset \alpha_{i_{n-k+1},i_{n-k+2},\ldots,i_{n-k+l}}.
\tag *
$$
\medskip

\noindent
Indeed, for $(L,-,N)\in A_I$, since
$rank(L\cap N)\geqslant k$,
the inequality $rank(N\cap W_r)\geqslant h$
implies $rank(L\cap W_r)\geqslant h-(n-k)$
for every $h,r$; this gives (*).
\bigskip

\noindent
{\bf {1)}} To prove this assertion we first use (*) (by the above
remark we can assume that
$l(I)>n-k$) and thus get
$$
codim_{_\Cal F}A_I-codim_{_X}p(A_I)\le (i_1+\ldots+i_{n-k+l})
-(i_{n-k+1}+\ldots+i_{n-k+l})
=i_1+\ldots+i_{n-k}.
$$
\medskip

\noindent
Then, since $I\not\supset(n,n-1,\ldots,k+1)$, we have
$$
i_1+\ldots+i_{n-k}<n+\ldots+(k+1),$$
\smallskip

\noindent
where the last number is the dimension of the fiber of $p$. Hence comparison
of the latter inequality with the former yields
$p_*A_I=0$.
\bigskip

\noindent
{\bf {2)}} To prove this, it suffices to show
$p(A_I)\subset\alpha_J$,
$dim\ A_I= dim\ \alpha_J$; and if $p_*A_I=d\cdot \alpha_J$ for some
$d\in {\Bbb Z}$ then $d=1$.
We have:

\medskip

$p(A_I)\subset\alpha_J$ : this is a direct consequence of (*).
\medskip

$dim\ A_I = dim\ \alpha_J$ : this results from comparison
of the following three formulas
$dim \ \Cal F = dim \ X + k(n-k) + (n-k)(n-k+1)/2$, \
$codim_{_X}\alpha_J = |J|$, \
and $ codim_{_\Cal F}A_I = n+\ldots+(k+1)+|J|$.

\bigskip

Therefore  $p_*A_I=d\cdot\alpha_J$  for some integer $d$.
To show $d=1$ it suffices to find an open subset
$U\subset \alpha_J$ such that $p|_{p^{-1}U}\ :\ p^{-1}U\to U$
is an isomorphism.
We define the open subset $U$ in question as \ $\alpha_J\smallsetminus
\Omega(n-k;W.)$. More explicitly,
$U$ is defined by the conditions:
$$
rank(L\cap W_{n+1-j_1})\ge 1,\ \ldots, \ rank(L\cap W_{n+1-j_l})\ge l \ \
\hbox{and} \ \
L\cap W_{n-k}=(0).
$$
Observe that these conditions really define an open nonempty subset of
$\alpha_J$
because
$\Omega(n+1-j_1,\ldots,n+1-j_l;W.)\not\subset
\Omega(n+1-(k+1);W.)$ for $J\subset \rho_k$. (Recall that for \
$I=(i_1>\ldots >i_l>0), \ J=(j_1>\ldots >j_{l'}>0)$  one has
$\Omega(n+1-i_1,\ldots,n+1-i_l;W.) \subset
\Omega(n+1-j_1,\ldots,n+1-j_{l'};W.)$
iff $I\supset J$. )
\medskip

Since our problem of showing that $d=1$ is of local nature, we can assume that
$B$ is
a point and deal with vector spaces instead of vector bundles.
Let us choose a basis $e_1,\ldots,e_n,f_1,\ldots,f_n$ such
that, denoting the symplectic form by $\Phi$, we have $\Phi(e_i,e_j)=
0=\Phi(f_i,f_j)$ and $\Phi(e_i,f_j)=-\Phi(f_j,e_i)=\delta_{i,j}$.
Assume that
$W_i$ is generated by the first $i$ vectors of $\{e_j\}$. Let $W^i$ be
the subspace generated by the last $i$ vectors of $\{e_j\}$. Moreover, let
$\widetilde W_{i}$ be the subspace generated by the first $i$ vectors of
$\{f_j\}$
and $\widetilde W^{i}$ be the subspace generated by the last $i$ vectors of
$\{f_j\}$.
\medskip

Observe that for a strict partition $\rho_n\supset I\supset(n,n-1,\ldots,k+1)$
a necessary
condition for "$(-,-,N)\in A_I$" \ is \ "$ N\supset W_{n-k}$".
(This corresponds
to the first $(n-k)$ Schubert conditions defining $A_I$.) \ On the other hand,
if $L\in U$ then $L\cap W_{n-k}=(0)$ and consequently
$L$ must contain
$\widetilde W_{n-k}$ (from the rest, i.e.
$W^k \oplus \widetilde W^k$,  we can get at most
$k-$dimensional isotropic subspace). Hence also
$|L\cap\bigl(W^k \oplus \widetilde W^k \bigr)|=k$ \ ($| - |$ denotes the
dimension). We conclude that a necessary choice for an $n$-dimensional
Lagrangian subspace $N$ such that $(L,M,N)\in A_I$ for some $M$, is
$$
N := W_{n-k}\oplus L\cap\bigl( W^k \oplus \widetilde W^k \bigr).
$$
It follows from the above discussion that $N$ is really a Lagrangian
subspace of dimension $n$ and it satisfies the first $(n-k)$ Schubert
conditions defining
$A_I$. $N$ also satisfies the
last $l \ (\le k)$ Schubert conditions defining $A_I$ :
since $|L\cap W_{n+1-j_h}|\ge h$
and
$L\cap W_{n-k}=(0)$, we have $|N\cap W_{n+1-j_h}|= |W_{n-k}|+h \ge n-k+h$ \ for
 \ $h=1,\ldots,l$.
\smallskip

Moreover, since $|L\cap N|=k$, the subspace $M$ above is determined uniquely:
$M = L\cap N$.

Summing up, we have shown that $d=1$; this ends the proof of
{\bf {2)}}.
\medskip

Thus the proposition has been proved.
\qed
\enddemo

\proclaim { Proposition 3.2} One has in $A^*(G)$,
$$
[\Omega(n-k+1,\ldots,n-1,n)] = \sum\limits_{strict\ I\subset\rho_k}
\widetilde Q_I(V_n\hak)_G\cdot \widetilde Q_{\rho_k\smallsetminus I}(R\hak)\ .
$$
\endproclaim

\demo {Proof}
This formula is obtained directly by pushing forward via $(p_1)_*$
the class of $Z$ in\ $A^*(G\Cal F)$ given by
$$
\sum\limits_{strict\ I\subset\rho_n}
\widetilde Q_I(D_{G{\Cal F}}\hak)\cdot \widetilde Q_{\rho_n\smallsetminus
I}(R_{G{\Cal F}}\hak)
$$
(see Corollary 2.6), with the help of Proposition 3.1.
\qed
\enddemo

\example{\bf Example 3.3} \ For successive $k$ (and any $n$) the formula reads
(with  $D=
D_{G{\Cal F}}, \ R=R_{G{\Cal F}}$ for brevity):
{\parindent=30pt
\smallskip
\item{k=1\ \ } $\widetilde Q_1D\hak+\widetilde Q_1R\hak$;
\smallskip
\item{k=2\ \ } $\widetilde Q_{21}D\hak+\widetilde Q_2D\hak\cdot \widetilde
Q_1R\hak+
               \widetilde Q_1D\hak\cdot \widetilde Q_2R\hak+\widetilde
Q_{21}R\hak$;
\smallskip
\item{k=3\ \ } $\widetilde Q_{321}D\hak+\widetilde Q_{32}D\hak\cdot \widetilde
Q_1R\hak+
               \widetilde Q_{31}D\hak\cdot \widetilde Q_2R\hak+\widetilde
Q_{21}D\hak\cdot \widetilde Q_3R\hak+
	       \widetilde Q_3D\hak\cdot \widetilde Q_{21}R\hak+   \hfill\break
	       \widetilde Q_2D\hak\cdot \widetilde Q_{31}R\hak+\widetilde
Q_1D\hak\cdot \widetilde Q_{32}R\hak+
	       \widetilde Q_{321}R\hak$.
\par}

\endexample

\medskip

In the odd orthogonal case, the analogs of Propositions 3.1 and 3.2 are
obtained by replacing $\widetilde Q$-polynomials
by $\widetilde P$-polynomials.

\proclaim { Proposition 3.4}
(i) Assume \ $a_.=(n-k+1,\ldots,n)$.
Let  $I\subset \rho_n$  be a strict partition. If
$(n,n-1,\ldots,k+1)\not\subset I$, then \ $p_*\widetilde P_ID\hak=0$~.
In the opposite case, write \ $I=(n,n-1,\ldots,
\break
k+1,j_1,\ldots,j_l)$, where $j_l>0$ and $l\leqslant k$~.
Then \ $p_*\widetilde P_ID\hak= \widetilde P_{j_1,\ldots,j_l}V_n\hak$~.
\smallskip

\noindent
(ii) One has in $A^*(OG_nV)$,
$$
[\Omega(n-k+1,\ldots,n-1,n)] = \sum\limits_{strict\ I\subset\rho_k}
\widetilde P_I(V_n\hak)_G\cdot \widetilde P_{\rho_k\smallsetminus I}(R\hak)\ .
$$
\endproclaim

Assertion (ii) follows from (i) like Proposition 3.2 follows
from Proposition 3.1.
The proof of (i) is essentially the same as the one
of Proposition 3.1.
More precisely, in the proof of (i), $\alpha_J$ and $A_I$
are defined in the same way as in the proof of this proposition.
Also, the whole reasoning is the same, word by word, except of
the following one point. To prove that $d=1$ one chooses now
a basis $e_1,\ldots,e_n,f_1,\ldots,f_n,g$ such that
denoting the orthogonal form by $\Phi$, we have $\Phi(e_i,e_j)=
\Phi(f_i,f_j)=\Phi(e_i,g)=\Phi(f_j,g)=0$, $\Phi(e_i,f_j)=\Phi(f_j,e_i)
=\delta_{i,j}$ and $\Phi(g,g)=1$.
Then $W_i$, $W^i$,
$\widetilde W_i$ and $\widetilde W^i$ defined like
in the proof of Proposition 3.1 allow us to
show that $d=1$ exactly in the same way as in the proof of this
proposition.

\bigskip

Let us pass now to the even orthogonal case.
So let $V\to X$ ($X$ is connected) be a rank $2n$ vector bundle endowed with
a nondegenerate quadratic form.
Fix an isotropic rank $n$ subbundle $V_n$ of $V$.
Recall that for $k\equiv n \ (mod \ 2)$ by $p':{\Cal F}'\to X$
we denote the flag bundles parametrizing flags $A_1\subset A_2$
of subbundles of $V$ such that $rank\ A_1=k$, $rank\ A_2=n$, $A_1\subset V_n$
and $A_2$ is isotropic with $dim(A_2\cap V_n)_x\equiv n\ (mod \ 2)$
for every $x\in X$. Similarly, for $k\equiv n+1 \ (mod \ 2)$
by $p'':{\Cal F}''\to X$
we denote the flag bundle parametrizing flags $A_1\subset A_2$
of subbundles of $V$ such that $rank\ A_1=k$, $rank\ A_2=n$, $A_1\subset V_n$
and $A_2$ is isotropic with $dim(A_2\cap V_n)_x\equiv n+1\ (mod \ 2)$
for every $x\in X$.
\medskip

\medskip

In the even orthogonal case the analog of Proposition 3.1 reads:

\proclaim{\bf Proposition 3.5}
Let $I\subset\rho_{n-1}$ be a strict partition.
If $(n-1,n-2,\ldots,k)\not\subset I$ then $p_*'\widetilde P_ID\hak=0$.
In the opposite case, write $I=(n-1,n-2,\ldots,k,j_1,\ldots,j_l)$, where
$j_l>0$ and $l\leqslant k-1$.
Then
$$
p_*'\widetilde P_ID\hak=\widetilde P_{j_1,\ldots,j_l}V_n\hak.
$$
The same formula is valid for $p_*''$.
\endproclaim

\demo{Proof}
We consider first the case of $p_*'$ i.e. $k\equiv n \ (mod \ 2)$.
It suffices to prove the formula for a rank $2n$ vector bundle $V\to B$
(we assume that $B$ is irreducible)
endowed with a nondegenerate orthogonal form, $X$ equal to ${O}G_n'V$
or ${O}G_n''V$ and $V_n$ equal to the tautological subbundle on $X$.
Then the variety \ff$'$ parametrizes triples $(L,M,N)$ such that
$dim(L\cap N)_b\equiv n\ (mod\ 2)$ for every $b\in B$
(i.e. $L$ and $N$ either belong together
to ${O}G_n'V$ or together to
${O}G_n''V$) and $M$ is a rank $k$
subbundle of $L\cap N$.

We will now prove the proposition for $X=OG_n'V$.
(Obvious modifications lead to a proof in the case $X=OG_n''V$.)
Since the strategy of proof is the same as in the Lagrangian case,
we will skip those parts of the reasoning which have appeared already in the
proof
of Proposition 3.1. Let $W.: W_1\subset W_2\subset \ldots \subset W_n$
be an isotropic
flag in $V$.

\smallskip

For $J=(j_1>\ldots>j_l>0)\subset\rho_{k-1}$ we define
\smallskip

\centerline {$\alpha_J=\Omega(n-j_1,\ldots,n-j_l;W.) \ \ if \ \ l\equiv n \
(mod\ 2)  \
and$}
\smallskip

\centerline {$\alpha_J=\Omega(n-j_1,\ldots,n-j_l,n;W.) \ \ if \ \
l\equiv n+1 \ (mod\ 2).$}
\smallskip

Similarly for $I=(i_1>\ldots>i_l>0)\subset \rho_{n-1},\ q:\Cal F'\to OG_n'V$
the projection on the third factor, we define
\smallskip

\centerline {$A_I=q^*\Omega(n-i_1,\ldots,n-i_l;W.) \ \
if \ \ l\equiv n \ (mod\ 2)  \ and$}
\smallskip

\centerline {$A_I=q^*\Omega(n-i_1,\ldots,n-i_l,n;W.) \ \
if \ \ l\equiv n+1 \ (mod\ 2).$}
\smallskip

It is known that
$\alpha_J$ and $A_I$ are irreducible subvarieties provided $B$ is.
The dominant terms
of the classes of $\alpha_J$ and $A_I$ are equal to
$\widetilde P_JV_n\hak$ and $\widetilde P_ID\hak$ respectively.
\smallskip

The proposition now follows from:

{\parindent=25pt
\item{1)} If $I\not\supset(n-1,n-2,\ldots,k)$ \ then \ $p'_*A_I=0$.
\item{2)} If $I\supset(n-1,n-2,\ldots,k)$ \ i.e $I=(n-1,n-2,\ldots,
k+1,k,j_1,\ldots,j_l)$,
where $j_l>0$ and $l\le k-1$, then \ $p'_*A_I=\alpha_J$  \ where
$J=(j_1,\ldots,j_l)$.
\par}
\smallskip

Assertion 1) (being obvious if $l(I)<n-k$) \ is a consequence of:
\smallskip
\underbar{Claim:} For every strict partition $I\subset\rho_{n-1}$, let
$l=card\{h| \ i_{n-k+h}\ne 0\}$. Assume that $l>0$.
Then one has
$$
p'(A_I)\subset\alpha_{i_{n-k+1},\ldots,i_{n-k+l}}.
\tag *
$$
\smallskip
Inclusion (*) also implies $p'(A_I)\subset\alpha_J$ in 2).
The equality $dim\ p'(A_I)=dim\ \alpha_J$ now follows from:
$dim\ \Cal F'=dim\ X+k(n-k)+(n-k)(n-k-1)/2$, $codim_{_X}\alpha_J=|J|$ and
$codim_{_{\Cal F'}}A_I=(n-1)+\ldots+k+|J|$.
\smallskip
Therefore $p_*'A_I=d\cdot \alpha_J$ for some integer $d$. To prove that
$d=1$ it is sufficient to show
an open subset $U\subset\alpha_J$ such that $p'|_{(p')^{-1}U}$:
$(p')^{-1}U\to U$ is an isomorphism.
The open subset $U$ in question parametrizes those $L\in \alpha_J$
for which $L\cap W_{n-k} = (0)$.

\smallskip
The problem being local, we can assume that $B$ is a point.
Let $e_1,\ldots,e_n,f_1,\ldots,f_n$ be a basis of $V$ such that denoting the
form by $\Phi$ we have $\Phi(e_i,e_j)=\Phi(f_i,f_j)=0$,
$\Phi(e_i,f_j)=\Phi(f_j,e_i)=\delta_{i,j}$ and $W_i$ is spanned by
$e_1,\ldots,e_i$.
Define $W^i,\widetilde W_i$ and $\widetilde W^i$
as in the proof of Proposition 3.1.

Now, given $L\in U$, the unique $N$ such that $(L,M,N)\in A_I$ for some $M$,
is defined also as in the proof of Proposition 3.1:
$N := W_{n-k}\oplus L\cap\bigl( W^k \oplus \widetilde W^k \bigr)$.

This $N$ is isotropic and satisfies the first $n-k$ Schubert conditions
because it contains $W_{n-k}$. Moreover, it satisfies the last
$l(\le k-1)$ Schubert conditions defining $A_I$:
since $|L\cap W_{n-j_h}|\ge h$ and $L\cap W_{n-k}=(0)$, we have
$|N\cap W_{n-j_h}|=|W_{n-k}|+h\ge n-k+h$ for $h=1,\ldots,l$.
Finally, the $M$ above is determined uniquely:
$M=L\cap N$, and $p'|_{(p')^{-1}U}$ is an isomorphism.

\medskip

We next consider the case of $p_*''$, i.e. $k\equiv n+1 \ (mod \ 2)$.
It suffices to prove the formula for a rank $2n$ vector bundle $V\to B$
($B$ is irreducible)
endowed with a nondegenerate orthogonal form, $X$ equal to ${O}G_n'V$
or ${O}G_n''V$ and $V_n$ equal to the tautological subbundle on $X$.
Then the variety \ff$''$ parametrizes triples $(L,M,N)$ such that
$dim(L\cap N)_b\equiv n+1\ (mod\ 2)$ for every $b\in B$
(i.e. $L$ and $N$  belong to different components
${O}G_n'V$ and ${O}G_n''V$) and $M$ is a rank $k$
subbundle of $L\cap N$.

We will prove the proposition for $X=OG_n''V$.
(Obvious modifications lead to a proof in the case $X=OG_n'V$.)
Let $W.: W_1\subset W_2\subset \ldots \subset W_n$ be an isotropic flag in $V$.
For $J=(j_1>\ldots>j_l>0)\subset\rho_{k-1}$ we define
\smallskip

\centerline {$\alpha_J=\Omega(n-j_1,\ldots,n-j_l;W.) \ \ if \ \ l\equiv n+1
 \ (mod\ 2)  \
and$}
\smallskip

\centerline {$\alpha_J=\Omega(n-j_1,\ldots,n-j_l,n;W.) \ \ if \ \
l\equiv n \ (mod\ 2).$}
\smallskip

Similarly for $I=(i_1>\ldots>i_l>0)\subset \rho_{n-1},\ q:\Cal F''\to OG_n'V$
the projection on the third factor, we define
\smallskip
\centerline {$A_I=q^*\Omega(n-i_1,\ldots,n-i_l,n;W.) \ \
if \ \ l\equiv n+1 \ (mod\ 2) \ and$}
\smallskip
\centerline {$A_I=q^*\Omega(n-i_1,\ldots,n-i_l;W.) \ \
if \ \ l\equiv n \ (mod\ 2).$}
\smallskip

The dominant terms
of the classes of $\alpha_J$ and $A_I$ are equal to
$\widetilde P_JV_n\hak$ and $\widetilde P_ID\hak$ respectively.
The proposition now follows from:

{\parindent=25pt
\item{1)} If $I\not\supset(n-1,n-2,\ldots,k)$ \ then \ $p_*''A_I=0$.
\item{2)} If $I\supset(n-1,n-2,\ldots,k)$ \ i.e $I=(n-1,n-2,\ldots,
k+1,k,j_1,\ldots,j_l)$,
where $j_l>0$ and $l\le k-1$, then \ $p_*''A_I=\alpha_J$  \ where
$J=(j_1,\ldots,j_l)$.
\par}
\smallskip
The proof of these assertions is analogous to the one above.
For every strict partition $I\subset\rho_{n-1}$ and
$l=card\{h| \ i_{n-k+h}\ne 0\}>0$, one has
$p''(A_I)\subset\alpha_{i_{n-k+1},\ldots,i_{n-k+l}}$,
which implies 1) and
$p''(A_I)\subset\alpha_J$ in 2); moreover, for the dimensions reasons we have
$p_* ''A_I=d\cdot \alpha_J$ for some integer $d$. One finishes the proof
like in the case of $p_*'$,
by showing that $p''|_{(p'')^{-1}U}$:  $(p'')^{-1}U\to U$ is an isomorphism,
where
an open subset $U\subset\alpha_J$
parametrizes those $L\in \alpha_J$ for which $L\cap W_{n-k}=(0)$.
Hence $d=1$ and the proof is complete.
\qed
\enddemo
\smallskip

\proclaim{\bf Proposition 3.6}
If $k\equiv n\ (mod\ 2)$ (resp. $k\equiv n+1\ (mod\ 2)$ ) then one has in
$A^*({O}G_n'V)$ (resp. in $A^*({O}G_n''V)$ ),
$$
[\Omega(n-k+1,\ldots,n-1,n)] =
\sum\limits_{strict\ I\subset\rho_{k-1}}\widetilde P_I(V_n\hak)_G\cdot
\widetilde P_{\rho_{k-1}\smallsetminus I}(R\hak).
$$
\endproclaim
\demo{Proof}
This formula is obtained directly by pushing forward via $p_*'$
(resp. $p_*''$) the class of $Z$ in $A^*(G\Cal F)$ given by
$$
\sum\limits_{strict\ I\subset\rho_{n-1}}\widetilde P_I
\bigl(D\hak_{G\Cal F}\bigr)
\cdot \widetilde P_{\rho_{n-1}\smallsetminus I}
\bigl(R\hak_{G\Cal F}\bigr),
$$
where $G\Cal F = G\times_{_X}{\Cal F}'$ (resp.
$G\Cal F = G\times_{_X}{\Cal F}''$)
, using Propositions 2.7
and 3.5.
\qed
\enddemo

\bigskip\medskip

\centerline{\bf 4. $\widetilde Q$-polynomials and their properties}

\bigskip
In this section we define a family of symmetric polynomials modelled
on Schur's $Q$-polynomials. In Schur's Pfaffian-definition (see [S]),
we replace $Q_i$ by $e_i$ -- the $i$-th elementary symmetric polynomial.
After this modification one gets an interesting family of
symmetric polynomials  $\widetilde Q_{I}$  (indexed by all partitions)
whose properties are studied in this section
and then applied in the next ones.
It turns out that  $\widetilde Q_{I}$ is the Young dual (in sense of the
involution $\omega$ of [Mcd1, I.2.(2.7)] to the Hall--Littlewood polynomial
 $\widetilde Q_I(Y;q)$ where the alphabet $Y$ is equal to $X_n / (1-q)$ in
the sense of $\lambda$-rings, specialized with $q=-1$
([L-L-T], [D-L-T]).
Though most of the
properties of the  $\widetilde Q_I$ given in this section can be deduced
from the theory of Hall--Littlewood polynomials, we give here their proofs
using the Pfaffian definition. The only exception is made for the Pieri-type
formula which is deduced from the one for Hall--Littlewood polynomials.

\medskip
Let $X=(x_1,x_2,\ldots)$ be a sequence of independent variables.
Denote by $X_n$ the subsequence $(x_1,\ldots,x_n)$.
We set $\widetilde Q_i(X_n):=e_i(X_n)$ -- the $i$-th elementary symmetric
polynomial in $X_n$.
Given two nonnegative integers $i,j$ we define
$$
\widetilde Q_{i,j}(X_n)=\widetilde Q_i(X_n)\widetilde Q_j(X_n) +
2 \sum\limits^j_{p=1}
(-1)^p\widetilde Q_{i+p}(X_n)\widetilde Q_{j-p}(X_n)
$$
Finally, for any (i.e. not necessary strict) partition
$I=(i_1\geqslant i_2\geqslant\ldots\geqslant i_k\ge 0)$,
with even $k$ (by putting $i_k=0$ if necessary), we set
$$
\widetilde Q_I(X_n) = Pf\Bigl(\widetilde Q_{i_p,i_q}(X_n)
\Bigr)_{1\leqslant p<q\leqslant k}.
$$
Equivalently (in full analogy to [S, pp.224--225]),
$\widetilde Q_I(X_n)$ is defined recurrently on $l(I)$,
by putting for odd $l(I)$,
$$
\widetilde Q_I(X_n) = \sum\limits_{j=1}^{l(I)} (-1)^{j-1}
\widetilde Q_{i_j}(X_n)\widetilde Q_{I\smallsetminus \{i_j\}}(X_n),
\leqno (*)
$$
and for even $l(I)$,
$$
\widetilde Q_I(X_n) = \sum\limits_{j=2}^{l(I)} (-1)^j \widetilde
Q_{i_1,i_j}(X_n)
\widetilde Q_{I\smallsetminus\{i_1,i_j\}}(X_n).
\leqno (**)
$$
The latter case, with $l=l(I)$, can be rewritten as
$$
\widetilde Q_I(X_n)=
\sum\limits_{j=1}^{l-1}(-1)^{j-1}\widetilde Q_{i_j,i_l}(X_n)
\widetilde Q_{I\backslash\{i_j,i_l\}}(X_n).
\leqno(***)
$$
Note that assuming formally $i_l=0$, the relation (***) specializes to (*).
We will refer to the above equations as to Laplace-type developements or
simply recurrent formulas.
(Invoking the raising operators $R_{ij}$
([Mcd1,I], [D-L-T]) the above definition is rewritten
$$\widetilde Q_I(X_n)=\prod_{i<j}{1-R_{ij}
\over 1+R_{ij}} e_I(X_n),$$
where $e_I(X_n)$ is the product of the elementary symmetric polynomials
in $X_n$ associated with the parts of $I$.)
\medskip
We start with a useful linearity-type formula for $\widetilde Q$-polynomials.
\smallskip
\proclaim{\bf Proposition 4.1}
For any strict partition $I$ one has
$$
\widetilde Q_I(X_n) = \sum\limits_{j=0}^{l(I)} x_n^j
\Bigl( \sum\limits_{|I|-|J|=j} \widetilde Q_J(X_{n-1})\Bigr),
$$
where the sum is over all (i.e. not necessary strict)
partitions $J\subset I$ such that $I/J$
has at most one box in every row.
\endproclaim

\demo{Proof}
We use induction on $l(I)$.
\smallskip
{\parindent =20pt
\item{$1^\circ$  } $l(I)=1$. \ Since we have:
$e_i(X_n) = e_i(X_{n-1}) + x_ne_{i-1}(X_{n-1})$,
the assertion follows.
\smallskip
\item{$2^\circ$  } $l(I)=2$. \
We have for $i>j>0$ and with $e_i=e_i(X_n)$, $\bar e_i=e_i(X_{n-1}),
\bar e_{-1}=0$,
\smallskip
$\widetilde Q_{i,j}(X_n)=e_ie_j+2\sum\limits^j_{p=1}(-1)^pe_{i+p}e_{j-p} =$
\smallskip
$=(\bar e_i+x_n\bar e_{i-1})(\bar e_j+x_n\bar e_{j-1})+2\sum\limits^j_{p=1}
(-1)^p(\bar e_{i+p}+x_n\bar e_{i+p-1})(\bar e_{j-p}+x_n\bar e_{j-p-1})$
\smallskip
$=\Bigl(\bar e_i\bar e_j + 2\sum\limits_{p=1}^j (-1)^p\bar e_{i+p}
\bar e_{j-p}\Bigr)
+x_n\Bigl[\Bigl(\bar e_{i-1}\bar e_j+2\sum\limits_{p=1}^j (-1)^p\bar e_{i-1+p}
\bar e_{j-p}\Bigr)+$
\smallskip
\ \ \ \
$+\Bigl(\bar e_i\bar e_{j-1}+2\sum\limits_{p=1}^{j-1} (-1)^p\bar e_{i+p}
\bar e_{j-1-p}\Bigr)\Bigr]
+x^2_n\Bigl(\bar e_{i-1}\bar e_{j-1}+2\sum\limits_{p=1}^{j-1} (-1)^p
\bar e_{i-1+p}
\bar e_{j-1-p}\Bigr)$
\smallskip
$= \widetilde Q_{i,j}(X_{n-1})+x_n\Bigl[\widetilde Q_{i-1,j}(X_{n-1})+
\widetilde Q_{i,j-1}(X_{n-1})\Bigr]
+x^2_n \widetilde Q_{i-1,j-1}(X_{n-1}).$
\medskip
\item{$3^\circ$  }
By the remarks before the proposition, to prove the assertion in
general it suffices to show it by using the recurrent relation (***).
(Note that the R.H.S. of the formula of the proposition specializes after
the formal replacement $i_l:=0\ (l=l(I))$ to the expression asserted for
$(i_1>i_2>\ldots>i_{l-1})$.
\item{}\ \ \ So, let us assume that $l$ is even and set
$\widetilde Q_I:=\widetilde Q_I(X_n),\ \bar Q_I=\widetilde Q_I(X_{n-1})$.
Moreover, let $\Cal P(I,j)$ be the set of all partitions $J\subset I$ such
that $I\backslash J$ has at most one box in every row and $|I|-|J|=j$.
We have by induction on $l$,
$$
\widetilde Q_{I\backslash\{i_j,i_l\}}=\sum\limits_{r=0}^{l-2}x_n^r
\Biggl(\sum\limits_{J\in\Cal P(I\backslash\{i_j,i_l\},r)}\bar Q_J\Biggr).
$$
Therefore, using $2^\circ$ we have
$$
\eqalign{
\widetilde Q_I= &\sum\limits_{j=1}^{l-1}(-1)^{j-1}
\Bigl[\bar Q_{j_j,i_l}+x_n(\bar Q_{i_j-1,i_l}+\bar Q_{i_j,i_l-1})
+x^2_n\bar Q_{i_j-1,i_l-1}\Bigr] \cr
&\times \ \Biggl[\sum\limits_{r=0}^{l-2} x_n^r\Bigl(
\sum\limits_{J\in\Cal P(I\backslash\{i_j,i_l\},r)}\bar Q_J\Bigr)\Biggr] \cr}
$$
On the other hand, apply the relation (***) to the R.H.S. of the formula
in the proposition.
One gets:
$$
\sum\limits^l_{j=0}x^j_n\Bigl(\sum\limits_{J\in\Cal P(I,j)}\bar Q_J\Bigr)
=\sum\limits^l_{j=0}x_n^j\Biggl[\sum\limits_{J\in\Cal P(I,j)}
\Bigl(\sum\limits_{q=1}^{l-1}(-1)^{q-1}\bar Q_{j_q,j_l}\cdot
\bar Q_{J\backslash\{j_q,j_l\} } \Bigr)\Biggr]   .
$$
It is straightforward to verify that both these sums contain
$2^l(l-1)$ terms of the form
$$
(-1)^s x^j \bar Q_{a,b} \bar Q_{c_1,\ldots,c_{l-2}},
$$
and such a term appears in both sums if and only if

$$(c_1,\ldots,c_s,a,c_{s+1},\ldots,c_{l-2},b)\in \Cal P(I,j).$$
Thus the assertion follows and the proof of the proposition is complete.
\qed
\par}
\enddemo

\medskip
\proclaim{\bf Proposition 4.2}: \ $\widetilde Q_{i,i}(X_n) =
e_i(x_1^2,\ldots, x_n^2)$.
\endproclaim

\demo{Proof}
By definition we have $\bigl(e_i=e_i(X_n)\bigr)$:
$$
\widetilde Q_{i,i}(X_n) = e_ie_i-2e_{i-1}e_{i+1}+2e_{i+2}e_{i-2}-\ \ldots
=\sum\limits_{p=0}^{2i} (-1)^{p+i}e_pe_{2i-p}.
$$
On the other hand, with an indeterminate $t$, we have
$$
(1+x_1t)\ldots(1+x_nt)(1-x_1t)\ldots(1-x_nt) = (1-x_1^2t^2)\ldots(1-x_n^2t^2),
$$
or equivalently,
$$
\Bigl(\sum e_pt^p\Bigr)\Bigl(\sum(-1)^qe_qt^q\Bigr) =
\sum(-1)^ie_i(x_1^2,\ldots,x_n^2)t^2,
$$
which implies
$$
(-1)^ie_i(x_1^2,\ldots,x_n^2) = \sum\limits_{p=0}^{2i}(-1)^pe_p e_{2i-p}.
$$
Comparison of the latter equation with the first one gives the assertion.
\qed
\enddemo

\medskip
\proclaim{\bf Proposition 4.3}
For partitions $I'=(i_1,i_2,\ldots,j,j,\ldots,i_{k-1},i_k)$ and $I=(i_1,
\ldots,i_k)$,
the following equality holds
$$
\widetilde Q_{I'}(X_n)=\widetilde Q_{j,j}(X_n)\widetilde Q_I(X_n).
$$
\endproclaim

\demo{Proof}
Write $\widetilde Q_I$ for $\widetilde Q_I(X_n)$.
We use induction on $k$.
For $k=0$, the assertion is obvious.
For $k=1$, we have $\widetilde Q_{i,j,j}=\widetilde Q_i\widetilde Q_{j,j}$
and $\widetilde Q_{j,j,i}=\widetilde Q_{j,j}\widetilde Q_i$
by the Laplace type developements, so the assertion follows.

\smallskip
In general, it suffices to show the assertion inductively, using the
relation (***), if the marked "$j$" does not appear on the last place;
and independently, to prove it (inductively) for
$I'=(i_1,\ldots,i_k,j,j)$.
In both instances $k$ is assumed to be even.

In the former case, using (***) we get
$$
\eqalign{
\widetilde Q_{I'} = &\widetilde Q_{i_1,j_k} \widetilde
Q_{i_2,\ldots,j,j,\ldots,i_{k-1}}-\ldots
\pm\widetilde Q_{j,i_k}\widetilde Q_{i_1,i_2,\ldots,j,\ldots,i_k}\cr
&\mp\widetilde Q_{j,i_k}\widetilde Q_{i_1,i_2,\ldots,j,\ldots,i_k}\pm\ldots
-\widetilde Q_{i_{k-1},i_k}\widetilde Q_{i_1,\ldots,j,j,\ldots,i_{k-2}} \cr}
$$
and the assertion follows from the induction assumption by using the relation
(***) w.r.t. $\widetilde Q_{i_1,\ldots,i_k}$ once again.

In the latter case we use the relation (**). We have
$$
\eqalign{
\widetilde Q_{i_1,\ldots,i_k,j,j}= \ &\widetilde Q_{i_1,i_2}
\widetilde Q_{i_3,\ldots,i_k,j,j}-\ldots+
\widetilde Q_{i_1,i_k}\widetilde Q_{i_2,\ldots,i_{k-1},j,j}     \cr
&-\widetilde Q_{i_1,j}\widetilde Q_{i_2,\ldots,i_k,j}+
\widetilde Q_{i_1,j}\widetilde Q_{i_2,\ldots,i_k,j}   \cr}
$$
and the assertion follows from the induction assumption by using the relation
(**) w.r.t. $\widetilde Q_{i_1,\ldots,i_k}$ once again.
\qed
\enddemo

\medskip
\proclaim{\bf Lemma 4.4}
Let $I=(i_1,i_2,\ldots ,i_k)$ be a partition.
If $i_1>n$ then $\widetilde Q_I(X_n)=0$.
\endproclaim

\demo{Proof}
We use induction on $l(I)$.
For $l(I)=1,2$ the assertion is obvious because $e_p(x_1, \ldots ,x_n)=0$
for $p>n$.
For bigger $l(I)$ one uses induction on the length and the recurrent formulas,
which immediately imply the assertion.
\qed
\enddemo

\example{\bf Example 4.5}
The following equalities hold:
(in $1^\circ$ and $2^\circ$ we set $\widetilde Q_I:=\widetilde Q_I(X_n)$
for brevity)

\bigskip\noindent
$1^\circ \ \ \ \widetilde Q_{5544441}=\widetilde Q_{55}\widetilde Q_{44441}
=\widetilde Q_{55}\widetilde Q_{44}\widetilde Q_{441}
=\widetilde Q_{55}\widetilde Q_{44}\widetilde Q_{44}\widetilde Q_1
=\widetilde Q_{55}\widetilde Q_{4444}\widetilde Q_1$;

\bigskip\noindent
$2^\circ \ \ \ \widetilde Q_{5554443331}
=\widetilde Q_{55}\widetilde Q_{44}\widetilde Q_{33}\widetilde Q_{5431}
=\widetilde Q_{554433}\widetilde Q_{5431}$;

\bigskip\noindent
$3^\circ$ \ \ \ Here, we set $\bar Q_I := \widetilde Q_I(x_1,x_2)$,
$ {\bar Q_I}' := \widetilde Q_I(x_1)$. Then

\smallskip
$\widetilde Q_{321}(x_1,x_2,x_3)=$
\smallskip \ \ \ \ \
$=x_3\bar Q_{221}+x^2_3(\bar Q_{211}+\bar Q_{22}) + x^3_3\bar Q_{21} =
x_3\bar Q_{22}\bar Q_1+x^2_3(\bar Q_{11}\bar Q_2+\bar Q_{22})+x^3_3\bar Q_{21}$
\smallskip \ \ \ \ \
$=x_3e_2(x^2_1,x^2_2)(x_1+x_2)+x^2_3\bigl[e_1(x^2_1,x^2_2)x_1x_2+
e_2(x^2_1,x^2_2)\bigr]+x^3_3\bigl(x_2{\bar Q_{11}}'+
x^2_2{\bar Q_1}'\bigr)$
\smallskip \ \ \ \ \
$=x_3(x^2_1x^2_2)(x_1+x_2)+x^2_3\bigl[(x^2_1+x^2_2)x_1x_2+x^2_1x^2_2\bigr]
+x^3_3\bigl(x_2x^2_1+x^2_2x_1\bigr)$.
\endexample

\medskip
By iterating the linearity formula for $\widetilde Q_I(X_n)$
(Proposition 4.1), we get the following
algorithm for decomposition of $\widetilde Q_I=\widetilde Q_I(X_n)$
into a sum of monomials:

\medskip
\noindent
1. \ If $I$ is not strict, we factorize
$$
\widetilde Q_I=\widetilde Q_{k_1,k_1}\cdot \widetilde Q_{k_2,k_2}\cdot\ldots
\cdot \widetilde Q_{k_l,k_l}\cdot \widetilde Q_L,
$$
where $L$ is strict (we use Proposition 4.3).

\smallskip
\noindent
2. \ We apply the linearity formula to $\widetilde Q_L(X_n)$ and $x_n$.
Also, we decompose
\smallskip
$$
\eqalign{
\widetilde Q_{k_p,k_p}(X_n)&=e_{k_p}\bigl(x^2_1,\ldots,x^2_n\bigr) \cr
&=e_{k_p}\bigl(x^2_1,\ldots,x^2_{n-1}\bigr)+
e_{k_p-1}\bigl(x^2_1,\ldots,x^2_{n-1}\bigr) x^2_n \cr
&=\widetilde Q_{k_p,k_p}(X_{n-1})+\widetilde Q_{k_p-1,k_p-1}(X_{n-1})x^2_n.\cr}
$$

\smallskip
\noindent
We then repeat 1 and 2 with the so obtained $\widetilde Q_I(X_{n-1})$'s,
thus extracting $x_{n-1}$;
then, we proceed similarly with the so obtained
$\widetilde Q_I(X_{n-2})$'s etc.

\medskip\noindent
Note that if we stop this procedure after extracting the variables
$x_n,x_{n-1},\ldots,x_{m+1}$ we get a development:
$$
\widetilde Q_I(X_n)=\sum\limits_J \widetilde Q_J(X_m) F_J(x_{m+1},...,x_n) \ ,
\tag *
$$
where the sum is over $J\subset I$ (this follows from the linearity formula;
$J$ are not necessary strict).
\smallskip

It follows from the above algorithm that $\widetilde Q_I(X_n)$ is a positive
sum of monomials. It is, in general, not a positive sum of Schur
$S$-polynomials in $X_n$ (we refer the reader to [Mcd1, I] and [L-S1]
for a definition and properties of Schur $S$-polynomials).
Here comes an example computed with the help
of SYMMETRICA [K-K-L].

\smallskip

\noindent
{\bf Example 4.6.} \ Let $n=5$, $\widetilde Q_I=\widetilde Q_I(X_5)$ and
$s_J=s_J(X_5)$. We have:

$$ \widetilde Q_{54}=s_{22221}  \ \ \ \  \widetilde Q_{53}=s_{22211} \ \ \ \
\widetilde Q_{52}=s_{22111}
\ \ \ \  \widetilde Q_{51}=s_{21111} $$
$$ \widetilde Q_{43}=s_{2221}-s_{22111} \ \ \ \ \widetilde
Q_{42}=s_{2211}-s_{21111} \ \ \ \
\widetilde Q_{41}=s_{2111}-s_{11111} $$
$$
\widetilde Q_{32}=s_{221}-s_{2111}+s_{11111} \ \ \ \ \widetilde
Q_{31}=s_{211}-s_{1111}  $$
$$
\widetilde Q_{21}=s_{21}-s_{111} $$
$$
\widetilde Q_{543}=s_{33321}-s_{33222} \ \ \ \ \widetilde
Q_{542}=s_{33221}-s_{32222} \ \ \ \
\widetilde Q_{541}=s_{32221}-s_{22222} $$
$$
\widetilde Q_{532}=s_{33211}-s_{32221}+s_{22222} \ \ \ \  \widetilde
Q_{531}=s_{32211}-s_{22221} $$
$$
\widetilde Q_{521}=s_{32111}-s_{22211} $$
$$
\widetilde Q_{432}=s_{3321}-s_{3222}-s_{33111} \ \ \
\widetilde Q_{431}=s_{3221}-s_{32111}
-s_{2222}$$
$$
\widetilde Q_{421}=s_{3211}-s_{31111}-s_{2221} $$
$$
\widetilde Q_{321}=s_{321}-s_{3111}-s_{222} $$
$$
\widetilde Q_{5432}=s_{44321}-s_{44222}-s_{43331} \ \ \ \  \widetilde
Q_{5431}=s_{43321}-
s_{43222}-s_{33331} $$
$$
\widetilde Q_{5421}=s_{43221}-s_{42222}-s_{33321} $$
$$
\widetilde Q_{5321}=s_{43211}-s_{42221}-s_{33311} $$
$$
\widetilde Q_{4321}=s_{4321}-s_{43111}-s_{4222}-s_{3331}+s_{32221}-2s_{22222}$$
$$
\widetilde Q_{54321}=s_{54321}-s_{54222}-s_{53331}-s_{44421}+s_{43332}-
2s_{33333} .$$
\bigskip

We denote by $S\Cal P(X_n)$ the ring of symmetric polynomials in $X_n$.

\proclaim {\bf Proposition 4.7}
The set $\bigl\{\widetilde Q_I(X_n)\bigr\}$ indexed by all partitions such
that $i_1\le n$ forms an additive basis of $S\Cal P(X_n)$. Moreover, for
any commutative ring $\Cal R$, the same set is a basis of a free
$\Cal R$-module $S\Cal P(X_n)\otimes \Cal R$.
\endproclaim
\demo{Proof}
Let us compare the family $\{\widetilde Q_I(X_n)\}$ with the $\Cal R$-basis
$\{e_I(X_n)\}$
of $S\Cal P(X_n)\otimes\Cal R$,
where $I$ runs over all partitions such
that $i_1\le n$ (and for a partition $I=(i_1,\ldots,i_k)$ we write $e_I(X_n)
=e_{i_1}(X_n)\cdot\ldots\cdot e_{i_k}(X_n)$).
Consider the reverse lexicographic ordering of [Mcd1, I] on the set
of such partitions,
inducing a linear ordering on the above set $\{e_I(X_n)\}$\,.
By the definition of $\widetilde Q_I(X_n)$, one has
$$\widetilde Q_I(X_n)=
e_I(X_n)+(\hbox{combination of earlier monomials in the\ }e_i(X_n)\hbox{'s}).
$$
Since $\{e_I(X_n)\}$ is a $\Cal R$-basis of $S\Cal P(X_n)\otimes\Cal R$\,,
$\{\widetilde Q_I(X_n)\}$ forms another $\Cal R$-basis of
$SP(X_n)\otimes\Cal R$\,.
\qed
\enddemo

\medskip

\proclaim{\bf Corollary-Definition 4.8}
For every $m\leqslant n$ and any partitions $J\subset I$,
there exist uniquely
defined polynomials
$\widetilde Q_{I/J}(x_{m+1},\ldots,x_n)\in
S\Cal P(x_{m+1},\ldots,x_n)$
such that the following equality holds
$$
\widetilde Q_I(X_n)=
\sum\limits_{J\subset I}\widetilde Q_J(X_m)\widetilde Q_{I/J}\bigl(x_{m+1},
\ldots,x_n\bigr).
$$
\endproclaim
\demo{Proof}
The existence of such polynomials $Q_{I/J}(x_{m+1},\dots ,x_n)\in
\Bbb Z[x_{m+1},\dots ,x_n]$ follows from
the above discussion and (*); we put \ $\widetilde Q_{I/J}:=F_J$.

Since $S\Cal P(X_n)\subset S\Cal P(X_m)\otimes S\Cal P(x_{m+1},\dots ,x_n)$ and
$\{\widetilde Q_J(X_m)| \ j_1\le m \}$ is a $\Bbb Z$-basis of $S\Cal P(X_m)$
(Proposition 4.7), we have the corresponding development:
$$
\widetilde Q_I(X_n)=\sum\limits_J \widetilde Q_J(X_m) G_J(x_{m+1},...,x_n) \ .
$$
Using Proposition 4.6 once again with $n$ replaced by $m$ and
$\Cal R=\Bbb Z[x_{m+1},\ldots, x_n]$, we infer that $F_J(=G_J)$ are symmetric
in $x_{m+1},\ldots, x_n$ (and defined uniquely).
\qed
\enddemo

\medskip

We will need also a family of $\widetilde P$-polynomials in
$S\Cal P(X_n)\otimes\Bbb Z[1/2]$ defined by
$\widetilde P_I(X_n):=2^{-l(I)}\widetilde Q_I(X_n)$ for a partition $I$\,.
Also, in analogy to the above, for every $m\le n$ and any partitions
$J\subset I$ there exists uniquely defined polynomials
$\widetilde P_{I/J}(x_{m+1},\ldots,x_n)\in S\Cal P(x_{m+1},\ldots,x_n)
\otimes\Bbb Z[1/2]$
such that
$$
\widetilde P_I(X_n)=\sum\limits_{J\subset I}\widetilde P_I(X_m)\,
\widetilde P_{I/J}(x_{m+1},\ldots,x_n)\,.
$$
$\widetilde P$-polynomials satisfy properties which can be automatically
gotten from the above established properties of $\widetilde Q$-polynomials.
For instance, an analogue of Proposition 4.1 for $\widetilde P$-polynomials
reads:
$$
\widetilde
P_I(X_n)=\sum\limits_{j=0}^{l(I)}x_n^j\Bigl(\,\sum\limits_{|I|-|J|=j}
2^{l(J)-l(I)}\widetilde P_J(X_{n-1})\Bigr)\,,
$$
the sum as in Proposition 4.1.

Given a rank $n$ vector bundle $E$
with the Chern roots $(r_1,\ldots,r_n)$
we set $\widetilde Q_IE:=\widetilde Q_I(X_n)$ and
$\widetilde P_IE:=\widetilde P_I(X_n)$
with
$x_i$ specialized to $r_i$.
Note that this notation is consistent with that used in Section 2 and 3.
Similarly, given a subbundle $E'\subset E$
with the Chern roots $(r_{m+1},\ldots, r_n)$,
we define $\widetilde Q_{I/J}E'=\widetilde Q_{I/J}(x_{m+1},\ldots, x_n)$
and  $\widetilde P_{I/J}E'=\widetilde P_{I/J}(x_{m+1},\ldots, x_n)$
with $x_i$ specialized to $r_i$.

In the next section we will need the following Pieri-type formula for the
$\widetilde Q_I$'s.

\smallskip
\proclaim{\bf Proposition 4.9} Let $I=(i_1,...,i_k)$ be a strict partition of
length $k$. Then

$$ \widetilde Q_I(X_n) \widetilde Q_r(X_n)=
\sum 2^{m(I,r;J)}\widetilde Q_J(X_n),$$
where the sum is over all partitions (i.e. not necessary strict) $J\supset I$
such that
$|J|=|I|+r$ and $J/I$ is a horizontal strip. Moreover, $m(I,r;J)
=card \{ 1\leq p \leq k | \ j_{p+1}<i_p<j_p \}$ or, equivalently, it is
expressed as the number of connected components of the strip $J/I$
not meeting the first column.
\endproclaim
\noindent
(A skew diagram $D$ is connected if each of the sets $\{i:\exists_j
(i,j)\in D\}$ and $\{j:\exists_i (i,j)\in D\}$ is an interval in $\Bbb Z$.)
\demo{Proof}
Let after [L-L-T], $Q_I'(X_n;q)$ denote the Hall--Littlewood
polynomial $ Q_I(Y;q)$ where the alphabet $Y$ is equal to $X_n/(1-q)$
(in the sense of $\lambda$-rings). Using raising operators $R_{ij}$
([Mcd1,I]
we have (see, e.g., [D-L-T])
$$Q_I'(X_n;q)= \prod_{i<j} (1-qR_{ij})^{-1}s_I(X_n)$$
Specialize $q=-1$ and invoke the well known Jacobi--Trudi formula
$$s_I(X_n)=\prod_{i<j} (1-R_{ij})h_I(X_n)$$
where $h_I(X_n)$ is the product of complete homogeneous symmetric
polynomials in $X_n$
associated with the parts of $I$. We have
$$ Q_I'(X_n;-1)=\prod_{i<j}{1-R_{ij} \over 1+R_{ij}} h_I(X_n).$$
Therefore, denoting by $\omega$ the Young duality-involution we get
$\widetilde Q_I(X_n)= \omega \bigl(Q_I'(X_n;-1)\bigr)$.

The required assertion now follows by an appropriate
specialization of the Pieri-type
formula for Hall--Littlewood polynomials  ([Mo], [Mcd1, III.3.(3.8)]).
\qed
\enddemo

\bigskip\smallskip

\centerline{\bf 5. Divided differences and isotropic Gysin maps;}
\centerline{\bf orthogonality of $\widetilde Q$-polynomials}

\medskip
Let $V\to X$ be a vector bundle of rank $2n$ endowed with a nondegenerate
symplectic form.
Let $\pi:{L}G_n(V)\to X$ and $\tau:{L}Fl(V)\to X$ denote respectively
the Grassmannian bundle parametrizing Lagrangian subbundles of $V$ and the flag
bundle parametrizing flags of rank $1$, rank $2$, ..., rank $n$
Lagrangian subbundles of $V$.
We have $\tau=\pi\circ\omega$ where $\omega:{L}Fl(V)\to{L}G_n(V)$
is the projection map. The main goal of this section is to derive several
formulas for the Gysin map $\pi_*:A^*({L}G_n(V))\to A^*(X)$.

We start by recalling the Weyl group $W_n$ of type $C_n$.
This group is isomorphic to $S_n\ltimes \zz_2^n$.
We write a typical element of $W_n$ as $w=(\sigma,\tau)$ where
$\sigma\in S_n$ and $\tau \in\zz_2^n$; so that if
$w'=(\sigma',\tau')$ is another element, their product in $W_n$ is
$w\cdot w'=(\sigma\circ\sigma',\delta)$ where "$\circ$" denotes the
composition of permutations and
$\delta_i=\tau_{\sigma'(i)}\cdot\tau_i'$.
To represent elements of $W_n$ we will use the standard "barred-permutation"
notation, writing them as permutations equipped with bars on those places
(numbered with "$i$") where $\tau_i=-1$.
Instead of using a standard system of generators of $W_n$ given by simple
reflections $s_i=(1,2,\ldots,i+1,i,\ldots,n), \ 1\leqslant i \leqslant n-1$,
and $s_n=(1,2,\ldots,n-1,\nkr n)$, we will use the following system of
generators
$S=\bigl\{s_o=(\nkr 1,2,\ldots,n),\ \ s_1,\ldots,s_{n-1}\bigr\}$
corresponding to the basis: $(-2\varepsilon_1), \varepsilon_1-\varepsilon_2,
\varepsilon_2-\varepsilon_3, ... ,\varepsilon_{n-1}-\varepsilon_n$.
It is easy to check that $(W_n,S)$ is a Coxeter system of type $C_n$.
This "nonstandard" system of generators has several advantages
over the standard one: it leads to easier reasonings by induction on $n$
and the divided differences associated with it produce "stable" symplectic
Schubert type polynomials (for the details concerning the latter
topic - consult a recent work of S. Billey and M. Haiman \cite {B-H}).
Let us record first the formula for the length of an element
$w=(\sigma,\tau)\in W_n$ w.r.t. $S$.
This formula can be proved by an easy induction on $l(w)$ and we leave this to
the reader.

\proclaim{\bf Lemma 5.1}
$l(w)=\sum\limits^n_{i=1}a_i+\sum\limits_{\tau_j=-1}(2b_j+1)$, where
$a_i:=card\bigl\{p | \  p>i\ \&\ \sigma_p<\sigma_i\bigr\}$ and
$b_j:=card\bigl\{p | \  p<j\ \&\ \sigma_p<\sigma_j\bigr\}$.
\endproclaim

In the sequel, whenever we will speak about the "length" of an element $w\in
W_n$,
we will have in mind the length w.r.t. $S$.
\smallskip

Let $X_n=(x_1,\ldots,x_n)$ be a sequence of indeterminates.

We now define symplectic divided differences
$\partial_i : \zz[X_n]\to \zz[X_n],\ i=0,1,\ldots,n-1$, setting

$\Odst\Odst \partial_0(f) = (f-s_0f)/(-2x_1)$ \ ,
\smallskip
$\Odst\Odst\partial_i(f)=(f-s_if)/(x_i-x_{i+1})\quad,\qquad i=1,\ldots,n-1
\quad,$
\medskip\noindent
where $s_0$ acts on $\zz[X_n]$ by sending $x_1$ to $-x_1$ and $s_i$ -- by
exchanging $x_i$ with $x_{i+1}$ and leaving the remaining variables
invariant.
For every $w\in W_n$, $l(w)=l$, let $s_{i_1}\cdot\ldots\cdot s_{i_l}$
be a reduced decomposition w.r.t. S.
Following the theory in [B-G-G] and [D1,2] we define
$\partial_w:=\partial_{s_{i_1}}\cdot\ldots\cdot\partial_{s_{i_l}}$.
By loc.cit. we get a well-defined operator of degree $-l(I)$ acting on
$\zz[X_n]$ (here, "well-defined" means: independent on the reduced
decomposition chosen).
\bigskip

We want first to study the operator $\partial_{w_o}$ where
$w_o=(\nkr 1,\nkr 2,\ldots,\nkr n)$
is the maximal length element of $W_n$.
To this end we need
some preliminary considerations.

Let $Q\Cal P(X_n)$ denote the ring of Schur's $Q$-polynomials in $X_n$.
We record the following (apparently new) identity in this ring.
In 5.2 -- 5.4 below we
will write: $e_i=e_i(X_n)$, $s_I=s_I(X_n)$, $Q_I=Q_I(X_n)$ and $\widetilde
Q_I=\widetilde Q_I(X_n)$ for brevity.

\medskip

\proclaim{\bf Proposition 5.2}
In $Q\Cal P(X_n)$,
$$Q_{\rho_k}=Det(a_{i,j})_{1\le i,j\le k},$$
where
$a_{i,j}=Q_{k+1+j-2i}$ if $k+1+j-2i\ne 0$
(with $Q_i=0$ for $i<0$) and $a_{i,j}=2$ if
$k+1+j-2i=0$.
\endproclaim

\demo{Proof}
We have from the theory of symmetric polynomials (see [P2] and the references
therein),
$$
Q_{\rho_k}=2^ks_{\rho_k}=Det\Bigl(2e_{k+1+j-2i}\Bigr)_{1\le
i,j\le k}.
$$

It follows from the Pieri formula (see [Mcd] and [LS1]) that
$$
\Biggl(2\sum\limits_{hooks\ I,\atop |I|=i-2}s_I\Biggr)\cdot s_2 + 2e_i =
2\sum\limits_{hooks\ J,\atop |J|=i}s_J.
$$
Hence,
by multiplying the $p$-th row by $s_2$ and adding it to the $(p-1)$-th
one successively for $p=k,k-1,\ldots ,3,2$, the latter determinant is
rewritten
in the form:
$$
Det\Biggl(2\sum\limits_{hooks\ I,\atop |I|=k+1+j-2i}s_I
\Biggr)_{_{1\le i,j\le k}}.
$$
Notice that the degree 0 entries in this determinant are equal to 2 and the
negative degree entries vanish.
Since $Q_i= 2\sum\limits_{hooks\ I,|I|=i}s_I$, the assertion
follows.
\qed
\enddemo
\bigskip

Let $\Cal J$ be the ideal in $S\Cal P(X_n)$ generated
by $e_i(x_1^2,\ldots,x_n^2)$, $1\le i \le n$.
We now invoke a corollary of [P2, Theorem 6.17] combined with [B-G-G,
Theorem 5.5] and [D2, 4.6(a)]: there is a ring isomorphism $S\Cal P(X_n)/
\Cal J \to Q\Cal P(X_n)/\oplus \Bbb Z Q_I(X_n)$, where $I$ runs over
all strict partitions $I\not\subset \rho_n$, given by $e_i(X_n)\mapsto
Q_i(X_n)$
(see the remark after Theorem 6.17 in [P2, pp.181-182]).
We thus get from the proposition:
\medskip

\proclaim{\bf Corollary 5.3}
In $S\Cal P(X_n)$, $\widetilde Q_{\rho_k}$ is congruent to
$Det(b_{i,j})_{1\le i,j\le k}$ modulo $\Cal J$, where
$b_{i,j}=e_{k+1+j-2i}$ \ if $k+1+j-2i\ne 0$ (with $e_i=0$ for $i<0$)
and $b_{i,j}=2$ \ if $k+1+j-2i=0$.
\endproclaim

\medskip

We now state:

\proclaim{\bf Lemma 5.4}
In $S\Cal P(X_n),\ \widetilde Q_{\rho_n}\equiv e_n e_{n-1}\ldots e_1
\equiv s_{\rho_n}$
(mod $\Cal J$).
\endproclaim

\demo{Proof}
By the corollary it is sufficient to prove that
$Det(b_{i,j})_{1\le i,j\le n}\equiv e_n e_{n-1}\ldots e_1 \equiv
s_{\rho_n}$ (mod $\Cal J$).
Recall that $s_{\rho_n}=Det(c_{i,j})_{1\le i,j\le n}$ where
$c_{i,j}=e_{n+1+j-2i}$ if $n+1+j-2i\ne 0$
and $c_{i,j}=1$ if $n+1+j-2i=0$, i.e. the matrices $(b_{i,j})$ and $(c_{i,j})$
are the same modulo the degree 0 entries.

Let us write the determinants $Det(b_{i,j})$ and $Det(c_{i,j})$
as the sums of the standard $n!$ terms
(some of them are zero).
It is easy to see that apart from the "diagonal" term $e_ne_{n-1}
\ldots e_1,$
every other term appearing in both sums is divisible by
$e_ne_{n-1}\ldots e_{p+1}e_p^2
$
for some $p\ge 1$. We claim that,
$e_ne_{n-1}\ldots e_{p+1}e_p^2\in \Cal J.
$
Indeed, $e_n^2\in\Cal J$ and suppose, by induction, that we have shown
$e_ne_{n-1}\ldots e_{q+1}e_q^2\in \Cal J
$
for $q>p$. Then
$$
e_ne_{n-1}\ldots e_{p+1}e_p^2=
e_ne_{n-1}\ldots e_{p+1}\Bigl[\widetilde Q_{p,p}
+2\sum\limits^p_{i=1}(-1)^{i-1}e_{p+i}e_{p-i}\Bigr]
$$
belongs to $\Cal J$ by the induction assumption,
because $\widetilde Q_{p,p}\in \Cal J$
(see Proposition 4.2).
This shows that
$$Det(b_{i,j})\equiv e_ne_{n-1}\ldots e_1\equiv Det(c_{i,j}) \
(mod \ \Cal J).$$

Thus the lemma is proved.
\qed
\enddemo

\bigskip

The following known result (see, e.g., [D1] where the result is given
also for other root systems) is accompanied by a proof
for the reader's convenience.

\proclaim{\bf Proposition 5.5}
One has for $f\in\zz[X_n]$,
$$
\partial_{w_o}(f)=(-1)^{n(n+1)/2}
\Bigl(2^nx_1\cdot \ldots \cdot x_n \prod\limits_{i<j}
(x^2_i-x^2_j)\Bigr)^{-1} \sum\limits_{w\in W_n}(-1)^{l(w)}w(f).
$$
\endproclaim

\demo{Proof}
By the definition of $\partial_{w_o}$ we infer that
$\partial_{w_o}=\sum\limits_{w\in W_n}\alpha_ww$ where the coefficients
$\alpha_w$ are rational functions in $x_1,\ldots,x_n$.
Since $w_o$ is the maximal length element in $W_n$,
$\partial_i\circ\partial_{w_o}=0$ for all $i=0,1,\ldots,n-1$.
Consequently $s_i\partial_{w_o}=\partial_{w_o}$ for $i=0,1,\ldots,n-1$
and hence $v\partial_{w_o}=\partial_{w_o}$ for all $v\in W_n$.
In particular, for every $v\in W_n$,
$\partial_{w_o}=\sum\limits_{w\in W_n}v(\alpha_w)vw$.
Thus $\alpha_{vw}=v(\alpha_w)$ for all $v,w\in W_n$, and we see that, e.g.,
$\alpha_{w_o}$ determines uniquely all the $\alpha_w$'s.

\smallskip
\noindent
\underbar{Claim} \
$\alpha_{w_o}=(-1)^{n(n-1)/2}\Bigl(2^nx_1\cdot \ldots \cdot
x_n\prod\limits_{i<j}
(x^2_i-x^2_j)\Bigr)^{-1}$.

\smallskip\noindent
{\it Proof of the claim}:
Denote now the maximal length element in $W_n$ by $w_o^{(n)}$.
We argue by induction on $n$.
For $n=1$, we have $\alpha_{w_o^{(1)}}={1\over2x_1}$.
We now record the following equality:
$$
w_o^{(k+1)}=s_k\cdot s_{k-1}\cdot\ldots\cdot s_1\cdot s_0\cdot s_1\cdot\ldots
\cdot s_{k-1}\cdot s_k\cdot w_o^{(k)},
$$
that implies
$$
\partial_{w_o^{(k+1)}}=\partial_k\circ\partial_{k-1}\circ\ldots
\circ\partial_1\circ\partial_0\circ\partial_1\circ\ldots
\circ\partial_{k-1}\circ\partial_k\circ\partial_{w_o^{(k)}}.
$$
It follows easily from the latter equality that
$$
\alpha_{w_o^{(k+1)}}=(-1)^{k}\Bigl(2x_{k+1}\prod\limits_{i\leqslant k}
(x_i-x_{k+1})
\prod\limits_{i\leqslant k}(x_i+x_{k+1})\Bigr)^{-1}\alpha_{w_o^{(k)}}.
$$

This allows us to perform the induction step $n\to n+1$,
thus proving the claim.

Finally, for arbitrary $w\in W_n$, $$\alpha_w=ww_o(\alpha_{w_o})=
(-1)^{n(n+1)/2+l(w)}\Bigl(2^nx_1\cdot \ldots \cdot x_n\prod\limits_{i<j}
(x^2_i-x^2_j)\Bigr)^{-1}$$
because, for $w=(\sigma,\tau)$, $l(\sigma,\tau)=\sum a_i+\sum_{\tau_j=-1}
(2b_j+1)\equiv
l(\sigma)+card\{p| \ \tau_p=-1\}$ $(mod \ 2)$ (see Lemma 5.1).
\qed
\enddemo

\proclaim{\bf Corollary 5.6}
(i) $\partial_{w_o}(x_1^{\alpha_1}x_2^{\alpha_2}\ldots x_n^{\alpha_n})=0$ if
$\alpha_p$ is even for
some $p=1,\ldots,n$.

\noindent
(ii) If all $\alpha_p$ are odd then
$$\partial_{w_o}(x_1^{\alpha_1}x_2^{\alpha_2}\ldots x_n^{\alpha_n})
= (-1)^{n(n+1)/2} s_{\rho_n}(X_n)^{-1}
\partial(x_1^{\alpha_1}x_2^{\alpha_2}\ldots x_n^{\alpha_n}),$$
where here and in the sequel $\partial$ denotes
the Jacobi symmetrizer $$\bigl(\sum_{\sigma \in S_{n}}(-1)^{l(\sigma)}
\sigma(-)\bigr)
/\prod\limits_{i<j}(x_i-x_j).$$
\endproclaim

\demo{Proof}
(i) Let us fix $\sigma\in S_n$ and look at all elements of $W_n$ of the form
$(\sigma,\tau)$ where $\tau \in\zz_2^{\strut n}$. Then, writing $x^\alpha$
for $x_1^{\alpha_1}\cdot \ldots \cdot x_n^{\alpha_n}$, we have
$$
\sum\limits_{\tau}(-1)^{l(\sigma,\tau)}(\sigma,\tau)x^\alpha=
(-1)^{l(\sigma)}\sigma(x^\alpha)\sum\limits_{\tau}(-1)^{card\{p| \ \tau_p=-1\}}
\tau_1^{\alpha_1}\ldots\tau_n^{\alpha_n},
$$
because $l(\sigma,\tau) \equiv
l(\sigma)+card\{p| \ \tau_p=-1\}$ $(mod \ 2)$.
Suppose that some numbers among $\alpha_1,\ldots,\alpha_n$ are even.
We will show that this implies
$$
\sum\limits_{\tau}(-1)^{card\{p| \ \tau_p=-1\}}
\tau_1^{\alpha_1}\ldots\tau_n^{\alpha_n}=0.
$$
We can assume that $\alpha_1,\ldots,\alpha_k$ are odd and $\alpha_{k+1},
\ldots,\alpha_n$ are
even for some $k<n$ (by permuting the $\tau_p$'s if necessary).
We have
\medskip
$
\sum\limits_{\tau}(-1)^{card\{p| \ \tau_p=-1\}}
\tau_1^{\alpha_1}\ldots\tau_n^{\alpha_n}=
$
\medskip\Odst\Odst
$
= \sum\limits_{\tau}(-1)^{card\{p| \ \tau_p=-1\}}
(-1)^{card\{p| \ \tau_p=-1,\ p\le k\}}
$
\medskip\Odst\Odst
$
= \sum\limits_{\tau}(-1)^{card\{p| \ \tau_p=-1,\ p>k\}}
$
\medskip\Odst\Odst
$
= 2^k\sum\limits_{i=0}^{n-k}(-1)^i{n-k\choose i}=2^k(1-1)^{n-k}=0.
$

\medskip

\noindent
(ii) Let us now compute $\partial_{w_o}(x_1^{\alpha_1}\ldots x_n^{\alpha_n})$
where all
$\alpha_p$ are odd. Then
$$
\sum\limits_{\tau}(-1)^{card\{j| \ \tau_j=-1\}}
\tau_1^{\alpha_1}\ldots\tau_n^{\alpha_n}=2^n,\hbox{ and}
$$
$$
\eqalign{
\partial_{w_o}(x^\alpha)
 &= (-1)^{n(n+1)/2} \Bigl(2^nx_1\ldots
x_n\prod\limits_{i<j}(x^2_i-x^2_j)\Bigr)^{-1}
 2^n\sum\limits_{\sigma\in S_n}(-1)^{l(\sigma)}\sigma(x^\alpha) \cr
 &= (-1)^{n(n+1)/2} s_{\rho_n}(X_n)^{-1}\partial(x^\alpha). \qed}$$
\enddemo

\bigskip

We now record the following properties of the operator
$\tk =\partial_{(\bar n,\ldots,\bar 2,\bar 1)}$.
In the following, let $\Cal I = \Cal J\Bbb Z[X_n]$.

\proclaim {\bf Lemma 5.7} (i) If $f\in S\Cal P(x^2_1,\ldots,x^2_n)$ then
$\tk(f\cdot g)=f\cdot\tk(g)$.

\noindent
(ii) $\tk\bigl(\widetilde Q_{\rho_n}(X_n)\bigr)=(-1)^{n(n+1)/2}$.
\endproclaim
\smallskip

\demo{Proof} (i) This assertion is clear because every polynomial in
$S\Cal P(x^2_1,\ldots,x^2_n)$ is $W_n$-invariant.
Observe that it implies that if $f\equiv g$ (mod $\Cal I$) then
$\tk(f)\equiv \tk(g)$ (mod $\Cal I$).
\smallskip
\noindent
(ii) (This can be also deduced from the Chow ring of the Lagrangian
Grassmannian. We present here a direct algebraic argument.)
In this part we will use the following properties of the Jacobi
symmetrizer $\partial$ \
(see [L-S2], [Mcd2]):
\smallskip

{\parindent=20pt
\item{1.} If $f\in S\Cal P(X_n)$, $g\in\Bbb Z[X_n]$ then
$\partial(f\cdot g)=f\cdot\partial(g)$.
\smallskip
\item{2.} For any $\alpha=(\alpha_1,\ldots,\alpha_n)\in\Bbb N^n$,
$\partial (x^{\alpha})=s_{\alpha-\rho_{n-1}}(X_n)$.
In particular, if $\alpha_i=\alpha_j$ for some $i\ne j$ then
$\partial (x^{\alpha})=0$.
\smallskip
\item{3.} $\partial=\partial_{(n,n-1,\ldots,1)}.$
\par}

\smallskip

Let $e_i=e_i(X_n)$.
Since $\widetilde Q_{\rho_n}(X_n)\equiv e_ne_{n-1}\ldots e_1$ (mod $\Cal I$)
 \ (by Lemma 5.4), we have
$$
\tk(\widetilde Q_{\rho_n}(X_n))=\tk(e_ne_{n-1}\ldots e_1)=
(\tk\circ\partial)(x^{\rho_{n-1}}e_ne_{n-1}\ldots e_1).
$$
by properties 1 and 2 above.
Since
$$
(\mathstrut\overline n,\overline{n-1},\ldots,\overline 1)\circ
(n,n-1,\ldots,1)=w_0,
$$
the latter expression equals
$\partial_{w_0}(x^{\rho_{n-1}}e_ne_{n-1}\ldots e_1)$ by property 3.
The degree of the polynomial $x^{\rho_{n-1}}e_ne_{n-1}\ldots e_1$ is $n^2$.
Assuming that $\alpha_1+\ldots+\alpha_n=n^2$, we have
$\partial_{w_0}(x^{\alpha})\ne 0$ only if
$$
x^{\alpha}=
x^{2n-1}_{\sigma(1)}x^{2n-3}_{\sigma(2)}\ldots x_{\sigma(n)}
$$
for some $\sigma\in S_n$.
Indeed, it follows from Corollary 5.6(i) that $\partial_{w_0}(x^{\alpha})\ne 0$
only if all the $\alpha_i$'s are odd.
Moreover, they must be all different; otherwise $\partial (x^{\alpha})=0$
(and consequently $\partial_{w_o}(x^{\alpha})=0$)
by property 2.
We conclude that $\{\alpha_1,\ldots,\alpha_n\}=\{2n-1,2n-3,\ldots,1\}$.
But there is only one such a monomial $x^{\alpha}$ in
$x^{\rho_{n-1}}e_ne_{n-1}\ldots e_1$, namely the one with
$(\alpha_1,\ldots,\alpha_n)=(2n-1,2n-3,\ldots,1)$.
Therefore
$$
\partial_{w_o}(x^{\rho_{n-1}}e_ne_{n-1}\ldots e_1)=\partial_{w_0}(x_1^{2n-1}
x_2^{2n-3}\ldots x_n)=(-1)^{n(n+1)/2}
$$
by Corollary 5.6(ii) and property 2.
\qed
\enddemo

\bigskip
We now pass to a geometric interpretation of the operator $\tk$.

\proclaim{\bf Proposition 5.8}
Specializing the variables $x_1,\ldots,x_n$
to the Chern roots $r_1,\ldots,r_n$ of the tautological
subbundle $R$ on ${L}G_nV$,  one has the equality
$$
\pi_*\bigl({f}(r_1,\ldots,r_n)\bigr)=
\Bigl(\partial_{(\nkr n,\nkr {n-1},\ldots,
\nkr 2,\nkr 1)}
f\Bigr)(r_1,\ldots,r_n),
$$
where ${f}(-)$ is a polynomial in $n$ variables.
\endproclaim
(A symmetrization operator variant of this proposition follows
also from a recent paper by M. Brion [Br]. We give here a short proof
using only divided differences interpretation of Gysin maps
for complete (usual and Lagrangian) flag bundles.)

\demo{Proof} \ We invoke a result saying that the Gysin map associated
with $\omega$ and $\tau$ is induced by the following
divided differences operators:
$$
\eqalign{
\tau_*\bigl({f}(r_1,\ldots,r_n)\bigr)
&=\Bigl(\partial_{(\nkr 1,\nkr 2,\ldots,\nkr n)}f\Bigr)
(r_1,\ldots,r_n) \ \hbox{ and}  \cr
\omega_*\bigl({f}(r_1,\ldots,r_n)\bigr)
&=\Bigl(\partial_{(n,n-1,\ldots,1)}f\Bigr)(r_1,\ldots,r_n).}
$$
As for the latter equality, see [P1, Sect.2] ,
as for the former
compare [Br] where the author gives a symmetrizing
operator expression for $G/B$-fibrations
(over a point, say, this expression was given in [A-C]).
The needed divided differences
interpretation of those symmetrizing operators follows, e.g., from [D1].

\noindent
Since
$$(\nkr 1,\nkr 2,\ldots,\nkr n)=
(\nkr n,\nkr {n-1},\ldots,\nkr 1)\circ(n,n-1,\ldots,1),$$
we get
$$
\partial_{(\nkr 1,\nkr 2,\ldots,\nkr n)}=
\partial_{(\nkr n,\nkr {n-1},\ldots,\nkr 1)}\circ
\partial_{(n,n-1,\ldots,1)}.
$$
\smallskip

Of course, $\tau_*=\pi_*\circ \omega_*$.
Since $\omega_*$ is surjective, comparison of the latter equation with
the former implies the desired assertion about $\pi_*$.
\qed
\enddemo

We now show how to compute the images via $\pi_*$ of $\widetilde
Q$-polynomials in the Chern classes of $R\hak$. Let us write
$X_n\hak=(-x_1,\ldots,-x_n)$ for brevity.

\medskip
\proclaim{\bf Proposition 5.9}
One has $\triangledown(\widetilde Q_I(X_n\hak)) \ne 0$ only if
the set of parts of
$I$ is equal to
\break
$\{1,2,\ldots,n\}$ and each number $p\ \ (1\leqslant p\leqslant n)$
appears in $I$ with an $\underline{odd}$\  multiplicity $m_p$.
Then, the following equality holds in $\Bbb Z[X_n]$,
$$
\triangledown\bigl(\widetilde Q_I(X_n\hak)\bigr) = \prod\limits^n_{p=1}
e_p(x_1^2,\ldots, x_n^2)^{(m_p-1)/2}.
$$
\endproclaim

\demo{Proof}
By Proposition 4.3 we can express $\widetilde Q_I(X_n\hak)$ as
$$
\widetilde Q_I(X_n\hak) = \widetilde Q_{j_1,j_1}(X_n\hak)\ \ldots \
\widetilde Q_{j_l,j_l}(X_n\hak) \widetilde Q_L(X_n\hak),
$$
where $L$ is a strict partition.
(We divide the elements of the multiset $I$ into pairs of equal elements
and the set $L$ whose elements are all different.)
Some of the $j_p$'s can be mutually equal.

By Proposition 4.2, $\widetilde Q_{j,j}(X_n\hak)=e_j(x_1^2,\ldots,x_n^2)$
is a scalar w.r.t. $\triangledown$.

By Lemma 4.4, $\widetilde Q_L(X_n\hak)\ne 0$ only if $L\subset\rho_n$.
On the other hand, for a strict partition $L\subset \rho_n$,
$\triangledown(\widetilde Q_L(X_n\hak))\ne 0$ only if $L=\rho_n$,
when it is equal to 1 (see Lemma 5.7(ii)).
\smallskip

Putting this information together, the assertion follows.
\qed
\enddemo

Consequently, specializing $(x_i)$ to the Chern roots $(r_i)$ of
the tautological subbundle on $LG_n(V)$ we have

\proclaim {\bf Theorem 5.10}
The element $\widetilde Q_I R\hak$
has a nonzero image under $\pi_*:A^*(LG_nV)\to A^*(X)$
only if each number
$p$, $1\leqslant p\leqslant n$, appears as a part of $I$ with an
odd multiplicity $m_p$. If the latter condition holds then
$$
\pi_* \widetilde Q_IR\hak
= \prod_{p=1}^n \bigl((-1)^p c_{2p}V\bigr)^{(m_p-1)/2}.
$$
\endproclaim

\demo {Proof}
This follows from Proposition 5.9 and the equality
$c_{2p}V = (-1)^p e_p(r_1^2,...,r_n^2)$.
\qed
\enddemo

Our next goal will be to show how to compute the images via $\pi_*$ of
$S$-polynomials in the Chern classes of the tautological Lagrangian bundle.
To this end
we record the following
identity of symmetric polynomials.
We have found this simple and remarkable identity during our work on isotropic
Gysin maps
and have not seen it in the literature.

\proclaim{\bf Proposition 5.11}
For every partition $I=(i_1,\ldots,i_n)$ and any positive integer $p$,
one has in $S\Cal P(X_n)$,
$$
s_I\bigl(x^p_1,\ldots,x^p_n\bigr)\cdot s_{(p-1)\rho_{n-1}}(X_n) =
s_{pI+(p-1)\rho_{n-1}}(X_n).
$$
Here, given a partition $I=(i_1,i_2,\ldots)$, we
write $pI=(pi_1,pi_2,\ldots)$.
\endproclaim

\demo{Proof}
We use the Jacobi presentation of a Schur polynomial as a ratio of two
alternants (see [Mcd1], [L-S1]).
We have:
\smallskip
$$
\eqalign{
s_I\bigl(x_1^p,\ldots,x_n^p\bigr)
&={Det\bigl(x_k^{(i_l+n-l)p}\bigr)_{1\leqslant k,l\leqslant n}
 \over
 Det\bigl(x_k^{p(n-l)}\bigr)_{1\leqslant k,l\leqslant n} } \cr
&\cr
&={
 Det\bigl(x_k^{pi_l+(n-l)(p-1)+(n-l)}\bigr)_{1\leqslant k,l\leqslant n}
 \over
 Det\bigl(x_k^{n-l}\bigr)_{1\leqslant k,l\leqslant n} \cdot
 \Biggl({
   Det\bigl(x_k^{(p-1)(n-l)+(n-l)}\bigr)_{1\leqslant k,l \leqslant n}
   \over
   Det\bigl(x_k^{n-l}\bigr)_{1\leqslant k,l\leqslant n}
 }\Biggr)
 } \cr
&\cr
&={s_{pI+(p-1)\rho_{n-1}}(X_n) \over s_{(p-1)\rho_{n-1}}(X_n)}. \qed \cr}
$$
\enddemo

\proclaim {\bf Corollary 5.12} \ For $p=2$ we get
$$
s_I(x_1^2,\ldots,x_n^2)\cdot s_{\rho_{n-1}}(X_n) = s_{2I+\rho_{n-1}}(X_n).
$$
\endproclaim

\smallskip
\noindent
(For another derivation of this identity with the help of Quaternionic
Grassmannians see Appendix A.)
\medskip

We now give a geometric translation of the latter formula,
or rather its consequence
$$
s_I(x_1^2,\ldots,x_n^2)\cdot s_{\rho_n}(X_n)=s_{\rho_n+2I}(X_n). \leqno(*)
$$
\medskip

\proclaim{\bf Theorem 5.13}
The element $s_IR\hak$ has a nonzero image under $\pi_*$ only if
the partition $I$ is of the form $2J+\rho_n$ for some partition $J$.
If $I=2J+\rho_n$ then
$$
\pi_*s_IR\hak = s_J^{^{[2]}}V \quad,
$$
where the right hand side is defined as follows:
if $s_J={P}(e.)$ is a unique presentation of $s_J$ as a polynomial
in the elementary symmetric functions $e_i$, $E \ -$ a vector bundle,
then $s_J^{[2]}(E) := {P}$ with $e_i$ replaced by $(-1)^ic_{2i}E\quad
(i=1,2,\ldots)$.
\endproclaim

\demo{Proof}
Since $s_IR\hak=\omega_*(q^{I+\rho_{n-1}})$
where $q=(q_1,\ldots,q_n)$ are the Chern roots of $R\hak$
(this is a familiar Jacobi-Trudi formula restated using the Gysin
map for the flag bundle -- see, e.g., [P3] and the references
therein),
we infer from Corollary 5.6(i) that
$s_IR\hak$ has a nonzero image under $\pi_*$ only if all parts of
$I+\rho_{n-1}$ are odd.
This implies that $l(I)=n$ and $I$ is strict thus of the form $I'+\rho_n$
for some partition $I'$.
Finally all parts of $I'+\rho_n+\rho_{n-1}$ are odd iff $I'=2J$ for some
partition $J$, as required.

Assume now that $I=2J+\rho_n$ and specialize the identity (*) by replacing
the variables $(x_i)$ by the Chern roots $(q_i)$.
The claimed formula now follows since: $s_I(q_1^2,\ldots,q_n^2)$ is
a scalar w.r.t. $\pi_*$, $\pi_*s_{\rho_n}(q_1,\ldots,q_n)=1$
by Lemma 5.7(ii) combined with Lemma 5.4; finally
$(-1)^ic_{2i}V=e_i(q_1^2,\ldots,q_n^2)$ because
of Lemma 1.1(2).
\qed
\enddemo
\medskip

Observe that the theorem contains an explicit calculation of the
ratio in Corollary 5.6(ii).

\bigskip

We now pass to the odd orthogonal case.
The Weyl group $W_n$ of type $B_n$.
is isomorphic to $S_n\ltimes \zz_2^n$ and its elements are
"barred-permutations".
We use the following system of
generators of $W_n$:
$S=\bigl\{s_o=(\nkr 1,2,\ldots,n),\ \ s_1,\ldots,s_{n-1}\bigr\}$
corresponding to the basis $(-\varepsilon_1), \varepsilon_1-\varepsilon_2,
\varepsilon_2-\varepsilon_3, ... ,\varepsilon_{n-1}-\varepsilon_n$.
Consequently, the divided differences
$\partial_i,\ i=1,\ldots,n-1$, are the same but
$\partial_0(f)=(f-s_0f)/(-x_1)$.

The odd orthogonal analog of Proposition 5.5 reads:
$$
\partial_{w_0}(f)=(-1)^{n(n+1)/2}
\Bigl(x_1\cdot ...\cdot x_n \prod\limits_{i<j}(x^2_i-x^2_j)
\Bigr)^{-1} \sum\limits_{w\in W_n}(-1)^{l(w)}w(f).
$$

Arguing essentially as
in the proof of Proposition 5.8 (with obvious modifications),
one shows that the Gysin map
associated with $\pi: OG_nV\to X$ is induced by the orthogonal divided
difference operator $\partial_{(\nkr n,\nkr {n-1},\ldots,\nkr 1)}$.

\smallskip
The odd orthogonal analog of Theorem 5.10 reads:
\proclaim{\bf Theorem 5.14}
The element $\widetilde Q_IR\hak$
has a nonzero image under $\pi_*:A^*({O}G_nV)\to
A^*(X)$ only if
each number $p$, $1\le p\le n$, appears as a part of $I$
with an odd multiplicity $m_p$.
If the latter condition holds then
$$
\pi_*\widetilde Q_IR\hak=2^{n}
\prod\limits^n_{p=1}\bigl((-1)^pc_{2p}V\bigr)^
{(m_p-1)/2}.
$$
\endproclaim

This holds because the calculation in Proposition 5.9 now goes as
follows: with the notation from the proof of Proposition 5.9,
the polynomial
$$
\widetilde Q_I(X_n\hak)=
2^{n} \widetilde Q_{j_1,j_1}(X_n\hak)\ \ldots \
\widetilde Q_{j_l,j_l}(X_n\hak) \widetilde P_{\rho_n}(X_n\hak)
$$
is mapped via $\partial_{(\nkr n,\nkr {n-1},\ldots,\nkr 1)}$
to
$$
2^{n}\prod\limits_{h=1}^l e_{j_h}(x_1^2,...,x_n^2),$$
since $\partial_{(\nkr n,\nkr {n-1},\ldots,\nkr 1)}
\bigl(\widetilde P_{\rho_n}(X_n\hak)\bigr)=1$. (The proof of the last statement
is the same as that of Lemma 5.7(ii).)
\medskip

Finally, the odd orthogonal analog of Theorem 5.13 reads:

\proclaim{\bf Theorem 5.15}
The element $s_IR\hak$ has a nonzero image under $\pi_*$ only if the
partition $I$ is of the form $2J+\rho_n$ for some partition $J$.
If $I=2J+\rho_n$ then
$$
\pi_*s_IR\hak=2^ns_J^{[2]}V,
$$
where $s_J^{[2]}(-)$ is defined as in Theorem 5.13.
\endproclaim

This holds because $s_{\rho_n}(X_n\hak)$ is congruent to
$2^n\widetilde P_{\rho_n}(X_n\hak)$
modulo $\Cal J$ (Lemma 5.4)
and $\pi_*\widetilde P_{\rho_n}R\hak=1$.
Also, we use Lemma 1.1(2).
\medskip

We now pass to the even orthogonal case.

In type $D_n$ the Weyl group $W_n$ is identified with the
subgroup of the group of "barred permutations" $(w_1,\ldots,w_n)$
whose elements have even number of bars only.
Consider a system $S$ of generators of $W_n$ consisting of \
$s_{\bar 1}=(\nkr 2,\nkr 1,3,\ldots,n)$ \ and \
$s_i=(1,2,\ldots,i-1,i+1,i,i+2,\ldots,n)$,
$i=1,2,\ldots,n-1$.
$(W_n,S)$ is a Coxeter system of type $D_n$ and the length function
w.r.t. $S$ is
$$
l(w)=\sum\limits^n_{i=1}a_i+\sum\limits_{\tau_j=-1}2b_j,
$$
where $a_i=card\{p\ | \ p>i\ \&\ w_p<w_i\}$ and
$b_j=card\{p\ | \ p<j\ \&\ w_p<w_j\}$.
The longest element $w_0$ in $W_n$ is equal to
$(\nkr 1, ... ,\nkr n)$
if $n$ is even and to $(1,\nkr 2,... ,\nkr n)$ if $n$ is odd.
Following \cite {B-G-G} and \cite {D1,2} one defines the operators
$\partial_w:\zz[X_n]\to\zz[X_n]$
for $w\in W_n$ ( here,
$$\partial_{\bar 1}f=\bigl(f-f(-x_2,-x_1,
x_3,...,x_n)\bigr)/(-x_1-x_2). \ )$$

\smallskip
The even orthogonal analog of Proposition 5.5 reads:
$$
\partial_{w_0}(f)=(-1)^{n(n-1)/2}
\prod\limits_{i<j}(x^2_i-x^2_j)
^{-1} \sum\limits_{w\in W_n}(-1)^{l(w)}w(f).
$$

The even orthogonal analog of Corollary 5.6 reads:

\proclaim{\bf Lemma 5.16}
(i) $\partial_{w_o}(x_1^{\alpha_1}x_2^{\alpha_2}\ldots x_n^{\alpha_n})=0$ if
$\alpha_p$ is odd for
some $p=1,\ldots,n$.

\noindent
(ii) If all $\alpha_p$ are even then
$$\partial_{w_o}(x^{\alpha})
= (-1)^{n(n-1)/2}2^{n-1}s_{\rho_{n-1}}(X_n)^{-1}
\partial(x_{\alpha}),$$
where $\partial$ denotes
the Jacobi symmetrizer.
\endproclaim

\demo{Proof}
(i) Let us fix $\sigma\in S_n$ and look at all elements of $W_n$ of the form
$(\sigma,\tau)$ where
$\tau\in \{+1,-1\}^n$ and $\prod_i \tau_i = 1$.
We have
$$
\sum\limits_{\tau}(-1)^{l(\sigma,\tau)}(\sigma,\tau)x^\alpha=
(-1)^{l(\sigma)}\sigma(x^\alpha)\sum\limits_{\tau}
\tau_1^{\alpha_1}\ldots\tau_n^{\alpha_n},
$$
because $l(\sigma,\tau) \equiv
l(\sigma)$ $(mod \ 2)$.
Suppose that some numbers among $\alpha_1,\ldots,\alpha_n$ are odd.
We can assume that $\alpha_1,\ldots,\alpha_k$ are odd and $\alpha_{k+1},
\ldots,\alpha_n$ are
even for some $k<n$ (by permuting the $\tau_p$'s if necessary).
We have
$$
\sum\limits_{\tau}
\tau_1^{\alpha_1}\ldots\tau_n^{\alpha_n}
=\sum\limits_{\tau}(-1)^{card\{p| \ \tau_p=-1,\ p>k\}}
=2^k\sum\limits_{i=0}^{n-k}(-1)^i{n-k\choose i}=0.
$$
\smallskip
\noindent
(ii) Let us now compute $\partial_{w_o}(x_1^{\alpha_1}\ldots x_n^{\alpha_n})$
where all
$\alpha_p$ are even. Then
$$
\sum\limits_{\tau}
\tau_1^{\alpha_1}\ldots\tau_n^{\alpha_n}=2^{n-1},\hbox{ and}
$$
$$
\eqalign{
\partial_{w_o}(x^\alpha)
 &= (-1)^{n(n-1)/2} \prod\limits_{i<j}(x^2_i-x^2_j)^{-1}
 2^{n-1}\sum\limits_{\sigma\in S_n}(-1)^{l(\sigma)}\sigma(x^\alpha) \cr
 &= (-1)^{n(n-1)/2}2^{n-1} s_{\rho_{n-1}}(X_n)^{-1}\partial(x^\alpha). \qed}$$
\enddemo

\bigskip

Let us now denote by $\Cal J$ the ideal in $S\Cal P(X_n)\otimes \Bbb Z[1/2]$
generated by $e_i(x_1^2,\ldots,x_n^2)$, $i=1,\ldots,n-1$, and $x_1\cdot\ldots
\cdot x_n$.
In the following analog of Lemma 5.4
we write $e_i=e_i(X_n)$, $s_I=s_I(X_n)$, $P_I=P_I(X_n)$ and $\widetilde
P_I=\widetilde P_I(X_n)$ for brevity.

\proclaim {\bf Lemma 5.17} \ In $S\Cal P(X_n)\otimes \Bbb Z[1/2]$,
$$
\widetilde P_{\rho_{n-1}}\equiv 2^{-(n-1)}e_{n-1}e_{n-2}\ldots e_1\equiv
2^{-(n-1)}s_{\rho_{n-1}} \ (mod \ \Cal J).$$
\endproclaim

\demo{Proof} Proposition 5.2 implies that
$P_{\rho_k} = Det (a_{i,j})_{1\le i,j,\le k}$ where
$a_{i,j}=P_{k+1+j-2i}$ if $k+1+j-2i\ne 0$
(with $P_i=0$ for $i<0$) and $P_{i,j}=1$ if
$k+1+j-2i=0$. Similarly as in Corollary 5.3, this implies that
in $S\Cal P(X_n)\otimes \Bbb Z[1/2]$, $\widetilde P_{\rho_k}$ is congruent to
$Det(b_{i,j})_{1\le i,j\le k}$ modulo $\Cal J$, where
$b_{i,j}=\widetilde P_{k+1+j-2i}$ \ if $k+1+j-2i\ne 0$ (
with $\widetilde P_i=0$ for $i<0$)
and $b_{i,j}=1$ \ if $k+1+j-2i=0$.
Thus it is sufficient to prove that
$Det(2b_{i,j})_{1\le i,j\le {n-1}}\equiv e_{n-1}\ldots e_1 \equiv
s_{\rho_{n-1}}$ (mod $\Cal J$).
Recall that $s_{\rho_{n-1}}=Det(c_{i,j})_{1\le i,j\le {n-1}}$ where
$c_{i,j}=e_{n+1+j-2i}$ if $n+1+j-2i\ne 0$
and $c_{i,j}=1$ if $n+1+j-2i=0$, i.e. the matrices $(2b_{i,j})$ and
$(c_{i,j})$ are the same modulo the degree 0 entries.

Let us write the determinants $Det(2b_{i,j})$ and $Det(c_{i,j})$
as the sums of the standard $n!$ terms (some of them are zero).
It is easy to see that apart from the "diagonal" term $e_{n-1}\ldots e_1,$
every other term appearing in both the sums is divisible either by $e_n$
or by
$e_{n-1}e_{n-2}\ldots e_{p+1}e_p^2
$
for some $p\ge 1$. We claim that,
$e_{n-1}e_{n-2}\ldots e_{p+1}e_p^2\in \Cal J.
$
To this end, it suffices to show that $e_{n-1}^2$ belongs to $\Cal J$
and argue as in the proof of Lemma 5.4.
The needed claim follows from the fact that $e_{n-1}^2-e_{n-1}(x_1^2,\ldots,
x_n^2)$ is divisible by $e_n$.

This shows that
$$Det(2b_{i,j})\equiv e_{n-1}\ldots e_1\equiv Det(c_{i,j}) \
(mod \ \Cal J).$$

Thus the lemma is proved.
\qed
\enddemo

The even orthogonal analog of Lemma 5.7 for the operator
$\tk =\partial_{(\bar n,\ldots,\bar 2,\bar 1)}$ if $n$ is even
and $\tk=\partial_{(\bar n,\ldots,\bar 2,1)}$ if $n$ is odd, reads
as follows.

\proclaim {\bf Lemma 5.18} (i) If $f\in S\Cal P(x^2_1,\ldots,x^2_n)
[x_1\cdot \ldots \cdot x_n]$ then
$\tk(f\cdot g)=f\cdot\tk(g)$.

\noindent
(ii) $\tk\bigl(\widetilde P_{\rho_{n-1}}(X_n)\bigr)=(-1)^{n(n-1)/2}$.
\endproclaim
\smallskip

\demo{Proof} (i) This assertion is clear because every polynomial in
$S\Cal P(x^2_1,\ldots,x^2_n)[x_1\cdot\ldots \cdot x_n]$ is $W_n$-invariant.
\smallskip
\noindent
(ii) In this part we will use the Jacobi symmetrizer $\partial$ \
(see the proof of Lemma 5.7). In the following, $\Cal I = \Cal J\Bbb Z[X_n]$.

\smallskip

Let $e_i=e_i(X_n)$.
Since $\widetilde P_{\rho_{n-1}}(X_n)\equiv 2^{-(n-1)} e_{n-1}\ldots e_1$
(mod $\Cal I$), we have
$$\aligned
\tk(\widetilde P_{\rho_{n-1}}(X_n))&=\tk(2^{-(n-1)}e_{n-1}\ldots e_1)=
(\tk\circ\partial)(2^{-(n-1)}x^{\rho_{n-1}}e_{n-1}\ldots e_1) \\
&=\partial_{w_0}(2^{-(n-1)}x^{\rho_{n-1}}e_{n-1}\ldots e_1). \\
\endaligned$$
The degree of the polynomial $x^{\rho_{n-1}}e_{n-1}\ldots e_1$ is $n^2-n$.
Assuming that $\alpha_1+\ldots+\alpha_n=n^2-n$, we have
$\partial_{w_0}(x^{\alpha})\ne 0$ only if
$$
x^{\alpha}=
x^{2n-2}_{\sigma(1)}x^{2n-4}_{\sigma(2)}\ldots x_{\sigma(n-1)}^2
x_{\sigma(n)}^0
$$
for some $\sigma\in S_n$.
Indeed, it follows from Lemma 5.16 that $\partial_{w_0}(x^{\alpha})\ne 0$
only if all the $\alpha_i$'s are even.
Moreover, they must be all different; otherwise $\partial (x^{\alpha})=0$
and consequently $\partial_{w_o}(x^{\alpha})=0$.
We conclude that $\{\alpha_1,\ldots,\alpha_n\}=\{2n-2,2n-4,\ldots,2,0\}$.
But there is only one such a monomial $x^{\alpha}$ in
$x^{\rho_{n-1}}e_{n-1}\ldots e_1$, namely the one with
$(\alpha_1,\ldots,\alpha_n)=(2n-2,2n-4,\ldots,2,0)$.
Therefore
$$
\partial_{w_0}(2^{-(n-1)}x^{\rho_{n-1}}e_{n-1}\ldots e_1)
=2^{-(n-1)}
\partial_{w_0}(x_1^{2n-2}x_2^{2n-4}\ldots x_{n-1}^2x_n^0)
=(-1)^{n(n-1)/2}
$$
by Lemma 5.16.
\qed
\enddemo
The even orthogonal analog of Proposition 5.9 reads (since $e_n(X_n\hak)$ is
a scalar for $\tk$, it suffices to evaluate the images via $\tk$
of $P_I(X_n\hak)$, where all $i_p\le n-1$):

\proclaim{\bf Proposition 5.19}\ Let I be a partition with all parts
not greater than $n-1$.
One has $\triangledown(\widetilde Q_I(X_n\hak)) \ne 0$ only if
the set of parts of
$I$ is equal to
$\{1,2,\ldots,n-1\}$ and each number $p\ \ (1\leqslant p\leqslant n-1)$
appears in $I$ with an odd multiplicity $m_p$.
Then, the following equality holds in $\Bbb Z[X_n]$,
$$
\triangledown\bigl(\widetilde Q_I(X_n\hak)\bigr)
= 2^{n-1}\prod\limits^{n-1}_{p=1}
e_p(x_1^2,\ldots, x_n^2)^{(m_p-1)/2}.
$$
\endproclaim

Arguing as in the proof of Proposition 5.8 one shows that the Gysin maps
$\pi_*$ associated with $\pi: OG_n'V\to X$ (resp. $\pi: OG_n''V\to X$)
are induced by the operator $\tk$.
The role of $LFl(V)$ is played now by the flag bundle parametrizing
flags of rank $1$, rank $2$, ..., rank $n$ isotropic subbundles of $V$
whose rank $n$ subbundle $E$ satisfies $dim(E\cap V_n)_x\equiv n \ (mod \ 2)$
(resp. $dim(E\cap V_n)_x\equiv n+1 \ (mod \ 2)$) for every $x\in X$.
Consequently, the proposition whose proof is the same as the
one of Proposition 5.9, has as its consequence:

\proclaim{\bf Theorem 5.20}
Let $I$ be a partition with all parts not greater than $n-1$.
The element $\widetilde Q_IR\hak$ has a nonzero image
under $\pi_*$ only if each number $p$, $1\le p\le n-1$, appears as a part of
$I$ with an odd multiplicity $m_p$.
If the latter condition holds then
$$
\pi_*\widetilde Q_IR\hak=2^{n-1}
\prod\limits_{p=1}^{n-1}\Bigl((-1)^pc_{2p}V\Bigr)^
{(m_p-1)/2}.
$$
\endproclaim

\proclaim{\bf Theorem 5.21}
The element  $s_IR\hak$ $(l(I)\le n-1)$ has a nonzero image under $\pi_*$
only if the partition $I$ is of the form $2J+\rho_{n-1}$ for some
partition $J$ $(l(J)\le n-1)$.
If $I=2J+\rho_{n-1}$, then
$$
\pi_*s_IR\hak=2^{n-1}s_J^{[2]}V,
$$
where $s_J^{[2]}(-)$ is defined as in Theorem 5.13.
\endproclaim
\demo{Proof}
Since $s_IR\hak=\omega_*(q^{I+\rho_{n-1}})$
where $q=(q_1,\ldots,q_n)$ are the Chern roots of $R\hak$,
we infer from Lemma 5.16 that
$s_IR\hak$ has a nonzero image under $\pi_*$ only if all parts of
$I+\rho_{n-1}$ are even.
This implies that $l(I)=n-1$ and $I$ is strict thus of the form $I'
+\rho_{n-1}$
for some partition $I'$.
Finally all parts of $I'+\rho_{n-1}+\rho_{n-1}$ are even iff $I'=2J$ for some
partition $J$, as required.

Assume now that $I=2J+\rho_{n-1}$ and specialize the identity
from Corollary 5.12
by replacing
the variables $(x_i)$ by the Chern roots $(q_i)$.
The claimed formula now follows since: $s_I(q_1^2,\ldots,q_n^2)$ is
a scalar w.r.t. $\pi_*$, $\pi_*s_{\rho_{n-1}}(q_1,\ldots,q_n)=2^{n-1}$
by Lemma 5.17 and 5.18; moreover
$2(-1)^ic_{2i}V=2e_i(q_1^2,\ldots,q_n^2)$ by Lemma 1.1(2).
\qed
\enddemo

\smallskip
\noindent
{\bf {Remark 5.22.}} \ 1. Our desingularizations of Schubert subschemes are
compositions of flag-- and isotropic Grassmannian bundles
(see Section 1). Therefore Corollary 2.6, the algebra of
$\widetilde Q$-polynomials together with formulas for Gysin push forwards
(Theorem 5.10 for Lagrangian Grassmannians and
a well known formula for projective bundles) give an
explicit algorithm
for calculation the fundamental classes of Schubert subschemes in the
Lagrangian Grassmannian bundles. One has analogous algorithms in the
orthogonal cases. Examples of such calculations are given in Section 6
and 7.

2. In case $X$ is singular, by interpreting polynomials in Chern
classes  as operators acting on Chow groups (see [F]) or singular
homology groups, the same formulas hold (after their obvious adaptation
to the operator setup).
\bigskip

We finish this section with the following important ``orthogonality''
property for the Gysin maps associated with isotropic Grassmannian
bundles.

\smallskip

\proclaim{\bf Theorem 5.23}
(i) For $\pi   :  LG_n V \to X$ and any strict partitions
$I,J \ (\subset \rho_n)$,
$$ \pi_* (\widetilde Q_I R\hak \cdot \widetilde Q_J R \hak) =
\delta_{I,\rho_n \smallsetminus J}.$$
\smallskip

(ii)  For $\pi  :  OG_n V \to X$  ($dimV=2n+1$) and
any strict partitions
$I,J \ (\subset \rho_n)$,
$$ \pi_* (\widetilde P_I R\hak \cdot \widetilde P_J R \hak) =
\delta_{I,\rho_n \smallsetminus J}.$$
\smallskip

(iii) For $\pi   :  OG_n' V \to X$  (resp. $ OG_n'' V \to X$),
 and any strict partitions
$I,J \ (\subset \rho_{n-1})$,
$$ \pi_* (\widetilde P_I R\hak \cdot \widetilde P_J R \hak) =
\delta_{I,\rho_{n-1} \smallsetminus J}.$$

(Here, $\delta_{.,.}$ is the Kronecker delta.)
\endproclaim
\smallskip

\demo{Proof}
We will prove first the Lagrangian case (i).
(In case (ii), the proof goes mutatis mutandis using the divided differences
operator $\partial_{(\nkr{n},\nkr{n-1},\ldots,\nkr{1})}$ for $SO(2n+1)$
instead of the operator $\tk$ for $Sp(2n)$\,.
Case (iii) will be discussed separately at the end of the proof.
\smallskip

Let $X_n=(x_1,...,x_n)$ be a sequence of variables.
We show that the operator $\tk  :  \zz[X_n] \to \zz[X_n]$,
satisfies the following formula for any strict
partitions $I,J \ (\subset \rho_n)$:
$$
\tk\bigl(\widetilde Q_I(X_n\hak) \cdot \widetilde Q_J(X_n\hak)\bigr) =
\delta_{I,\rho_n \smallsetminus J}.
$$
Since $\pi_*$ is induced by $\tk$ (Proposition 5.8), this implies the
assertion. Observe that for the degree reasons
$\tk(\widetilde Q_I \cdot \widetilde Q_J) = 0$
for $|I|+|J| < n(n+1)/2$
(here and in the rest of the proof, \ $\widetilde Q_I=\widetilde
Q_I(X_n\hak)$).
Also, because of the universality of the formula
for $\pi_*$ (see e.g.
Theorem 5.10), we know (Lemma 2.3, Lemma 2.4 and 5.8)
that for  $|I|+|J| = n(n+1)/2$,
$\tk(\widetilde Q_I \cdot \widetilde Q_J) = 0$
unless $J=\rho_n \smallsetminus I$, when
$\tk(\widetilde Q_I \cdot \widetilde Q_J) = 1$.
So it remains to show that for $|I|+|J|>n(n+1)/2$,
$\tk(\widetilde Q_I \cdot \widetilde Q_J) = 0$.
The proof is by double induction whose first parameter is $l(I)$ and the
 second one is $i_l$ where $l=l(I)$
(i.e. the shortest part of $I$).
\medskip

Assume first that $I=(i)$ and use the Pieri-type formula  from
Proposition 4.9.
A general partition $J'$ indexing the R.H.S. of the formula from
Proposition 4.9 stems from $J$ by adding a
horizontal strip of length $i$. Since $|J|+i > n(n+1)/2$, the only
possibility for getting $\tk(\widetilde Q_{J'})\neq 0$ is the following
(Theorem 5.10): there exist two equal parts $p$ in $J'$ such that
after factoring out
$\widetilde Q_{p,p}$ from $\widetilde Q_{J'}$ (Proposition 4.3)
we obtain $\widetilde Q_{\rho_n}$  (recall that $\widetilde Q_{p,p}$
is a scalar w.r.t. $\tk$). But $l(J') \leq l(J)+1 \leq n+1$,
so after factoring out
the length of the so-obtained partition is not greater than $n-1$, i.e.
this partition is not $\rho_n$.
\bigskip

To perform the induction step write $I'= (i_1,...,i_{l-1})$ and $r=i_l$ where
we assume that $l=l(I) \geq 2$. Using the Pieri-type formula again, we have:
$$\widetilde Q_{I} \cdot \widetilde Q_{J}
= (\widetilde Q_{I'} \cdot \widetilde Q_{r}) \cdot \widetilde Q_{J}-
(\sum_M 2^{m(I',r;M)}\widetilde Q_{M}) \cdot \widetilde Q_{J}
= \widetilde Q_{I'} \cdot (\widetilde Q_{J} \cdot \widetilde Q_{r})-
(\sum_M 2^{m(I',r;M)}\widetilde Q_{M}) \cdot \widetilde Q_{J}$$
$$= \widetilde Q_{I'} \cdot (\sum_N 2^{m(J,r;N)}\widetilde Q_{N})-
(\sum_M 2^{m(I',r;M)}\widetilde Q_{M}) \cdot \widetilde Q_{J}.$$

\noindent
Here $M$ runs over all partitions different from $I$ which contain $I'$ with
$M/I'$ being a horizontal strip of length $r$. Observe that either $l(M)<l(I)$
or
$l(M)=l(I)$ but $m_l<i_l=r$, so we can apply the induction assumption
to $M$. The partitions $M$ and $N$ can have equal parts; if so,
using the factorization property, we write:
\medskip
\centerline
{$\widetilde Q_{M} = \widetilde Q_{p_1,p_1} \cdot\ldots\cdot \widetilde
Q_{p_s,p_s}
\cdot \widetilde Q_{M_1}$ \ \ \
and \ \ \  $\widetilde Q_{N} = \widetilde Q_{q_1,q_1} \cdot \ldots \cdot
\widetilde Q_{q_t,q_t}
\cdot \widetilde Q_{N_1},$}
\medskip
\noindent
where $M_1 ,N_1$ are strict partitions
and $p_1>\ldots>p_s$, $q_1>\ldots>q_t$ are positive integers.
Using the induction assumption or because of the degree
reasons we see that the only possibility to get in the first
sum a summand (corresponding to $N$) which
is not anihilated by $\tk$ is: after adding to $J$ a horizontal
strip of length $r$
and factoring out all pairs of equal rows, we obtain
the partition $N_1 = \rho_n \smallsetminus I'$.
Similarly, the only possibility to get
in the second sum a summand
(corresponding to $M$) which is not anihilated by $\tk$ is:
after adding to $I'$ a horizontal strip of length $r$
and factoring out all pairs of equal rows, we obtain the partition
$M_1 = \rho_n\smallsetminus J$.
\medskip

Therefore to conclude the proof it is sufficient to define, for a fixed
pair of strict partitions $I', J$ and fixed positive integers
$r$ and $p. : p_1>\ldots >p_s$,  a bijection between the sets of partitions
(with parts not exceeding $n$):
\smallskip

{\parindent=35pt
\item{$\Cal N = \bigl\{$} $N\ |\ N\supset J;\ N/J$ is a horizontal strip of
length $r$;
$N$ has exactly $s$ parts occuring twice, equal to $p.$; after subtraction
from $N$ the parts $p.$ one obtains $\rho_n \smallsetminus I' \bigr\}$
\item{and}
\smallskip
\item{$\Cal M = \bigl\{$} $M\ |\ M\supset I';M/I'$ is a horizontal strip of
length $r$;
$M$ has exactly $s$ parts occuring twice, equal to $p.$; after subtraction
from $M$ the parts $p.$ one obtains $\rho_n \smallsetminus J \bigr\}$
\par}
\smallskip\noindent
which preserves the cardinality of the connected components of the strip,
not meeting the first column (compare the Pieri-type formula used).
\medskip

In order to define the bijection $\Phi:\Cal N\to \Cal M$ we first invoke
the diagramatic presentation of
the $\rho_n$-complementary partition from [P2, p.160]: for example $n=9$,
$I=(9,6,3,2)$, $\rho_9\smallsetminus I=(8,7,5,4,1)$,
\medskip

\vbox{\settabs 30 \columns
\+ &&&&&&&&&& &&  & & & & & & & & $\bullet$ &&&&&&&&& \cr
\+ &&&&&&&&&& &&  & & & & & & &$\bullet$ & $\bullet$ &&&&&&&&& \cr
\+ &&&&&&&&&& &&  & & & & & &$\bullet$ &$\bullet$ & $\bullet$ &&&&&&&&& \cr
\+ &&&&&&&&&& &&  & & & & &$\bullet$ &$\bullet$ &$\bullet$ & $\bullet$
&&&&&&&&& \cr
\+ &&&&&&&&&& &&  & & & &$\circ$ &$\bullet$ &$\bullet$ &$\bullet$ & $\bullet$
&&&&&&&&& \cr
\+ &&&&&&&&&& &&  & & &$\circ$ &$\circ$ &$\circ$ &$\circ$ &$\bullet$ &
$\bullet$ &&&&&&&&& \cr
\+ &&&&&&& Fig.1 &&& &&  & &$\circ$ &$\circ$ &$\circ$ &$\circ$ &$\circ$
&$\bullet$ & $\bullet$ &&&&&&&&& \cr
\+ &&&&&&&&&& &&  &$\circ$ &$\circ$ &$\circ$ &$\circ$ &$\circ$ &$\circ$
&$\circ$ & $\bullet$ &&&&&&&&& \cr
\+ &&&&&&&&&& &&$\circ$ &$\circ$ &$\circ$ &$\circ$ &$\circ$ &$\circ$ &$\circ$
&$\circ$ & $\bullet$ &&&&&&&&& \cr
}
\medskip
\noindent
(the collection of ``$\bullet$'' gives the shifted diagram of $I$
(appropriately placed); the collection of ``$\circ$'' gives the shifted
diagram of $\rho_9\smallsetminus I$).
The map $\Phi:\Cal N\to \Cal M$ is defined as follows.
Having an element $N\in \Cal N$, i.e. a strict partition $J$ with an
added horizontal strip of length $r$, e.g.
$J=(9,6,3,2), r=5, N=(9,8,3,3,2), s=1, p.: 3$ \ ( and $I'=(7,6,5,4,3,1)$ ):
\bigskip

\vbox{\settabs 30 \columns
\+ &&&&&&&&&& &&$\circledast$&$\circledast$& & & & & & &  &&&&&&&&& \cr
\+ &&&&&&&&&& &&$\bullet$&$\bullet$&$\circledast$ & & & & & &  &&&&&&&&& \cr
\+ &&&&&&&&&& &&$\bullet$&$\bullet$&$\bullet$& & & & & & &&&&&&&&& \cr
\+ &&&&&&& Fig.2 &&&
&&$\bullet$&$\bullet$&$\bullet$&$\bullet$&$\bullet$&$\bullet$&$\circledast$
&$\circledast$ &  &&&&&&&&& \cr
\+ &&&&&&&&&&
&&$\bullet$&$\bullet$&$\bullet$&$\bullet$&$\bullet$&$\bullet$&$\bullet$
&$\bullet$ &$\bullet$  &&&&&&&&& \cr
}
\bigskip
\noindent
we remove the $s$ bottom rows in all pairs of equal rows (in the example,
the third row) and place the shift of the so-obtained diagram as in Fig.1
to get the diagram $\widehat N$, say.
In our example we get the diagram in Fig.3 \ :
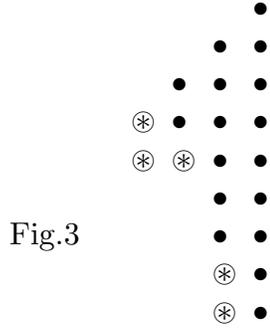
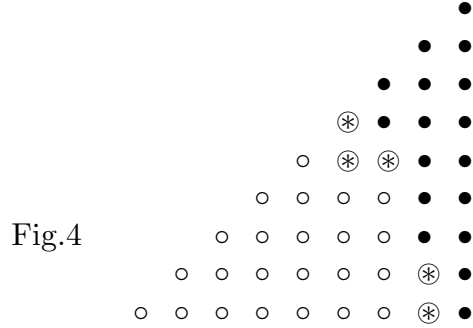
\bigskip

\vbox{\settabs 30 \columns
\+ &&&&& & & & & $\bullet$ &&&&&&&& &  & & & & & & & & $\bullet$ &&& \cr
\+ &&&&& & & &$\bullet$ & $\bullet$ &&&&&&&& &  & & & & & & &$\bullet$ &
$\bullet$ &&&\cr
\+ &&&&& & &$\bullet$ &$\bullet$ & $\bullet$ &&&&&&&& &  & & & & & &$\bullet$
&$\bullet$ & $\bullet$ &&& \cr
\+ &&&&& &$\circledast$&$\bullet$ &$\bullet$ & $\bullet$ &&&&&&&& &  & & & &
&$\circledast$ &$\bullet$ &$\bullet$ & $\bullet$ &&& \cr
\+ &&&&&  &$\circledast$ &$\circledast$ &$\bullet$ & $\bullet$ &&&&&&&& &  & &
& &$\circ$ &$\circledast$ &$\circledast$ &$\bullet$ & $\bullet$ &&& \cr
\+ &&&&&  & & &$\bullet$ & $\bullet$ &&&&&&&& &  & & &$\circ$ &$\circ$ &$\circ$
&$\circ$ &$\bullet$ & $\bullet$ &&& \cr
\+ &&& Fig.3 &&  & & &$\bullet$ & $\bullet$ &&&&&& Fig.4 && &  & &$\circ$
&$\circ$ &$\circ$ &$\circ$ &$\circ$ &$\bullet$ & $\bullet$ &&& \cr
\+ &&&&&  & & &$\circledast$ & $\bullet$ &&&&&&&& &  &$\circ$ &$\circ$ &$\circ$
&$\circ$ &$\circ$ &$\circ$ &$\circledast$ &$\bullet$ &&& \cr
\+ &&&&&  & & &$\circledast$& $\bullet$ &&&&&&&& &$\circ$ &$\circ$ &$\circ$
&$\circ$ &$\circ$ &$\circ$ &$\circ$ &$\circledast$ & $\bullet$ &&& \cr
}

\bigskip
\noindent
(We know, by the definition of $\Cal N$,
that if we would also remove from $\widehat N$ the remaining parts of lengths
$p.$
then the resulting partition will be $\rho_n\smallsetminus I'$. We preserve
these parts, however, because we need them for the construction of $\Phi(N)$.)
Then we construct the complement of the so-obtained diagram in $\rho_n$.
In our example, using ``$\circ$'' to visualize
the complementary diagram we get the diagram in Fig.4.
By reshifting the so-obtained complementary diagram plus {\it the same}
horizontal
strip (now added to this complementary diagram) - call it ${\Phi(N)}_0$,
and inserting $s$ rows
of lengths
$p.$, we get the needed partition $\Phi(N)$.
Observe that :

\smallskip
\noindent
1) Since at the last stage we have inserted rows of lengths
$p.$, $\Phi(N)$ consists of the diagram $I'$ with an added horizontal strip
of length $r$.
\smallskip
\noindent
2) $\Phi(N)$ has exactly $s$ parts occuring twice, equal
to $p.$ (apart from the parts inserted at the last stage, the remaining $s$
parts are the rows whose the rightmost boxes are precisely the lowest boxes
of the rows of length $p.$ in $\widehat N$).
\smallskip
\noindent
3) After removing from $\Phi(N)$
the $2s$ parts equal to $p.$, we get $\rho_n\smallsetminus J$ (this is the same
as removing from ${\Phi(N)}_0$ the $s$ parts equal to $p.$ - but
${\Phi(N)}_0$ minus $s$ parts equal to $p.$ complements precisely $J$
in $\rho_n$).
\smallskip
\noindent
Therefore $\Phi(N) \in \Cal M$. Also, the cardinality of the
connected components of the strip not meeting the first column is
preserved by $\Phi$.
In our example, we obtain
\bigskip

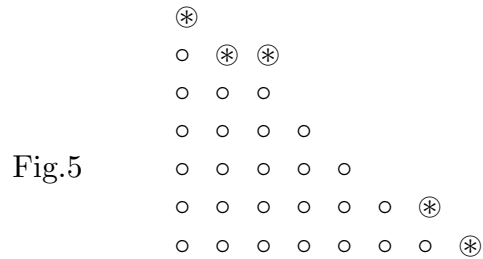
\vbox{\settabs 30 \columns
\+ &&&&&&&&&& &$\circledast$& & & & & & & &   &&&&&&&&&& \cr
\+ &&&&&&&&&& &$\circ$&$\circledast$ &$\circledast$ & & & & & &   &&&&&&&&&&
\cr
\+ &&&&&&&&&& &$\circ$&$\circ$ &$\circ$ & & & & & &   &&&&&&&&&& \cr
\+ &&&&&&&&&& &$\circ$&$\circ$ &$\circ$ &$\circ$ & & & & &   &&&&&&&&&& \cr
\+ &&&&&&& Fig.5 &&& &$\circ$&$\circ$ &$\circ$ &$\circ$ &$\circ$ & & & &
&&&&&&&&&& \cr
\+ &&&&&&&&&& &$\circ$&$\circ$ &$\circ$ &$\circ$ &$\circ$ &$\circ$
&$\circledast$ & &   &&&&&&&&&& \cr
\+ &&&&&&&&&& &$\circ$ &$\circ$ &$\circ$ &$\circ$ &$\circ$ &$\circ$ &$\circ$
&$\circledast$ &   &&&&&&&&&& \cr
}
\medskip
\noindent
i.e. $\Phi(N)=(8,7,5,4,3,3,1)$.
\smallskip

Let us now define, by reversing the roles of $J$ and $I'$,
the map $\Psi :\Cal M \to \Cal N$.
If we define, by a complete analogy to the above, the partitions
$\widehat M$ and $\Psi(M)_0$, then we have $\widehat N=\Psi(M)_0$ and
$\Phi(N)_0=\widehat M$; and clearly
$\Psi \circ \Phi = id_{\Cal N}$ and $\Phi \circ \Psi = id_{\Cal M}$.
\smallskip

This proves the orthogonality theorem in the Lagrangian case.
\smallskip

Essentially the same proof, with $\widetilde Q$-polynomials replaced
by $\widetilde P$-polynomials (for which a Pieri-type formula is
given below), works in the odd orthogonal case (ii).

\smallskip

In the even orthogonal case the proof goes as follows. Let $\tk :\Bbb Z[X_n]
\to \Bbb Z[X_n]$ be the even orthogonal divided differences operator
inducing $\pi_*$.

We show that the operator $\tk  :  \zz[X_n] \to \zz[X_n]$,
satisfies the following formula for any strict
partitions $I,J \ (\subset \rho_{n-1})$:
$$
\tk\bigl(\widetilde P_I(X_n\hak) \cdot \widetilde P_J(X_n\hak)\bigr) =
\delta_{I,\rho_{n-1} \smallsetminus J}.
$$
Since $\pi_*$ is induced by $\tk$, this implies the
assertion. Observe that for the degree reasons
$\tk(\widetilde P_I \cdot \widetilde P_J) = 0$
for $|I|+|J| < n(n-1)/2$
(here and in the rest of the proof, \ $\widetilde P_I=\widetilde
P_I(X_n\hak)$).
Also, because of the universality of the formula
for $\pi_*$
that for  $|I|+|J| = n(n-1)/2$,
$\tk(\widetilde P_I \cdot \widetilde P_J) = 0$
unless $J=\rho_{n-1} \smallsetminus I$, when
$\tk(\widetilde P_I \cdot \widetilde P_J) = 1$.
So it remains to show that for $|I|+|J|>n(n-1)/2$,
$\tk(\widetilde P_I \cdot \widetilde P_J) = 0$.
The proof is by double induction whose first parameter is $l(I)$ and the
 second one is $i_l$ where $l=l(I)$
(i.e. the shortest part of $I$).
\medskip

Assume first that $I=(i)$ and use the Pieri-type formula for
$\widetilde P$-polynomials (see Proposition 4.9;
a Pieri-type formula for $\widetilde P$-polynomials
reads similarly:
$$
\widetilde P_J\cdot \widetilde P_i
= \sum 2^{m'(J,i;J')}\widetilde P_{J'},
$$
the only difference being the exponent $m'(J,i;J')$
equal to $m(J,i;J')$ if $J'/J$ meets the first column and $m(J,i;J')-1$
 -- if not.)

A general partition $J'$ indexing the R.H.S. of the Pieri formula
stems from $J$ by adding a
horizontal strip of length $i$. Since $|J|+i > n(n-1)/2$, the only
possibility for getting $\tk(\widetilde P_{J'})\neq 0$ is the following:
there exists a (single) row of length $n$ or
there exist two equal parts $p$ in $J'$ ($1\le p \le n-1$) such that
after factoring out $\widetilde P_n$ and
$\widetilde P_{p,p}$ from $\widetilde P_{J'}$ (Proposition 4.3)
we obtain $\widetilde P_{\rho_{n-1}}$.
As in the proof of the Lagrangian case we see that it is
impossible to get $\widetilde P_{\rho_{n-1}}$ after factoring out
$\widetilde P_{p,p}$. Also, it is impossible to get
$\widetilde P_{\rho_{n-1}}$ by factoring out $\widetilde P_n$.
Indeed, using the Pieri-type formula, we should add
one box to each of the first $n$ columns which is impossible
because $i\le n-1$.
\medskip

To perform the induction step write $I'= (i_1,...,i_{l-1})$ and $r=i_l$ where
we assume that $l=l(I) \geq 2$. Using the Pieri-type formula, we have:
$$\widetilde P_{I} \cdot \widetilde P_{J}
= \widetilde P_{I'} \cdot (\sum_N 2^{m'(J,r;N)}\widetilde P_{N})-
(\sum_M 2^{m'(I',r;M)}\widetilde P_{M}) \cdot \widetilde P_{J}.$$

\noindent
Here $M$ runs over all partitions different from $I$ which contain $I'$ with
$M/I'$ being a horizontal strip of length $r$. Observe that either $l(M)<l(I)$
or
$l(M)=l(I)$ but $m_l<i_l=r$, so we can apply the induction assumption
to $M$. The partitions $M$ and $N$ can have equal parts; if so,
using the factorization property, we write:
\medskip
\centerline
{$\widetilde P_{M} = \widetilde P_{p_1,p_1} \cdot\ldots\cdot \widetilde
P_{p_s,p_s}
\cdot \widetilde P_{M_1}$ \ \ \
and \ \ \  $\widetilde P_{N} = \widetilde P_{q_1,q_1} \cdot \ldots \cdot
\widetilde P_{q_t,q_t}
\cdot \widetilde P_{N_1},$}
\medskip
\noindent
where $M_1 ,N_1$ are strict partitions
and $p_1>\ldots>p_s$, $q_1>\ldots>q_t$ are positive integers.
Moreover $M$ and $N$ can contain a single row of length $n$, and if so,
then the polynomial $\widetilde P_n$ can be factored out by a property
of the operator $\tk$.
Using the induction assumption or because of the degree
reasons we see that the only possibility to get in the first
sum a summand (corresponding to $N$) which
is not anihilated by $\tk$ is: after adding to $J$ a horizontal
strip of length $r$
and factoring out all pairs of equal rows and the row of length $n$ (if any),
we obtain the partition $N_1 = \rho_{n-1} \smallsetminus I'$.
Similarly, the only possibility to get
in the second sum a summand
(corresponding to $M$) which is not anihilated by $\tk$ is:
after adding to $I'$ a horizontal strip of length $r$
and factoring out all pairs of equal rows and the row of length $n$,
if any, we obtain the partition
$M_1 = \rho_{n-1}\smallsetminus J$.
\medskip

Therefore to conclude the proof it is sufficient to give two
bijections.

The data of the first bijection are: a pair of strict partitions
$I', J\subset \rho_{n-1}$ and fixed positive integers  $r$
and $p. : p_1>\ldots >p_s$. One needs a bijection between
the sets of partitions (with parts
not exceeding $n-1$):
\smallskip

{\parindent=35pt
\item{$\Cal N = \bigl\{$} $N\ |\ N\supset J;\ N/J$ is a horizontal strip of
length $r$;
$N$ has exactly $s$ parts occuring twice, equal to $p.$; after subtraction
from $N$ the parts $p.$ one obtains $\rho_{n-1} \smallsetminus I' \bigr\}$
\item{and}
\smallskip
\item{$\Cal M = \bigl\{$} $M\ |\ M\supset I';M/I'$ is a horizontal strip of
length $r$;
$M$ has exactly $s$ parts occuring twice, equal to $p.$; after subtraction
from $M$ the parts $p.$ one obtains $\rho_{n-1} \smallsetminus J \bigr\}$
\par}
\smallskip\noindent
which preserves the property that the strip meets or not the first column
and preserves the cardinality of the connected components of the strip,
not meeting the first column -- compare the Pieri-type formula used).

Here the bijection $\Phi: \Cal N\to \Cal M$ from the proof
of the Lagrangian case with $n$
replaced by $n-1$ does the job
(i.e. to construct $\Phi(N)_0$ for $N\in \Cal N$
we take the complement in $\rho_{n-1}$). Note that $\Phi$ preserves
the property that the strip meets or not the first column by the
construction.
\smallskip

The data of the second bijection are also: a pair of strict partitions
$I', J\subset \rho_{n-1}$
and fixed positive integers  $r$ and $p. : p_1>\ldots
>p_s$. One needs a bijection between the sets of partitions (with parts
not exceeding $n$):
\smallskip

{\parindent=35pt
\item{$\Cal N' = \bigl\{$} $N\ |\ N\supset J;\ N/J$ is a horizontal strip of
length $r$; $N$ has a single part equal to $n$;
$N$ has exactly $s$ parts occuring twice, equal to $p.$; after subtraction
from $N$ the parts $p.$ and $n$ one obtains
$\rho_{n-1} \smallsetminus I' \bigr\}$
\item{and}
\smallskip
\item{$\Cal M' = \bigl\{$} $M\ |\ M\supset I';M/I'$ is a horizontal strip of
length $r$; $M$ has a single part equal to $n$;
$M$ has exactly $s$ parts occuring twice, equal to $p.$; after subtraction
from $M$ the parts $p.$ and $n$ one obtains
$\rho_{n-1} \smallsetminus J \bigr\}$
\par}
\smallskip\noindent
which preserves the property that the strip meets or not the first
column and preserves cardinality of the connected components of the strip,
not meeting the first column.

Here we also use the map $\Phi$ from the proof
of case (i)
(to construct $\Phi(N)_0$ for $N\in \Cal N'$
we take the complement in $\rho_n$).
We ilustrate the map $\Phi$ on the following example.

Let $n=10$\,, $J=(8,7,4,2)$\,, $N=(10,8,4,4,2)$ and $I'=(9,7,6,5,4,3,1)$\,.
\bigskip
\vbox{\settabs 30 \columns
\+ &&&& && &&&&&&&& &&&      &&&&&&$\bullet$ &&&&&& \cr
\+ &&&& && &&&&&&&& &&&      &&&&&$\bullet$&$\bullet$ &&&&&& \cr
\+ &&&& &&$\circledast$&$\circledast$& & & & & & &  &&&
&&&&$\bullet$&$\bullet$&$\bullet$ &&&&&&\cr
\+ &&&& &&$\bullet$&$\bullet$&$\circledast$&$\circledast$ & & & & &  &&&
&&&$\circledast$&$\bullet$&$\bullet$&$\bullet$ &&&&&&\cr
\+ &&&& $N=$&&$\bullet$&$\bullet$&$\bullet$&$\bullet$& & & & & &&&
&$\widehat N=$&&$\circledast$&$\circledast$&$\bullet$&$\bullet$ &&&&&& \cr
\+ &&&& &&$\bullet$&$\bullet$&$\bullet$&$\bullet$&$\bullet$&$\bullet$&$\bullet$
&$\circledast$ &  &&&
&&&& $\circledast$&$\bullet$&$\bullet$ &&&&&& \cr
\+ &&&& &&$\bullet$&$\bullet$&$\bullet$&$\bullet$&$\bullet$&$\bullet$&$\bullet$
&$\bullet$ &$\circledast$&$\circledast$&&
&&&&&$\bullet$&$\bullet$ &&&&&& \cr
\+ &&&& && &&&&&&&& &&&       &&&&&$\bullet$&$\bullet$ &&&&&& \cr
\+ &&&& && &&&&&&&& &&&       &&&&&$\circledast$&$\circledast$ &&&&&& \cr
\+ &&&& && &&&&&&&& &&&       &&&&&&$\circledast$ &&&&&& \cr
}
\smallskip
\vbox{\settabs 30 \columns
\+ && &&&&& &&&&& &&&&& & $\bullet$ &&& &&&&& &&& \cr
\+ && &&&&& &&&&& &&&&& $\bullet$ & $\bullet$ &&& &&&&& &&&\cr
\+ && &&&&& &&&&& &&&& $\bullet$ &$\bullet$ & $\bullet$ &&& &&&&& &&& \cr
\+ && &&&&& &&&&& &&& $\circledast$ &$\bullet$ &$\bullet$ & $\bullet$ &&& &&&&&
&&& \cr
\+ && &&&&& &&&&& && $\circ$ &$\circledast$ &$\circledast$ &$\bullet$ &
$\bullet$ &&& &&&&& &&&\cr
\+ && &&&&& &&&&& & $\circ$ &$\circ$ &$\circ$ &$\circledast$ &$\bullet$ &
$\bullet$ &&& &&&&& &&&\cr
\+ && &&&&& &&&&&  $\circ$ &$\circ$ &$\circ$ &$\circ$ &$\circ$ &$\bullet$ &
$\bullet$ &&& &&&&& &&& \cr
\+ && &&&&& &&&&  $\circ$ &$\circ$ &$\circ$ &$\circ$ &$\circ$ &$\circ$
&$\bullet$ &$\bullet$ &&& &&&&& &&& \cr
\+ && &&&&& &&&  $\circ$ &$\circ$ &$\circ$ &$\circ$ &$\circ$ &$\circ$ &$\circ$
&$\circledast$ & $\circledast$ &&& &&&&& &&& \cr
\+ && &&&&& &&  $\circ$ &$\circ$ &$\circ$ &$\circ$ &$\circ$ &$\circ$ &$\circ$
&$\circ$ &$\circ$ & $\circledast$ &&& &&&&& &&& \cr
}
\bigskip\bigskip\medskip
\vbox{\settabs 30 \columns
\+ &&& &&&&&& $\circledast$ &&&&& &&&&& $\circledast$ &&&&&&&&& & \cr
\+ &&& &&&&& $\circ$&$\circledast$&$\circledast$ &&&&
&&&&& $\circ$&$\circledast$&$\circledast$ &&&&&&& & \cr
\+ &&& &&&& $\circ$&$\circ$&$\circ$&$\circledast$ &&&&
&&&&& $\circ$&$\circ$&$\circ$&$\circledast$ &&&&&& & \cr
\+ & $\phi(N)_0=$ &&&&&$\circ$&$\circ$&$\circ$&$\circ$&$\circ$ &&&&
&$\phi(N)=$ &&&& $\circ$&$\circ$&$\circ$&$\circ$ &&&&&& & \cr
\+ &&& && $\circ$&$\circ$&$\circ$&$\circ$&$\circ$&$\circ$ &&&&
&&&&& $\circ$&$\circ$&$\circ$&$\circ$&$\circ$& &&&&& & \cr
\+ &&& &
$\circ$&$\circ$&$\circ$&$\circ$&$\circ$&$\circ$&$\circ$&$\circledast$&$\circledast$ &&
&&&&& $\circ$&$\circ$&$\circ$&$\circ$&$\circ$&$\circ$ &&&& & \cr
\+ &&&
$\circ$&$\circ$&$\circ$&$\circ$&$\circ$&$\circ$&$\circ$&$\circ$&$\circ$&$\circledast$ &&
&&&&&
$\circ$&$\circ$&$\circ$&$\circ$&$\circ$&$\circ$&$\circ$&$\circledast$&$\circledast$ & & \cr
\+ &&& &&&&&& &&&&& &&&&&
$\circ$&$\circ$&$\circ$&$\circ$&$\circ$&$\circ$&$\circ$&$\circ$&$\circ$&$\circledast$ & \cr
}
\bigskip\smallskip\noindent
$(N=\{10,8,4,4,2\})\smallsetminus\{10,4,4\}=\{9,8,7,6,5,4,3,2,1\}
\smallsetminus (\{9,7,6,5,4,3,1\}=I')$
\smallskip\noindent
$(\Phi(N)=\{10,9,6,5,4,4,3,1\})\smallsetminus\{10,4,4\}
=\{9,8,7,6,5,4,3,2,1\}\smallsetminus (\{8,7,4,2\}=J)$

\smallskip
We have $\Phi(\Cal N')\subset \Cal M'$. Indeed, if $N\in \Cal N'$ then
$\Phi(N)$ has a part equal to $n$ by the construction. Moreover,
the equations:
$$
N\smallsetminus \{n,p.\} = \rho_{n-1} \smallsetminus I' \ ,
 \ \Phi(N)\smallsetminus \{n,p.\} = \rho_{n-1} \smallsetminus J
$$
are equivalent to the equations:
$$
N\smallsetminus \{p.\} = \rho_{n} \smallsetminus I' \ ,
 \ \Phi(N)\setminus \{p.\} = \rho_{n} \smallsetminus J,
$$
so the assertion follows from the proof in case (i) above.
\smallskip

By reversing the roles of $J$ and $I'$, one defines (as in the proof of
case (i))
the map $\Psi :\Cal M' \to \Cal N'$ which satisfies:
$\Psi \circ \Phi = id_{\Cal N'}$ and $\Phi \circ \Psi = id_{\Cal M'}$.
\smallskip

This ends the proof of the theorem.
\qed
\enddemo

\bigskip\medskip

\centerline{\bf 6. Single Schubert condition}

\bigskip

We consider first the Lagrangian case $G={L}G_nV$ and
follow the notation introduced in Section 1.

\proclaim{\bf Proposition 6.1}
The class of $\Omega(a)$ in $A^*(G)$, where $a=n+1-i$, is given by the formula
$$
\bigl[\Omega(a)\bigr] = \sum\limits^i_{p=0} c_pR\hak\cdot
s_{i-p}(V_{a}\hak).
$$
\endproclaim

\demo{Proof}
The desingularization $\Cal F$ of $\Omega(a)\subset G$ is given
by the composition (recall that $Fl(a.)$ from Section 1 is here $\Bbb P(V_a)$
and $C$ is the tautological line bundle on it):
$$
\Cal F=LG_{n-1}(C^{\perp}/C) @>\pi_1>> \Bbb P(V_{a}) @>\pi_2>> G\ ,
$$
where $\pi_1$ and $\pi_2$ denote the corresponding projection maps.
By Corollary 2.6 we have
$$
[Z] = \sum\limits_{strict\ I\subset\rho_n}
\widetilde Q_ID\hak \cdot \widetilde Q_{\rho_n\smallsetminus I}R\hak .
\tag *
$$
Let $S$ be the tautological rank $n-1$ bundle on $\Cal F$; \ $S=D/C_{\Cal F}$.
Let $c=c_1(C\hak)$. Then, by Proposition 4.1,
$$
\widetilde Q_ID\hak = \sum\limits^n_{k=0} (\pi_1^*c)^k\cdot \sum\limits_J
\widetilde Q_JS\hak\ ,
\tag **
$$
the sum over all partitions $J\subset I$
of weight $|J|=|I|-k$ and $I/J$ has at most
one box in each row.
By Theorem 5.10 the only $I$'s in (*)
for which $(\pi_1)_*\widetilde Q_ID\hak \ne 0$,
are those containing $\rho_{n-1}$, i.e. $I$ must be equal to one of the
partitions of the form
$I_p = \bigl(n,n-1,\ldots,p+1,p-1,\ldots,1\bigr)$
for some $p=0,1,\ldots,n$.
For $I=I_p$ the only term in (**) which contributes
after applying $(\pi_1)_*$
is the one
with $J=\rho_{n-1}$ and $k=n-p$.

Since, by a well-known push forward formula for projective bundles,
we have
$$
(\pi_2)_*\bigl(c^{n-p}\bigr) = s_{n-p-(n-i)}\bigl(V_{a}\hak\bigr) =
s_{i-p}\bigl(V_{a}\hak\bigr)\ ,
$$
we infer that only $p=0,1,\ldots,i$ \ give a nontrivial contribution from (**)
(with $k=n-p$).
Finally, we get
$$
\bigl[\Omega(a)\bigr] = (\pi_2\pi_1)_* [Z] =
\sum\limits^i_{p=0} \widetilde Q_pR\hak \cdot s_{i-p}(V_{a}\hak)=
\sum\limits^i_{p=0} c_pR\hak \cdot s_{i-p}(V_{a}\hak)
$$
as asserted.
\qed
\enddemo

\bigskip

Essentially the same computation gives the following formula
in for $G=OG_nV$ where $dim\,V=2n+1$\,.
\smallskip

\proclaim{\bf Proposition 6.2}
The class of $\Omega(a)$ in $A^*(G)$\,, where $a=n+1-i$ is given by the formula
\medskip

\Odst\Odst\Odst
$\bigl[\Omega(a)\bigr]=1/2\cdot\sum\limits^i_{p=0}c_pR\hak\,
s_{i-p}(V_a\hak)\,.$\footnote{\eightrm Observe that though "1/2" appears
in the formula, the integrality property of the class obtained holds true
(i.e. we get the class in the Chow group with the integer coefficients).
This follows directly from our way of computing it.
Indeed, the (odd) Orthogonal version of Proposition 2.5 and consequently
also of Corollary 2.6 holds true over integers.
Also, the integrality is preserved by the Gysin maps in the odd
Orthogonal analogs of Propositions 3.1 and 3.2.
The same remark applies to Proposition 6.3 and Theorem 7.9 below.}
\endproclaim
\bigskip

Consider finally the even orthogonal case where the computation is slightly
different.

\proclaim{\bf Proposition 6.3}
The class of $\Omega(a)$ in $A^*(OG'_nV)$ for odd $n$, or
in $A^*(OG''_nV)$ for even $n$, is given by the following expression
where $i=n-a$:
$$
\bigl[\Omega(a)\bigr]=1/2\cdot \sum\limits^i_{p=0}\bigl(c_pR\hak+c_pV_n
\bigr)\,s_{i-p}(V_a\hak).
$$
\endproclaim
\demo{\bf Proof}
Suppose that $n$ is odd.
Then the desingularization $\Cal F'$ of $\Omega(a)\subset G=OG_n'V$ is given
by the composition:
$$
\Cal F'=OG_{n-1}'(C^{\perp}/C) @>\pi_1>> \Bbb P(V_{a}) @>\pi_2>> G\ ,
$$
where $\pi_1$ and $\pi_2$ denote the corresponding projection maps.
\smallskip

If $n$ is even then
the desingularization $\Cal F''$ of $\Omega(a)\subset G=OG_n''V$ is given
by the composition:
$$
\Cal F''=OG_{n-1}''(C^{\perp}/C) @>\pi_1>> \Bbb P(V_{a}) @>\pi_2>> G\ ,
$$
where $\pi_1$ and $\pi_2$ denote the corresponding projections. In the
following we denote by $\Cal F$ both $\Cal F'$ and $\Cal F''$ for brevity.
\medskip

By Proposition 2.7 we have in both the cases:
$$
[Z] = \sum\limits_{strict\ I\subset\rho_{n-1}}
\widetilde P_ID\hak \cdot \widetilde P_{\rho_{n-1}\smallsetminus I}R\hak .
$$
Let $S$ be the tautological rank $n-1$ bundle on $\Cal F$; \ $S=D/C_{\Cal F}$.
Let $c=c_1C\hak$. Then, by Proposition 4.1 interpreted
now in terms of $\widetilde P$-polynomials we have
$$
\widetilde P_ID\hak = \sum\limits^n_{k=0} (\pi_1^*c)^k\cdot \Bigl(\sum\limits_J
2^{l(J)-l(I)}
\widetilde P_JS\hak \Bigr)\ ,
\tag***
$$
the sum over all partitions $J\subset I$
of weight $|J|=|I|-k$ and $I/J$ has at most
one box in each row.
By Theorem 5.20, if
$(\pi_1)_*\widetilde Q_ID\hak \ne 0$
then $I\supset \rho_{n-2}$, so $I$ must be equal to the partition
$I_p = \bigl(n-1,n-2,\ldots,p+1,p-1,\ldots,1\bigr)$
for some $p=0,1,\ldots,n-1$.
More precisely, the only terms in (***) which contribute nontrivially
after applying $(\pi_1)_*$ correspond to the following two instances:
\medskip
\noindent
(1) \ $I=I_0=\rho_{n-1}$, $k=0$ and $J=\rho_{n-1}$ -- this gives
a difference between the odd
orthogonal case and the present one.
\medskip
\noindent
(2) \ $I=I_p$, $k=n-p-1$ and $J=\rho_{n-2}$; here $p=0,1,\ldots, n-1$ but
we will see that only $p=0,1,\ldots,i$  \ give a nontrivial contribution.
\medskip

Let us first compute the contribution of (1).
We claim that
$$
(\pi_1)_*P_{\rho_{n-1}}S\hak = 1/2\cdot v,
$$
where $v=c_{n-1}(V_n/C)$.
Indeed, if $n$ is odd then
$$(\pi_1)_*P_{\rho_{n-1}}S\hak = 1/2\cdot (\pi_1)_*\bigl(c_{n-1}S \cdot
P_{\rho_{n-2}}S\hak \bigr)
=1/2\cdot (\pi_1)_*\bigl(v\cdot
P_{\rho_{n-2}}S\hak \bigr)
= 1/2\cdot v
$$
by Lemma 5.18 and a theorem of Edidin-Graham [E-G] asserting, in this case,
the equality $c_{n-1}S=v$.
\medskip

If $n$ is even then
$$(\pi_1)_*P_{\rho_{n-1}}S\hak = -1/2\cdot (\pi_1)_*\bigl(c_{n-1}S \cdot
P_{\rho_{n-2}}S\hak \bigr)
= -1/2\cdot (\pi_1)_*\bigl(-v\cdot
P_{\rho_{n-2}}S\hak \bigr)
= 1/2\cdot v
$$
because the theorem of Edidin-Graham now asserts that $c_{n-1}S=-v$.
\medskip

Therefore the contribution of (1) is equal to
$$\aligned
(\pi_2\circ \pi_1)_* P_{\rho_{n-1}} S\hak
&=1/2\cdot(\pi_2)_* c_{n-1}(V_n/C) \\
&=1/2\cdot(\pi_2)_* \bigl(\sum_{p=0}^{n-1}(-1)^{n-p-1}
c_pV_n\cdot s_{n-p-1}C\bigr) \\
&=1/2\cdot(\pi_2)_* \bigl(\sum_{p=0}^{n-1}c_pV_n\cdot c^{n-p-1}\bigr)
=1/2\cdot \bigl(\sum_{p=0}^{n-1}c_pV_n\cdot s_{i-p}C\hak\bigr) \\
\endaligned$$
\bigskip
On the other hand the contribution of (2) is equal to
$$
1/2\cdot(\pi_2\circ\pi_1)_*\bigl(c_pR\hak
\cdot c^{n-p-1}\cdot P_{\rho_{n-2}}S\hak\bigr)
=1/2\cdot c_pR\hak\cdot(\pi_2)_*(c^{n-p-1})
=1/2\cdot c_pR\hak\cdot s_{i-p}(V_a\hak).
$$

Summing up the contributions of (1) and (2) we infer
$$
[\Omega(a)]=1/2\cdot \sum_{p=0}^i \bigl(c_pV_n+c_pR\hak\bigr)\cdot
s_{i-p}(V_a\hak),
$$
which is the asserted formula.
\qed
\enddemo

\bigskip\bigskip

\centerline {\bf 7. Two Schubert conditions}

\medskip

In this section we treat the classes of Schubert subschemes defined by two
Schubert conditions in the Lagrangian and odd orthogonal cases.
\smallskip

We consider first the Lagrangian case.
Our desingularization of $\Omega(n+1-i,n+1-j)$ in $G=LG_nV$
is given by the composition
(we use the notation of Section 1, $rank \ C=2$):
$$
\Cal F={L}G_{n-2}\bigl(C^{\perp}/C\bigr) @>\pi_1>> Fl(V_a\subset V_b) @>\pi_2>>
G,
$$
where $(a,b)=(n+1-i,n+1-j)$ and the element to be push forwarded via
$(\pi_2\pi_1)_*$ is
$\sum \widetilde Q_ID\hak\cdot \widetilde Q_{\rho_n\smallsetminus I}R\hak$,
the sum over all strict $I\subset\rho_n$.
Let $S$ be the tautological rank $(n-2)$ bundle on ${L}G_{n-2}
\bigl(C^{\perp}/C\bigr)$.
Using $[D\hak]=[S\hak]+[C\hak_{\Cal F}]$
and the linearity formula from Proposition 4.1 together with the factorization
property from Proposition 4.3, we have $(\pi_1)_*\widetilde Q_I D\hak \ne 0$
only if $(\pi_1)_*\widetilde Q_J S\hak \ne 0$ for some $J\subset I$.
By virtue of Theorem 5.10 (applied to $S\hak$),
$(\pi_1)_* \widetilde Q_JS\hak\ne 0$ only if $J\supset \rho_{n-2}$.
Consequently,
the unique strict $I$'s for which
$(\pi_1)_* \widetilde Q_ID\hak\ne 0$ must contain $\rho_{n-2}$,
i.e. they are of the form:
$I=\rho_n,\ I=(n,n-1,\ldots,\hat p,\ldots,1) =: I_p,\
I=(n,n-1,\ldots,\hat p,\ldots,\hat q,\ldots,1) =: I_{p,q}$
(here, $p$ and $q$ run over $\{1,\ldots,n\}$ and
the symbol \ "\ $\hat{\ }$\ " \
indicates the corresponding omission).
\medskip

We need the following technical lemma.

\proclaim{\bf Lemma 7.1}
If $rank \ C = 2$ then
{\parindent=30pt
\medskip
\item{(i) \ \ } $\widetilde Q_{I_p/\rho_{n-2}}(C\hak) = s_{n-1,n-p}(C\hak);$
\medskip
\item{(ii) \ } For $q<p$ , \
$
\widetilde Q_{I_{p,q}/\rho_{n-2}}(C\hak)=s_{n-q-1,n-p}(C\hak);
$
\medskip
\item{(iii) } For $0\leqslant v \leqslant n-2$, \
$
\widetilde Q_{\rho_n/(\rho_{n-2}+(2)^v)\widetilde{\ }}(C\hak)
= s_{n-v,n-v-1}(C\hak).
$
\par}
\endproclaim
\demo{Proof}
The proof is an easy application of the linearity formula from
Proposition 4.1 and is given here in case (i)
(the proofs of (ii) and (iii) being similar).

Denote the Chern roots of $C\hak$ by $x_1,x_2$.
Consider the skew Ferrers' diagram of $I_p/\rho_{n-2}$ and fill up
with "1" the
boxes, whose subtraction correspond to the summands in Proposition 4.1
applied to $x_1$ instead of $x_n$.
Then fill up with "2" the boxes, whose subtraction correspond to the
summands in
Proposition 4.1 applied to $x_2$ instead of $x_n$.
Of course it is impossible to have two "1" or two "2" in one row.
Also, the following configuration cannot appear:
\smallskip
$$
\CD
2&\ \ &   \cr
\hbox{x}&\ \ &1  \cr
\endCD
$$
\bigskip
\noindent
where the box "x" belongs to $D_{\rho_{n-2}}$ (Having two equal rows ending
with
${2\atop x}$ we use Proposition 4.3, thus we must subtract both boxes
instead of the higher one only). For example, for $n=6,\ p=3$ we get two
Ferrers' diagrams, one contained in another (depicted with "$\bullet$" and
"x"):
\bigskip
\vbox{\settabs 30\columns
\+&&&&&&&&&&&&  $\bullet$    &      &      &      &    &   &&&&&&&&&&&& \cr
\+&&&&&&&&&&&& x     & $\bullet$   &      &    &    &     &&&&&&&&&&&&\cr
\+&&&&&&&&&&&&  x    &  x    &  $\bullet$  &  $\bullet$  &      &
&&&&&&&&&&&&\cr
\+&&&&&&&&&&&&   x   &  x &  x    &  $\bullet$    &  $\bullet$    &
&&&&&&&&&&&&\cr
\+&&&&&&&&&&&& x   & x     &  x    &  x    &  $\bullet$    &$\bullet$
&&&&&&&&&&&&\cr
}
\bigskip
\noindent
and we have 3 possibilities:
\smallskip

\vbox{\bigskip
\settabs 30\columns
\+&&&&2  &  &  &  &  & &&&1 &  &  &  & & &&& 1 &  &  &  & & &&&&\cr
\+&&&&  & 2 &  &  &  & &&&  & 2&  &  &  &  &&&  &1  &  & & &  &&&&\cr
\+&&&&  &  &2 &1 &  &  &&&  &  &2 & 1&  &  &&&  &   &2 &1 &  &  &&&&\cr
\+&&&&  & &  & 2 & 1&  &&&  &  &  &2 &1  &  &&&  &  &  &2 & 1&  &&&&\cr
\+&&&&  &  &  &  &2 &1 &&&  &  &  &  & 2 &1 &&& &  &  &  &2 & 1 &&&&\cr
}
\bigskip
\noindent
giving $Q_{I_3/\rho_4}(x_1,x_2)=(x_1 x_2)^3(x_1^2+x_1x_2+x_2^2)=
s_{5,3}(x_1,x_2)$.
In general, arguing in the same way, we get
$$
\eqalign{
Q_{I_p/\rho_{n-2}}(x_1,x_2)
 &=(x_1x_2)^{n-p}(x_1^{p-1}+x_1^{p-2}x_2+\ldots+x_2^{p-1})= \cr
 &=e_2(x_1,x_2)^{n-p}s_{p-1}(x_1,x_2)=s_{n-1,n-p}(x_1,x_2). \qed \cr}
$$
\enddemo
\medskip

\proclaim{\bf Lemma 7.2}
With the above notation we have:
{\parindent=30pt
\medskip
\item{(i) \ \ } $(\pi_1)_*\bigl(\widetilde Q_{I_p}D\hak\bigr)
= s_{n-1,n-p}(C\hak)$;
\medskip
\item{(ii) \ \ } For $q<p$, \ \ $(\pi_1)_*\bigl(\widetilde
Q_{I_{p,q}}D\hak\bigr)
=s_{n-q-1,n-p}(C\hak)$;
\medskip
\item{(iii) \ \ } $(\pi_1)_*\bigl(\widetilde Q_{\rho_n}D\hak\bigr) =
\sum\limits_{k=0}^{n-2}(-1)^k c_{2k}V \cdot \Bigl[s_{n-k,n-k-1}(C\hak)-
s_{n-k+1,n-k-2}(C\hak)+\ldots$
\par\rightline{$\ldots+(-1)^{n-k}s_{2(n-k-1),1}(C\hak)\Bigr]$.}
\par}
\endproclaim
\demo{Proof}
Assertions (i) and (ii) follow immediately from Lemma 7.1(i),(ii) and
Theorem 5.10.
As for (iii), we have \ ( \ in the following, $(\pi_1)_*( \ other \ terms \
)=0$ \ ):
\medskip
\noindent
$(\pi_1)_*\bigl(\widetilde Q_{\rho_n}D\hak\bigr) =$
\medskip\Odst
$=(\pi_1)_*\Bigl[\sum\limits_{v=0}^{n-2}
\widetilde Q_{(\rho_{n-2}+(2)^v)\widetilde{\ }}
(S\hak)\cdot \widetilde Q_{\rho_n/(\rho_{n-2}+(2)^v)\widetilde{\ }}(C\hak_{\Cal
F})+ \ (other \ terms)
\Bigr]$
\medskip\Odst
$=\sum\limits_{v=0}^{n-2}(-1)^v c_{2v}(C^{\perp}/C)\cdot
\widetilde Q_{\rho_n/(\rho_{n-2}+(2)^v)\widetilde{\ }}(C\hak)$
\medskip\Odst
$=\sum\limits_{v=0}^{n-2}(-1)^v \Bigl[\sum\limits_{k+l=v}c_{2k}V\cdot
s_{2l}(C\oplus C\hak)\Bigr]\cdot s_{n-v,n-v-1}(C\hak)$
\medskip\Odst
$=\sum\limits_{k=0}^{n-2}(-1)^k c_{2k}V\cdot \Bigl[\sum\limits_{l=0}^{n-2-k}
(-1)^l s_{2l}(C\oplus C\hak)\cdot s_{n-k-l,n-k-l-1}(C\hak)\Bigr]$
\bigskip\Odst
$=\sum\limits_{k=0}^{n-2}(-1)^k c_{2k}V\cdot \Bigl[s_{n-k,n-k-1}(C\hak)-
s_{n-k+1,n-k-2}(C\hak)+\ldots$
\medskip\Odst
 \ \ \ \ \ \ \ \ \ \ \ \ \ \ \ \ \ \ \ \ \ \ \ \ \ \ \ \ \ \ \ \ \ \ \ \ \ \ \
 \ \ \ \ \ \ \ $\ldots+(-1)^{n-k}s_{2(n-k-1),1}(C\hak)\Bigr],$
\medskip
\noindent
where the above equalities follow from: Theorem 5.10, Lemma 1.1
and Pieri's formula
([Mcd1], [L-S1]); recall that $rank\ C=2$.
\qed
\enddemo
\bigskip

\proclaim{\bf Lemma 7.3}
Let $0<a<b$ and $k\ge l\ge 0$ be integers.
Let $C$ be the rank 2 tautological (sub)bundle of $\tau:Fl(a,b)\to X$.
Then
$$
\tau_*s_{k,l}(C\hak)=s_{k-(b-2)}(V_b\hak)\cdot s_{l-(a-1)}(V_a\hak)
-s_{k-(a-2)}(V_a\hak)\cdot s_{l-(b-1)}(V_b\hak).
$$
where we assume $s_h(-)=0$ for $h<0$.
\endproclaim

\demo{Proof}
Let $C_1\subset C_2=C$ be the tautological subbundles on $Fl(a,b)$,
$C_1\subset V_a$, $C_2\subset V_b$; $rank\ C_h=h$, $h=1,2$.
Let $x_1=c_1(C_1\hak)$ and $x_2=c_1\bigl((C_2/C_1)\hak\bigr)$.
The flag bundle $\tau:Fl(a,b)\to X$ is equal to the
composition:
$$
\Bbb P\bigl(V_b/C_1 \bigr) @>\tau_2>> \Bbb P(V_a) @>\tau_1>> X.
$$
We have
$$
\tau_*s_{k,l}(C\hak)=\tau_*\Bigl[(x_1x_2)^l\bigl(x_1^{k-l}+x_1^{k-l-1}x_2+\ldots
+x_1x_2^{k-l-1}+x_2^{k-l}\bigr)\Bigr].
$$
The assertion now follows by applying to all summands the well known formulas:
$$
\eqalign{
(\tau_2)_*(x_2^p)
&=s_{p-(b-2)}(V_b / C_1)\hak=s_{p-(b-2)}(V_b\hak)-s_{p-(b-2)-1}(V_b\hak)
\cdot x_1,  \cr
(\tau_1)_*(x_1^p)
&=s_{p-(a-1)}(V_a\hak)}
$$
and simplifying.
\qed
\enddemo
\bigskip\smallskip

\proclaim{\bf Theorem 7.4}
For $n\ge i>j>0$ one has in $A^*(G)$ with $a=n+1-i$, $b=n+1-j$,
\bigskip

$\bigl[\Omega(a,b)\bigr]=\sum\limits_{p>q\ge 0\atop p\le i, q\le j}
\widetilde Q_{p,q}R\hak\cdot
\bigl(s_{i-p}(V_a\hak)\cdot s_{j-q}(V_b\hak)-s_{i-q}(V_a\hak)\cdot
s_{j-p}(V_b\hak)\bigr)+$
\medskip

$+\sum\limits_{p=0}^{i-1}\sum\limits_{t\ge 1}(-1)^{p+t-1}c_{2p}V\cdot
\bigl(s_{i-p-t}(V_a\hak)
\cdot s_{j-p+t}(V_b\hak)-s_{i-p+t}(V_a\hak)\cdot s_{j-p-t}(V_b\hak)\bigr),$

\bigskip

\noindent
where we assume $s_h(-)=0$ for $h<0$.
\endproclaim

\bigskip

\demo{Proof}
It follows from Lemma 7.2 that

\bigskip\noindent
$\bigl[\Omega(a,b)] =
 \sum\limits_{0\le q<p}(\pi_2)_*\bigl(s_{n-q-1,n-p}(C\hak)\bigr) \cdot
\widetilde Q_{p,q}R\hak +$
\bigskip
$ \ \ +\sum\limits_{p=0}^{n-2}(-1)^p c_{2p}V \cdot (\pi_2)_*
\Bigl[s_{n-p,n-p-1}(C\hak)-
s_{n-p+1,n-p-2}(C\hak)+\ldots$
\bigskip
\rightline{$\ldots +(-1)^{n-p}s_{2(n-p-1),1}(C\hak)\Bigr].$}
\bigskip\noindent
Applying Lemma 7.3 to $\pi_2:Fl(a,b)\to X$, the assertion follows.
\qed
\enddemo
\bigskip\medskip

\example{\bf Example 7.5} \ 1. For $i=2, j=1$ and any $n$ the formula reads:
\bigskip
\noindent
$
\widetilde Q_{21}R\hak+\widetilde Q_{2}R\hak \cdot s_1V_n\hak+
\widetilde Q_{1}R\hak \cdot\bigl(s_1V_{n-1}\hak \cdot s_1V_n
\hak-s_2V_{n-1}\hak\bigr)+$

\medskip
\noindent$
\ \ \ \ \ \ \ \ \ \ \ \ \ \ \ \ \ \ \ \ \ \
+\bigl(s_1V_{n-1}\hak\cdot s_2V_n\hak-s_3V_{n-1}\hak-
s_3V_n\hak-
c_2V\cdot s_1V_n\hak\bigr)=$

\medskip
\noindent$
 \ \ \ \ \ \ \ \ \ \ \ \ \ \ \ \ \ \ \ \ \ \ \ \ \ \ \ \ \ \ \ \ \ \ \ \ \ \ \
\ \ \ \ \ \ \ \ \
=\widetilde Q_{21}R\hak+\widetilde Q_{2}R\hak \cdot \widetilde Q_1V_n\hak+
\widetilde Q_{1}R\hak \cdot \widetilde Q_2V_n\hak+\widetilde Q_{21}V_n\hak.$
\bigskip\smallskip

\noindent
2. For $i=3, j=1$ and any $n$ one obtains, with $\widetilde Q_{p,q}=\widetilde
Q_{p,q}R\hak, \ s_k=s_k(V_{n-2}\hak)$ and $s_k'=s_k(V_n\hak)$, the expression:
\bigskip

\noindent
$\widetilde Q_{31}+\widetilde Q_3\cdot s_1'+\widetilde Q_{21}\cdot s_1+
\widetilde Q_2\cdot s_1\cdot s_1'+
\widetilde Q_1\cdot (s_2\cdot s_1'-s_3)+$
\smallskip

\ \ \ \ \ \ \ \ \ \ \ \ \ \ \ \ \ \ \ \ \ \ \ \ \ \ \ \ \ \ \ \ \ \ \ \ \ \ \ \
\ \ \ \ \ \ \
$+s_2\cdot s_2'-s_4-s_1\cdot s_3'+s_4'
-c_2V\cdot (s_1\cdot s_1'-s_2')+c_4V.$
\bigskip

\noindent
3. For $i=3, j=2$ and any $n$ one obtains, with $\widetilde Q_{p,q}=\widetilde
Q_{p,q}R\hak$ and $s_{k,l}=s_{k,l}(V_{n-1}\hak)$, the expression:
\medskip

\noindent
$\widetilde Q_{32}+\widetilde Q_{31}\cdot s_1+\widetilde Q_3\cdot s_2+
\widetilde Q_{21}\cdot s_{11}+\widetilde Q_2\cdot s_{21}+\widetilde Q_1\cdot
s_{22}+s_{32}-s_{41}+s_5
-c_2V\cdot (s_{21}-s_3)+c_4V\cdot s_1.$
\endexample

\medskip

More generally we have:

\proclaim{\bf Corollary 7.6} With the above notation and $j=i-1, \ s_{k,l}=
s_{k,l}(V_{n+2-i}\hak)$, the class $[\Omega(a,b)]$ equals
\medskip

$\sum\limits_{i\ge p>q\ge 0} \widetilde Q_{p,q}R\hak \cdot
s_{i-1-q,i-p}+\sum\limits_{p=0}^{i-1}(-1)^p c_{2p}V\cdot \sum\limits
_{h=0}^{i-1-p}(-1)^h s_{i-p+h,i-1-p-h}.$

\endproclaim

\bigskip\medskip

Consider now the odd orthogonal case.
Our desingularization in case $a.=(a,b):=(n+1-i,n+1-j)$
is given by the composition ($rank \ C=2$):
$$
{O}G_{n-2}\bigl(C^{\perp}/C\bigr) @>\pi_1>> Fl(V_a\subset V_b) @>\pi_2>> G.
$$

Then the analog of Lemma 7.1 reads (with the notation explained before this
lemma):

\proclaim{Lemma 7.7}:
\smallskip
\noindent
(i) \quad $\widetilde P_{I_p/\rho_{n-2}}(C\hak)=1/2\cdot s_{n-1,n-p}(C\hak).$
\qquad (ii) \quad
$\widetilde P_{I_{p,q}/\rho_{n-2}}(C\hak)=s_{n-q-1,n-p}(C\hak).$
\smallskip \noindent
(iii) \Odst\Odst $\widetilde P_{\rho_n/(\rho_{n-2}+(2)^v)\widetilde{\
}}(C\hak)=
s_{n-v,n-v-1}(C\hak)\,,\quad 0<v\le n-2\,,$
\smallskip
\centerline{and\quad $\widetilde P_{\rho_n/\rho_{n-2}}(C\hak)=1/4\cdot
s_{n,n-1}(C\hak).$}
\endproclaim

\bigskip

The element to be push forwarded via $(\pi_2\pi_1)_*$ is
$\sum \widetilde P_ID\hak\cdot \widetilde P_{\rho_n\smallsetminus I}R\hak$,
the sum over all strict $I\subset\rho_n$ (the notation as in Section 1).
The analog of Lemma 7.2 reads:

\proclaim{\bf Lemma 7.8}:
\smallskip
\noindent
(i) $(\pi_1)_*(\widetilde P_{I_p}D\hak)=1/2\cdot s_{n-1,n-p}(C\hak).$
\quad (ii)
$q<p\quad (\pi_1)_*(\widetilde P_{I_{p,q}}D\hak)=s_{n-q-1,n-p}(C\hak).$
\smallskip
\noindent
(iii) \qquad $(\pi_1)_*(\widetilde P_{\rho_n}D\hak)=
1/4\cdot\sum\limits_{k=0}^{n-2}(-1)^kc_{2k}\,V\cdot
\bigl[s_{n-k,n-k-1}(C\hak)-s_{n-k+1,n-k-2}(C\hak)+$
\smallskip
\rightline{$\ldots (-1)^{n-k}\cdot s_{2(n-k-1),1}(C\hak)\bigr]\,.$}
\endproclaim
\medskip

Consequently, the analog of Theorem 7.4 now reads:

\proclaim{\bf Theorem 7.9}
For $n\ge i>j>0$ one has in $A^*(G)$ with $a=n+1-i$, $b=n+1-j$,
\medskip
\noindent
$\bigl[\Omega(a,b)\bigr]=\sum\limits_{p>q>0 \atop p\le i,q\le j}
\widetilde P_{p,q}\,R\hak\cdot
\Bigl(s_{i-p}(V_a\hak)\cdot s_{j-q}(V_b\hak)-s_{i-q}(V_a\hak)\cdot
s_{j-p}(V_b\hak)\Bigr)+$
\medskip
\centerline{$1/2\cdot\sum\limits_{p>0\atop p\le i} \widetilde P_p\,R\hak\cdot
\Bigl(s_{i-p}(V_a\hak)\cdot s_j(V_b\hak)-s_i(V_a\hak)\cdot s_{j-p}(V_b\hak)
\Bigr)+$}
\medskip
\rightline{$+1/4\cdot\sum\limits_{p=0}^{i-1}\,\sum\limits_{t\ge 1}
(-1)^{p+t-1}c_{2p}\,V\cdot\Bigl(s_{i-p-t}(V_a\hak)\cdot s_{j-p+t}(V_b\hak)
-s_{i-p+t}(V_a\hak)\cdot s_{j-p-t}(V_b\hak)\Bigr).$}
\endproclaim

\medskip
For instance, invoking Example 7.5, the formula
reads for $i=2, j=1$ and any $n$:
\smallskip
$$
[\Omega(n-1,n)] = \widetilde P_{21}R\hak+\widetilde P_{2}R\hak \cdot \widetilde
P_1V_n\hak+
\widetilde P_{1}R\hak \cdot \widetilde P_2V_n\hak+\widetilde P_{21}V_n\hak.
$$

\bigskip\smallskip

\centerline {\bf 8. An operator proof of Proposition 3.1}

\medskip

The goal of this section is to provide another proof
of Proposition 3.1 and its odd orthogonal analogue
by using divided differences operators.
We start with the Lagrangian case.
Let $X_n=(x_1,\ldots,x_n)$ be a sequence of indeterminates.
Recall (see Section 5) that the symplectic Weyl group $W_n$ is isomorphic to
$S_n\ltimes\zz_2^n$ and the
elements of $W_n$ are identified with "barred permutations":
if $w=(\sigma,\tau),\ \sigma\in S_n,\ \tau \in\zz_2^n$ then
we write $w$ as the sequence $(w_1,\ldots,w_n)$
endowed with bars on places where $\tau_i=-1$.
In particular, $w_0=(\nkr1,\nkr2,\ldots,\nkr{n})$ is the longest  element of
$W_n$.
Consider in $W_n$ the poset $W^{(n)}$ of minimal length left coset
representatives of
$W_n$ modulo its subgroup generated by reflections corresponding
to the simple roots $\varepsilon_1-\varepsilon_2,\ldots,\varepsilon_{n-1}-
\varepsilon_n$ (in the standard notation):
$$
W^{(n)}=\Bigl\{(\nkr z_1>\nkr z_2>\ldots>\nkr z_l;y_1<\ldots<y_{n-l})\in
W_n,\ l=0,1,\ldots,n\Bigr\}.
$$
The assignment $w=(\nkr z_1,\ldots,\nkr z_l;y_1,\ldots,y_{n-l})\mapsto
I=(z_1,\ldots,z_l)$ establishes a bijection between the poset $W^{(n)}$
and the poset of all strict partitions contained in $\rho_n$.

\smallskip
One has divided differences
$\partial_w:\zz[X_n]\to\zz[X_n]$ \ ($w\in W_n$) \
i.e. operators of degree $-l(w)$, whose definition has been explained in
Section 5.
\medskip

Fix now an integer $0<k<n$ and denote:
$$
w^{(k)}:=(\nkr{n},\nkr{n-1},\ldots,\nkr{k+1};1,2,\ldots,k).
$$

Observe first that for a strict partition $I\subset\rho_n$ of length $l(I)$,
$\partial_{w^{(k)}}\widetilde Q_I(X_n)\ne 0$ only if $l(I)\ge n-k$.
(This is because $\partial_{w^{(k)}}$ decreases the degree by
$l(w^{(k)})=n+(n-1)+\ldots+(k+1)$.)
More precisely, writing $X_n\hak =(-x_1,\ldots,-x_n)$, we have:
\smallskip

\proclaim{\bf Proposition 8.1}
For a strict partition $I$ of length $\ge n-k$, $\partial_{w^{(k)}}
\widetilde Q_I(X_n\hak)
\ne 0$ iff $I\supset(n,n-1,\ldots,k+1)$.
In the latter case, writing $I=(n,n-1,\ldots,k+1,j_1,\ldots,j_l)$,
where $j_l>0$ and $l\le k$,
one has in $\zz[X_n]$,
$$
\partial_{w^{(k)}}\widetilde Q_I(X_n\hak) =
\widetilde Q_{j_1,\ldots,j_l}(X_n\hak).
$$
\endproclaim

\demo{Proof}
Let $I$ be a strict partition of length $h\ge n-k$.
Let
$$
w_I=(\nkr \sigma_1,\nkr \sigma_2,\ldots,\nkr
\sigma_h;\sigma_{h+1},\ldots,\sigma_n)
$$
be the element of $W^{(n)}$ corresponding to $I$. Then taking into account
that
$$
(w^{(k)})^{-1}=(n-k+1,n-k+2,\ldots,n;\nkr {n-k},\nkr {n-k-1},\ldots,\nkr1),
$$
we get \ \
$w_I\circ(w^{(k)})^{-1}=$\bigskip\noindent$
\bigl(\nkr \sigma_{n-k+1}>\nkr \sigma_{n-k+2}>\ldots>\nkr \sigma_h,
\sigma_{h+1}<\sigma_{h+2}<\ldots<\sigma_n,
\sigma_{n-k}<\sigma_{n-k-1}<\ldots<\sigma_1\bigr).
$\medskip

We have $l(w_I)=\sigma_1+\ldots+\sigma_h, \ \ l(w^{(k)})=n+(n-1)+\ldots+(k+1)$,
\ and
$$
l\bigl(w_I\circ(w^{(k)})^{-1}\bigr)=\sigma_{n-k+1}+\sigma_{n-k+2}+\ldots+\sigma_h+
\sum\limits_{j=1}^{n-h}card\bigl\{1\le p\le n-k\ | \
\sigma_p<\sigma_{h+j}\bigr\}
$$
by Lemma 5.1.
Thus, denoting the above sum $\sum\limits_{j=1}^{n-h}(\ldots)$ by
$\sum$, we get:
\smallskip
$ l(w_I)-l(w^{(k)})-l\bigl(w_I\circ(w^{(k)})^{-1}\bigr)=$
\medskip
\hfill$=\sigma_1+\ldots+\sigma_{n-k}-\bigl(n+(n-1)+\ldots+(k+1)\bigr)-\sum$.
\medskip
Now, a necessary condition for $\partial_{w^{(k)}}\widetilde Q_I(X_n\hak)\ne0$
is:
$$
\sigma_1+\ldots+\sigma_{n-k}-\bigl(n+(n-1)+\ldots+(k+1)\bigr)-\sum=0,
$$
which implies \
$(\sigma_1,\ldots,\sigma_{n-k})=(n,n-1,\ldots,k+1)$ \ \
and \ \ $\sum=0$, i.e., $\sigma_n<\sigma_{n-k}$.
(Using the theory from [B-G-G], [D1,2] and the result from [P2] recalled
in Theorem 2.1, the just proved assertion easily implies that
$\partial_{w^{(k)}}\widetilde Q_I(X_n\hak) =
\widetilde Q_{j_1,\ldots,j_l}(X_n\hak) \ (mod \ \Cal I)$).
We will now prove directly that for $I=(n,n-1,\ldots,k+1,j_1,\ldots,j_l)$ one
has
$\partial_{w^{(k)}}\widetilde Q_I(X_n\hak) =
\widetilde Q_{j_1,\ldots,j_l}(X_n\hak)$ already in $\zz[X_n]$. Observe
that
$$\partial_{w^{(k)}}=(\partial_{k}...\partial_1 \partial_0)...
(\partial_{n-2}...\partial_1 \partial_0)
(\partial_{n-1}...\partial_1 \partial_0).$$
The proof is by induction on $n-k-1$. For $n-k-1=0$,
one has ($J=(j_1,...,j_l))$:

$$\partial_{n-1}...\partial_1 \partial_0\bigl(\widetilde Q_{n,J}
(X_n\hak)\bigr) =
\partial_{n-1}...\partial_1 \partial_0\bigl(e_n(X_n \hak)\cdot
\widetilde Q_{J}(X_n\hak)\bigr)$$
$$=\partial_{n-1}...\partial_1 \bigl((-x_2)...(-x_n)
\widetilde Q_{J}(X_n\hak) -
e_n(X_n\hak)\cdot \partial_0\widetilde Q_{J}(X_n\hak)\bigr)$$
$$=\widetilde Q_{J}(X_n\hak) - e_n(X_n\hak)\cdot
\partial_{n-1}...\partial_1 \partial_0
\bigl(\widetilde Q_{J}(X_n\hak)\bigr) = \widetilde Q_J(X_n \hak),$$

\noindent
where the vanishing of the second summand in the latter difference
follows from the just proved first assertion.
\smallskip

The induction step goes as follows. By the equality proved above,

$$(\partial_{k}...\partial_1 \partial_0)...
(\partial_{n-2}...\partial_1 \partial_0)
(\partial_{n-1}...\partial_1 \partial_0)
\bigl(\widetilde Q_{n,n-1,...,k+1,J}(X_n\hak)\bigr)$$
$$(\partial_{k}...\partial_1 \partial_0)...
(\partial_{n-2}...\partial_1 \partial_0)
\bigl(\widetilde Q_{n-1,n-2,...,k+1,J}(X_n\hak)\bigr)$$
$$(\partial_{k}...\partial_1 \partial_0)...
(\partial_{n-2}...\partial_1 \partial_0)
\bigl(\sum_{i\geq 0}(-x_n)^i\sum \widetilde Q_{I}(X_{n-1}\hak)\bigr)$$

\noindent
where the sum is over all partitions
$I \subset (n-1,n-2,...,k+1,J)$ such that the diagram
$(n-1,n-2,...,k+1,J)/I$ is of weight $i$ and has at most one box in
every row (use the linearity formula, i.e. Proposition 4.1).
Each time we get two equal parts $p$ in a partition $I$ such that
$\widetilde Q_{I}(X_{n-1}\hak)$ appears in the expression, we factor out
$\widetilde Q_{p,p}(X_{n-1}\hak)$ by Proposition 4.3.
The last sum can be rewritten in the form:

$$\sum_{i\geq 0}\sum_M (-x_n)^i\widetilde Q_{M}(X_{n-1}\hak) f_M +
\sum_{i\geq 0}\sum_N(-x_n)^i\widetilde Q_{N}(X_{n-1}\hak) g_N,$$
where $M$ (resp. $N$) runs over the so-obtained partitions contained in
the partition
$(n-1,n-2,...,k+1,J)$  where some box is removed from the first
$n-k-1$ places (resp. no box is subtracted from the first $n-k-1$ places),
and $f_M$ (resp. $g_N$) denotes the corresponding monomial
in the elements $\widetilde Q_{p,p}(X_{n-1}\hak)$
obtained by factoring out. By the first assertion (applied to $X_{n-1}\hak$)
we know that our operator anihilates the first sum. By the
induction assumption
we get (with $N=(n-1,n-2,...,k+1,J')$ and $g_{J'} = g_N$)

$$(\partial_{k}...\partial_1 \partial_0)...(\partial_{n-2}...
\partial_1 \partial_0)
\bigl(\sum_{i\geq 0}\sum_N (-x_n)^i \widetilde Q_{N} (X_{n-1}\hak) g_N\bigr)$$
$$=\sum_{i\geq 0}\sum_{J'} (-x_n)^i \widetilde Q_{J'} (X_{n-1}\hak) g_{J'} =
\widetilde Q_{J}(X_n\hak)$$

\noindent
by the factorization property and the linearity formula, now used backwards.
\qed
\enddemo
\medskip

We now pass to a geometric interpretation of the proposition.
The setup and the notation is the same as in the proof of Proposition 3.1:
$V\to B$ - rank $2n$ vector bundle endowed with a nondegenerate
symplectic form,
$X=LG_n V$, $V_n$ denotes here
the tautological subbundle on $X$ and
$p:\Cal F \to X$ is the composition (see Section 1):
$$
LG_{n-k}(C^{\perp}/C) @>\pi_1>>G_k(V_n)@>\pi_2>>X,
$$
where $C$ is the tautological rank $k$ bundle on $G_k(V_n)$.
The tautological rank $n-k$ subbundle $S$ for $LG_{n-k}(C^{\bot}/C)$
is identified with $D/C_{\Cal F}$ where $D$ is rank $n$ tautological
subbundle on  $\Cal F$.
Let $r_1,...,r_n$ be the Chern roots of $V_n$ and $d_1,...,d_n$ - the Chern
roots of $D$.
Since $C_{\Cal F}\subset (V_n)_{\Cal F}$ and $C_{\Cal F}\subset D$,
we can assume that $r_1=d_1,...,r_k=d_k$ are the Chern roots of $C$.

\smallskip
\noindent
\underbar{Claim}: For any symmetric polynomial $f$ in $n$ variables,
$$
(\pi_1)_*\left(f(d_{k+1},...,d_n,d_1,...,d_k)\right)=
(\partial_vf)(r_{k+1},...,r_n,r_1,...,r_k),
$$
where $v=(\overline{n-k},\overline{n-k-1},...,\overline{1},n-k+1,...,n)$.
\smallskip

Indeed, for the Chern roots $d_{k+1},...,d_n$ of $S$ one has by
Proposition 5.8,
$$
(\pi_1)_*\left(f(d_{k+1},...,d_n,d_1,...,d_k)\right)=
(\pi_1)_*\left(f(d_{k+1},...,d_n,r_1,...,r_k)\right)$$
$$
=(\partial_vf)(d_{k+1},...,d_n,r_1,...,r_k).
$$
We know by Proposition 5.9 that $\partial_vf$ is a polynomial symmetric
in the squares of the first $n-k$ variables.
By Lemma 1.1 we have
$$
\eqalign{
[S]+[S\hak]&=[(C^{\perp}/C)_{\Cal F}]
=[V_{\Cal F}]-[C_{\Cal F}]-[C\hak_{\Cal F}] \cr
&=[(V_n)_{\Cal F}]+[(V_n\hak)_{\Cal F}]-[C_{\Cal F}]-[C\hak_{\Cal F}]
=[(V_n)_{\Cal F}/C_{\Cal F}]+[((V_n)_{\Cal F}/C_{\Cal F})\hak].      \cr
}
$$
Hence, for the Chern roots $r_{k+1},...,r_n$ of $V_n/C$,
$$
(\partial_vf)(d_{k+1},...,d_n,r_1,...,r_k)=
(\partial_vf)(r_{k+1},...,r_n,r_1,...,r_k)
$$
and the claim is established.
\medskip

We are now in position to give
\demo{\bf Another proof of Proposition 3.1}
\smallskip

By virtue of the previous proposition it suffices to show that for every
symetric polynomial $f$ in $n$ variables $p_*\bigl( f(d_1,...,d_n)\bigr)=
\bigl(\partial_{w^{(k)}}f\bigr)(r_1,...,r_n)$.
For a polynomial $g$ symmetric in
the first $n-k$ - and in the last $k$ variables, one has
$$
(\pi_2)_*\left(g(r_{k+1},...,r_n,r_1,...,r_k)\right)=
(\partial_ug)(r_1,...,r_n),
$$
where $u=(k+1,...,n,1,2,...,k)$ (see [L2], [P2] and [Br]).
(This can be proved using a reasoning similar to the one in the proof
of Proposition 5.8 above.)
Since $w^{(k)}=u\circ v$ and $l(w^{(k)})=l(u)+l(v)$,
we thus have, invoking the claim:

$$
\eqalign{
p_*\bigl(f(d_1,...,d_n)\bigr)
&=p_*\bigl(f(d_{k+1},...,d_n,d_1,...,d_k)\bigr)
=\pi_{2*}\Bigl(\pi_{1*}\bigl(f(d_{k+1},...,d_n,d_1,...,d_k)\bigr)\Bigr) \cr
&=\pi_{2*}\bigl((\partial_vf)(r_{k+1},...,r_n,r_1,...,r_k)\bigr)
=\bigl(\partial_u(\partial_vf)\bigr)(r_1,...,r_n) \cr
&=\bigl((\partial_u\circ\partial_v)f\bigr)(r_1,...,r_n)
=(\partial_{w^{(k)}}f)(r_1,...,r_n)                         \cr
}$$
\noindent
which is the desired assertion.
\qed
\enddemo
\smallskip

In the odd orthogonal case, by replacing $\widetilde Q$-polynomials
by $\widetilde P$-polynomials
and arguing in the same way as above, one proves the following
proposition.

\proclaim{\bf Proposition 8.2}
For a strict partition $I$ of length $\ge n-k$, $\partial_{w^{(k)}}
\widetilde P_I(X_n\hak)
\ne 0$ iff $I\supset(n,n-1,\ldots,k+1)$.
In the latter case, writing $I=(n,n-1,\ldots,k+1,j_1,\ldots,j_l)$,
where $j_l>0$ and $l\le k$,
one has in $\zz[X_n]$,
$$
\partial_{w^{(k)}}\widetilde P_I(X_n\hak) =
\widetilde P_{j_1,\ldots,j_l}(X_n\hak).
$$
\endproclaim

Let $V\to B$ be a rank $2n+1$ vector bundle endowed with a nondegenerate
orthogonal form, $X=OG_n V$ and $V_n$ denote the tautological
subbundle on $X$. Then, by an appropriate interpretation of the
Gysin map associated with
the composition:
$$
OG_{n-k}(C^{\perp}/C) @>\pi_1>>G_k(V_n)@>\pi_2>>X,
$$
where $C$ is the tautological rank $k$ bundle on $G_k(V_n)$, one gets
another proof of Proposition 3.4.
\smallskip

We refer the reader to [L-P-R] for another operator treatment of
$\widetilde Q$- and $\widetilde P$-polynomials and their generalizations.
\bigskip\smallskip

\centerline {\bf 9. Main results in the generic case}
\medskip

Let $V$ be a rank $2n$ vector bundle over a smooth pure-dimensional
scheme $X$ endowed with
a nondegenerate symplectic form.
Let $E$ and $F.$ : $F_1\subset F_2\subset ... \subset F_n=F$ be Lagrangian
subbundles of $V$ with $rank \ F_i=i$ and $rank \ E=n$. For a given sequence
$a.=(1\le a_1< ...< a_k\le n)$, we are interested
in a locus

$$
D(a.):=\Bigl\{\ x\in X | \ dim\bigl(E\cap F_{a_p}\bigr)_x\geqslant p, \
p=1,...,k \Bigr\}.
$$

Let $G=LG_nV$ and let $R\subset V_G$ be the tautological rank $n$ subbundle on
$G$.
By a well known universality property of Grassmannians there exists a
morphism $s: X \to G$ such that $E=s^*R$.
Therefore (in the set-theoretic sense) we have:
$$
D(a.)=s^{-1}(\Omega (a.;F.)),
$$
where
$$
\Omega (a.;F.)=
\Bigl\{\ g\in G | \ dim(R\cap F_{a_p})_g\geqslant p, \ p=1,...,k \ \Bigr\}.
$$
We take this equality as the definition of a scheme structure on $D(a.)$,
i.e., $D(a.)$ is defined in $X$ by the inverse image ideal sheaf (see [Ha,
p.163]):
$s^{-1}\Cal I(\Omega (a.;F.))\cdot \Cal O_X$ where $\Cal I(\Omega (a.;F.)$
is the ideal sheaf defining $\Omega$ in $G$.
It follows from the main theorem of [DC-L] that
$\Omega (a.;F.)$ is a Cohen-Macaulay scheme.
Hence, by [K-L, Lemma 9] we get
$[D(a.)]=s^*[\Omega (a.;F.)]$ provided $D(a.)$ is
either empty or of pure codimension equal to the codimension
of $\Omega (a.;F.)$ in $G$. Therefore, having a formula for the
fundamental class of $\Omega (a.;F.)$ given by a polynomial $P$ in
$c.(R)$ and $c.(F_{a_p})_G$, $p=1,...,k$, the formula for $D(a.)$ becomes
$P\bigl(c.(E), c.(F_{a_p})_{p=1,...,k}\bigr)$. Moreover,
by using the Chow groups
for singular schemes and a technique from [F] one can prove the following
refinement of the above. If $X$ is a pure-dimensional
Cohen-Macaulay scheme and $D(a.)$
is either empty or of pure codimension equal to the codimension
of $\Omega (a.;F.)$ in $G$ then the class of $D(a.)$ in the Chow
group of $X$ equals $P\bigl(c.(E), c.(F_{a_p})_{p=1,...,k}\bigl)\cap [X]$.
This reasoning (with obvious modifications) also applies, word by word,
to the case of rank $2n+1$ vector bundle endowed with a nondegenerate
orthogonal form.
\smallskip

In particular,
for $a.=(n-k+1,...,n)$ we have by Proposition 3.2:

\proclaim{\bf Theorem 9.1}
If $X$ is a pure-dimensional Cohen-Macaulay
scheme and the subscheme
$$D^k=\bigl\{x\in X| \ dim(E\cap F)_x \geqslant k\bigr\}$$ is
either empty
or of pure
codimension $k(k+1)/2$ in $X$, then the class of $D^k$ (endowed with
the above scheme structure) in the Chow group
of $X$ equals
$$
[D^k] = \Bigl(\sum \widetilde Q_IE\hak\cdot
\widetilde Q_{\rho_k\smallsetminus I}F\hak\Bigr)\cap[X],
$$
where the sum is over all strict partitions $I\subset\rho_k$.
\endproclaim

\smallskip

\example{\bf Example 9.2}
The expressions giving the classes for successive $k$ are:
{\parindent=30pt
\smallskip
\item{k=1\ \ } $\widetilde Q_1E\hak+\widetilde Q_1F\hak$;
\smallskip
\item{k=2\ \ } $\widetilde Q_{21}E\hak+\widetilde Q_2E\hak\cdot \widetilde
Q_1F\hak+
               \widetilde Q_1E\hak\cdot \widetilde Q_2F\hak+\widetilde
Q_{21}F\hak$;
\smallskip
\item{k=3\ \ } $\widetilde Q_{321}E\hak+\widetilde Q_{32}E\hak\cdot \widetilde
Q_1F\hak+
               \widetilde Q_{31}E\hak\cdot \widetilde Q_2F\hak+\widetilde
Q_{21}E\hak\cdot \widetilde Q_3F\hak+
	       \widetilde Q_3E\hak\cdot \widetilde Q_{21}F\hak+   \hfill\break
	       \widetilde Q_2E\hak\cdot \widetilde Q_{31}F\hak+\widetilde
Q_1E\hak\cdot \widetilde Q_{32}F\hak+
	       \widetilde Q_{321}F\hak$.
\par}

\endexample

\smallskip

For $a.=(n+1-i)$ we get:

\proclaim {\bf Theorem 9.3}
Let $X$ be a pure-dimensional Cohen-Macaulay scheme and assume that the
subscheme $S^i=
\{ x\in X | \ dim(E\cap F_{n+1-i})_x \geqslant 1\}$ is either empty or
of pure codimension $i$ in $X$. Then
$$
[S^i] = \bigl(\sum_{p=0}^i c_pE\hak \cdot
s_{i-p} F_{n+1-i}\hak\bigr)\cap [X].
$$
\endproclaim
\medskip
\example{\bf Example 9.4}
The expressions giving the classes for successive $i$ are:
{\parindent=30pt
\smallskip
\item{i=1\ \ } $c_1E\hak + s_1F\hak$;
\smallskip
\item{i=2\ \ } $c_2E\hak + c_1E\hak s_1F_{n-1}\hak + s_2F_{n-1}\hak$;
\smallskip
\item{i=3\ \ } $c_3E\hak + c_2E\hak s_1F_{n-2}\hak + c_1E\hak s_2F_{n-2}\hak +
s_3F_{n-2}\hak$.
\par}

\endexample
\smallskip
The theorem is a globalization to degeneracy loci of Proposition 6.1.
Also other formulas from Sections 6 and 7 admit analogous globalizations.
We concentrate ourselves on a solution to J. Harris' problem for
Mumford-type degeneracy loci mentioned in the Introduction.

\smallskip

The odd orthogonal analog of Theorem 9.1 is a consequence of
Proposition 3.4 and reads as follows:

\proclaim{\bf Theorem 9.5}
Let $X$ be a pure-dimensional Cohen-Macaulay scheme over a field of
characteristic different from 2. Suppose that $V$ is a rank $2n+1$ vector
bundle endowed with a nondegenerate orthogonal form. Let $E$ and $F$
be two rank $n$ isotropic subbundles of $V$. If
the subscheme
$$D^k=\bigl\{x\in X | \ dim(E\cap F)_x\ge
k\bigr\}$$ is either empty or of pure
codimension $k(k+1)/2$ in $X$, then the class of $D^k$ in the Chow group
of $X$ equals
$$
\Bigl(\sum \widetilde P_IE\hak\cdot
\widetilde P_{\rho_k\smallsetminus I}F\hak\Bigr)\ \cap\ [X],
$$
where the sum is over all strict partitions $I\subset\rho_k$.
\footnote {\eightrm Observe that though the $\widetilde P$-polynomials of
a vector bundle
are defined only when the Chern classes of the vector bundle are divisible
by $2$, the integrality property of the classes obtained in the theorem
holds true. One argues as in the preceding footnote, taking into account that
the base change argument [K-L, Lemma 9] preserves the integrality
too. The same remark applies to the even orthogonal case (Theorem 9.6 below).}

\endproclaim
\medskip

Let now $V$ be a rank $2n$ vector bundle over a connected pure-dimensional
scheme $X$
endowed with a nondegenerate orthogonal form.
Let $E$ and $F.$ : $F_1\subset F_2\subset\ldots\subset F_n=F$ be isotropic
subbundles of $V$ with $rank\ F_i=i$ and $rank\ E=n$.
One should be careful here with the definition of $D(a.)$.
For a given sequence $a_.=(1\le a_1<\ldots<a_k\le n)$, where $k$ is such that
$dim(E\cap F)_x\equiv k (mod \  2)$ if $a_k=n$,
we are interested in
the locus
$$
D(a_.)=\Bigl\{x\in X | \ dim(E\cap F_{a_p})_x\ge p \ , \ p=1,\ldots,k \Bigr\}.
$$
There is a morphism $s=(s',s''): X \to OG_n'V\cup OG_n''V$ such that
$s^*R=E$ where $R$ is the tautological rank $n$ subbundle on
$OG_n'V\cup OG_n''V$.
We have (in the scheme -- theoretic sense) that if $k\equiv n \ (mod \  2)$
then
$$D(a_.)=(s')^{-1}\Omega\bigl(a_.;(F_.)_{OG_n'V}\bigr);$$
and if
$k\equiv n+1 (mod \  2)$ then $$D(a_.)=(s'')^{-1}\Omega
\bigl(a_.;(F_.)_{OG_n''V}\bigr).$$

\smallskip
The even orthogonal analog of Theorem 9.1 reads as follows:
\smallskip
\proclaim{\bf Theorem 9.6}
If $X$ is a connected pure-dimensional Cohen-Macaulay scheme over a field of
characteristic different from 2 and the subscheme
$$
D^k=\bigl\{x\in X | \ dim(E\cap F)_x\ge k \bigr\},
$$
defined for $k$ such that $k\equiv dim(E\cap F)_x \ (mod \ 2)$ where
$x\in X$,
is either empty or is of pure
codimension $k(k-1)/2$
in $X$, then the class of $D^k$ in the Chow group of $X$ equals
$$
\Bigl(\sum \widetilde P_IE\hak\cdot \widetilde
P_{\rho_{k-1}\smallsetminus I}F\hak\Bigr)\cap[X],
$$
where the sum is over all strict partitions $I\subset\rho_{k-1}$.
\endproclaim

\smallskip
\example{\bf Example 9.7}
The expressions giving the classes for successive $k$ are:
{\parindent=30pt
\smallskip
\item{k=1\ \ } $1$;
\smallskip
\item{k=2\ \ } $\widetilde P_1E\hak+\widetilde P_1F\hak$;
\smallskip
\item{k=3\ \ } $\widetilde P_{21}E\hak+\widetilde P_2E\hak\cdot \widetilde
P_1F\hak+
               \widetilde P_1E\hak\cdot \widetilde P_2F\hak+\widetilde
P_{21}F\hak$;
\smallskip
\item{k=4\ \ } $\widetilde P_{321}E\hak+\widetilde P_{32}E\hak\cdot \widetilde
P_1F\hak+
               \widetilde P_{31}E\hak\cdot \widetilde P_2F\hak+\widetilde
P_{21}E\hak\cdot \widetilde P_3F\hak+
	       \widetilde P_3E\hak\cdot \widetilde P_{21}F\hak+   \hfill\break
	       \widetilde P_2E\hak\cdot \widetilde P_{31}F\hak+\widetilde
P_1E\hak\cdot \widetilde P_{32}F\hak+
	       \widetilde P_{321}F\hak$.
\par}

\endexample

\medskip
\noindent
{\bf Remark 9.8.} All the formulas stated in this section in the Chow
groups have their direct analogs in topology. Perhaps the simplest version
is the following. Assume that $X$ is a compact complex manifold, the bundles
$E$, $F_i$ are holomorphic and the morphism $s$ from $X$ to
$LG_nV$
above is transverse to the smooth locus of the Schubert
variety $\Omega(a.;F.)$. Then the cohomology fundamental classes of $D(a.)$
are evaluated by the corresponding (given above) expressions in the
Chern classes of $E$ and $F_i$. The same applies to the orthogonal case.

\bigskip\smallskip

\centerline {\bf Appendix A : Quaternionic Schubert calculus}
\medskip

Let $\Bbb H$ denote the (skew) field of quaternions.
Let $\Bbb P^{\strut n}_{\Bbb H}$ be the projective space that is identified
with $(\Bbb H^{n+1} \smallsetminus \{0\}) / \thicksim$, where
$(h_1,\ldots,h_{n+1})\thicksim(h_1',\ldots,h_{n+1}')$ iff there is
$0\ne h\in \Bbb H$ such that $h_i=h\cdot h_i'$ for every $i$.
It is a compact, oriented manifold over $\Bbb R$ of dimension $4n$.
Let us recall after Hirzebruch [H1], that, in general, this real manifold
does not admit a
structure of a complex analytic manifold.

Let $G_k(\Bbb H^{\strut n})$ be the set of all $k-$dimensional subspaces
\footnote{\eightrm \ the word "(sub)space" means always a
"left $\Bbb H$-(sub)space".}
of $\Bbb H^{\strut n}$.
$G_k(\Bbb H^{\strut n})$ has a natural structure of $4k(n-k)$-dimensional,
compact, oriented manifold over $\Bbb R$.
Of course $G_1(\Bbb H^{\strut n+1})=\Bbb P_{\Bbb H}^{\strut n}$.
\smallskip

Let \ $Fl_{k_1, ...,k_r}(\Bbb H^n)$ \ be the set of all \ flags of
subspaces \ of \ consecutive \ dimensions

\noindent
$(k_1,\ldots, k_r)$ over $\Bbb H$.
It is also a compact, oriented manifold over $\Bbb R$.
One has (see [B], [Sl]),
$Fl_{k_1, ...,k_r}(\Bbb H^n) =
Sp(n)/\prod_{i=0}^r Sp(k_{i+1}-k_i)$ \ (here,
$k_0=0$ and $k_{r+1}=n$).
Of course $Fl_{k_1}(\Bbb H^n)=G_{k_1}(\Bbb H^n)$.

\smallskip

\proclaim{\bf 10.1}{\rm ([B, 31.1 p.202]) }
Let $y_1,\ldots,y_n$ be a sequence of independent variables with $deg\ y_i=4$.
Then
$$
H^*\Bigl(Fl_{k_1, ...,k_r}(\Bbb H^n),\Bbb Z\Bigr)\cong
S\Cal P(y_1,\ldots,y_n)/I_{k_1,\ldots,k_r},
$$
where $I_{k_1, \ldots, k_r}$ is the ideal generated by polynomials symmetric
in each of the sets
\hfill\break
$\{y_{k_i+1},\ldots,y_{k_{i+1}}\}$, \ $i=0,1,...,r$, \
separately \ ($k_0=0$, $k_{r+1}=n$).
\endproclaim
\smallskip

For instance (all cohomology groups are taken with coefficients in $\Bbb Z$),
$$
H^*\bigl(\Bbb P_{\Bbb H}^{\strut n}\bigr)=
\Bbb Z[y]/(y^{n+1}), \ deg\ y=4;
$$
$$
H^*\Bigl(G_k\bigl(\Bbb H^{\strut n}\bigr)\Bigr)=
S\Cal P(y_1,\ldots,y_n)/I_{k},\  deg\ y_i=4.
$$
We see that these cohomology rings are double-degree isomorphic with the
cohomology rings of their complex analogues.

Fix now a flag $V.\ :\ V_1\subset V_2\subset\ldots\subset V_n$ of subspaces
of $\Bbb H^n$
with $dim_{\Bbb H}V_i=i$.
For every partition $I\subset (n-k)^k$ we set
$$
\sigc(I)=\Bigl\{L\in G_k(\Bbb H^{\strut n})\ | \ dim_{\Bbb H}
(L\cap V_{n-k+p-i_p})=p \ , \ p=1,\ldots,k \Bigr\}.
$$
The so defined $\sigc(I)\ \bigl(I\subset(n-k)^k\bigr)$ give a cellular
decomposition of $G_k(\Bbb H^{\strut n})$ and the codimension of $\sigc(I)$
is $4|I|$. Now define

\noindent
$$
\sigma(I)=\sigma(I,V.)=
\Bigl\{L\in G_k(\Bbb H^{\strut n})\ | \ dim_{\Bbb H}
\bigl(L\cap V_{n-k+p-i_p}\bigr)\ge p \ , \ p=1,\ldots,k \Bigr\}.
$$
The cohomology classes of $\sigma(I,V.)$, in fact, do not depend on the
flag V. chosen
and will be denoted by the same symbol $\sigma(I)$.
We record:
\smallskip

\proclaim{\bf 10.2}
(Pieri-type formula) \ In $H^*(G_k(\Bbb H^n)$ one has
$$\sigma(I)\cdot\sigma(r)=\sum\sigma(J),$$ where the sum is
over $J$ such that $i_p\le j_p\le i_{p-1}$ and $|J|=|I|+r$.
\endproclaim
\smallskip

Not all proofs of the Pieri formula for Complex Grassmannians can be
extended to the quaternionic case. However, the proof
in [G-H, pp.198-204] has this advantage. As a matter of fact,
$G_k(\Bbb H^{\strut n})$ is an oriented compact manifold and thus its
cohomology ring is endowed with
the Poincar\'e duality. Moreover, one checks by direct
examination that
$$
\sigma(I)\cdot\sigma\bigl(n-k-i_k,\ldots,n-k-i_1\bigr)=\sigma\bigl((n-k)^k
\bigr)=[pt].
$$
Then the proof in loc.cit. goes through mutatis mutandis also in the
quaternionic case.

\smallskip
We can restate this information about the multiplicative structure
in $H^*\bigl(G_k(\Bbb H^{\strut n})\bigr)$ as follows:

\smallskip

\proclaim{\bf 10.3}
Let $Y=(y_1,\ldots,y_k)$ be independent variables of degree 4.
The assignment \break
$s_I(y_1,\ldots,y_k)\mapsto\sigma(I)$ for $I\subset(n-k)^k$,
and 0 -otherwise, is a ring homomorphism, and allows one
to identify $H^*\bigl(G_k(\Bbb H^{\strut n})\bigr)$ with a quotient of $S\Cal
P(Y)$
modulo the ideal $\oplus \Bbb Z s_I(Y)$, the sum over $I\not\subset(n-k)^k$.
\endproclaim
\smallskip

This result has a number of useful consequences. For example, it implies
immediately
that the signature of the Complex Grassmannian
(see [H, p.163] and [H-S, Formula (23) p.336]
is the same as the
one of the Quaternionic Grassmannian - a result proved originally
in \cite {Sl} using different methods.
\smallskip

We now describe a certain fibration which makes the Quaternionic
Grassmannians useful in study of the Grassmannians of non-maximal Lagrangian
subspaces (which are not Hermitian symmetric spaces).

Let $V=\Bbb C^{\strut 2n}$ be endowed with a nondegenerate symplectic
form $\Phi$ given by the matrix
\smallskip
$$
A=\left( \CD
0&\ \ &I_n\\
-I_n&\ \ &0
\endCD\right).
$$
where $I_n$ is the $(n\times n)$-identity matrix.
\smallskip
\noindent
Having in mind the standard notation associated with
$\Bbb H$
we endow $V$ with a structure of
$\Bbb H$-space setting  {\bf j}$\cdot v= A\overline{v}$, where
"$ \ ^{-} \ $" denotes the complex conjugation (note
that $A^2=-id_V$).
\smallskip

\proclaim{\bf 10.4}
If $U\subset V$ is $k-$dimensional Lagrangian $\Bbb C$-subspace of $V$ then
$dim_{\Bbb H}({\Bbb H}\cdot U)=k$.
Moreover, the restriction of the symplectic form $\Phi$ to any $\Bbb
H$-subspace
of $V$, is nondegenerate.
\endproclaim

To show this consider the standard Hermitian scalar product $<,>$ on
$V=\Bbb C^{\strut 2n}$.
Now given $U$, we take a $\Bbb C-$basis $u_1,\ldots,u_k$ such that
$<u_p,u_q>=\delta_{p,q}$.
We claim that $u_1,\ldots,u_k,\hbox{\bf j}u_1,\ldots,\hbox{\bf j}u_k$ are
linearly independent over $\Bbb C$ (which implies
$dim_{\Bbb H}(\Bbb H\cdot U)=k)$.
This claim follows immediately from
$\Phi(u_p,u_q)=0=\Phi(\hbox{\bf j}u_p,\hbox{\bf j}u_q)$ and
$\Phi(u_p,\hbox{\bf j}u_q)=u^t_pA(A\overline{u}_q)=-<u_p,u_q>=-\delta_{p,q}$.

\smallskip
Suppose now a $\Bbb H-$subspace $W\subset V$ is given with
$dim_{\Bbb H}W=k$, say.
We can always find $\Bbb C-$linearly independent vectors
$w_1,\ldots,w_k\in W$ such that $\Phi(w_p,w_q)=0$ and $<w_p,w_q>=\delta_{p,q}$.
Then $\hbox{\bf j}w_1,\ldots,\hbox{\bf j}w_k$ also belong to $W$.
It follows from $\Phi(w_p,w_q)=0=\Phi(\hbox{\bf j}w_p,\hbox{\bf j}w_q)$ and
$\Phi(w_p,\hbox{\bf j}w_q)=-\delta_{p,q}$ that
$w_1,\ldots,w_k,\hbox{\bf j}w_1,\ldots,\hbox{\bf j}w_k$ form a $\Bbb C-$basis
of $W$ and the form $\Phi$ restricted to $W$ is nondegenerate.

\smallskip

We infer from the above
\medskip
\proclaim{\bf 10.5}
The assignment $U\mapsto\Bbb H\cdot U$, defines a locally trivial fibration of
$LG_k\bigl(\Bbb C^{\strut 2n}\bigr)$ over $G_k\bigl(\Bbb H^{\strut n}\bigr)$
with the fiber $LG_k\bigl(\Bbb C^{\strut 2k}\bigr)$.
\endproclaim

\smallskip

In other words, denoting by $S$ the tautological (sub)bundle
over $G_k\bigl(\Bbb H^{\strut n}\bigr)$,
$rank_{\Bbb H}S=k$,
we have an identification
$LG_k\bigl(\Bbb C^{\strut 2n}\bigr)\cong LG_k(S)$, where the latter
symbol denotes (the total space of) the corresponding Grassmannian bundle.
\smallskip

This identification can be used in reduction of some problems about
Grassmannians of non-maximal Lagrangian subspaces to the problems about
the Grassmannians of maximal ones.
For example, we get from 10.5 the following identity of Poincar\'e series:
$$
P_{LG_k(\Bbb C^{2n})}(t)=
P_{G_k(\Bbb H^{n})}(t)\cdot
P_{LG_k(\Bbb C^{2k})}(t),
$$
thus reproving the result from \cite {P-R2, Corollary 1.7}.
\bigskip

Similar fibrations exist for flag varieties. Let $LFl_{k_1,\ldots, k_r}
(\cc^{2n})$ be the variety paramet-
\noindent
rizing Lagrangian ( w.r.t. $\Phi$ ) flags of
dimensions $(k_1,\ldots, k_r)$ in $\cc^{2n}$.
\medskip

\proclaim{\bf 10.6}
The assignment \ ($dim_{\cc}U_i=k_i, \ i=1,\ldots, r$):
$$\Bigl(U_1\subset U_2\subset \ldots \subset U_r\Bigr) \mapsto
\Bigl(\Bbb H\cdot U_1 \subset \Bbb H\cdot U_2 \subset \ldots
\subset \Bbb H\cdot U_r\Bigr)
$$
is a locally trivial fibration of $LFl_{k_1,\ldots, k_r}(\cc^{2n})$ over
$Fl_{k_1,\ldots, k_r}(\Bbb H^{n})$. If $\Bbb C^{\strut 2k_1}\subset
\Bbb C^{\strut 2k_2}\subset\ldots\subset\Bbb C^{\strut 2k_r}$ is
a (part of) the standard flag, then the fiber of this fibration is the
variety parametrizing Lagrangian flags $W_1\subset W_2\subset\ldots\subset W_r$
such
that $W_i\subset\Bbb C^{\strut 2k_i}$
and $dim_{\Bbb C}W_i=k_i,\ i=1,\ldots,r$.
\endproclaim
\smallskip

Therefore the fiber is a composition of Lagrangian Grassmannian
bundles of maximal subspaces.
In particular, we obtain the following formula for the
Poincar\'e series of $LFl_{k_1,\ldots,k_r}(\Bbb C^{\strut 2n})$:
$$
P_{LFl_{k_1,\ldots,k_r}(\Bbb C^{2n})}(t)=
P_{Fl_{k_1,\ldots,k_r}(\Bbb H^{n})}(t)\cdot\prod\limits^r_{i=1}
P_{LG_{k_i-k_{i-1}}(\Bbb C^{2(k_i-k_{i-1})})}(t),
$$
where $k_0=0$. Since explicit expressions for the factors on the right-hand
side
are known (see (10.1)), this gives an explicit formula for
$P_{LFl_{k_1,\ldots,k_r}(\Bbb C^{2n})}(t)$.
\smallskip

\noindent
{\bf 10.7.} Finally, we show an algebro-topological interpretation
(as well as another proof) of the identity:
$$
s_I(x^2_1,\ldots,x^2_n)\cdot s_{\rho_{n}}(x_1,\ldots,x_n)=
s_{2I+\rho_{n}}(x_1,\ldots,x_n)
$$
from Section 5.
To this end we show two different ways of constructing
$LFl:=LFl\bigl(\Bbb C^{\strut 2n}\bigr)$.
The first way is given by taking the total space of the flag bundle
$Fl(R)\to LG_n\bigl(\Bbb C^{\strut 2n}\bigr)$ where $R$ is the
tautological vector bundle on $LG_n\bigl(\Bbb C^{\strut 2n}\bigr)$.
The second way relies on the following observation:
$LFl$ can be interpreted as the variety of flags
$W_1\subset W_2\subset\ldots\subset W_{2n}$ such that $dim_{\Bbb C}W_j=j$
and each $W_{2j}$ is a $\Bbb H$-subspace.
This realization is given by the assignment:
$$
\Bigl(V_1\subset V_2\subset\ldots\subset V_n\Bigr)\mapsto
\Bigl(V_1\subset\Bbb H\cdot V_1\subset\Bbb H\cdot V_1+V_2\subset
\Bbb H\cdot V_1+\Bbb H\cdot V_2\subset\ldots\Bigr)
$$
Equivalently, using the tautological sequence
$S_1\subset S_2\subset\ldots\subset S_n$, $rank_{\Bbb H}S_i=i$, on
$Fl_{\Bbb H}$, this corresponds to taking the total space of the product
of projective bundles
$$
\Bbb P:=\Bbb P(S_2/S_1)\times_{Fl_{\Bbb H}}\ldots\times_{Fl_{\Bbb H}}
\Bbb P(S_n/S_{n-1})\to Fl_{\Bbb H}
$$
where $S_{i+1}/S_i$, $i=1,\ldots,n$, are considered as rank $2$
complex  bundles.
\smallskip
The same holds in the relative situation, i.e. given a rank $2n$ vector
bundle $V\to X$ endowed with a symplectic form we get a commutative diagram

$$
\CD
\Bbb P &=\!= &LFl(V) &=\!= &Fl(R)\\
@V\tau_1VV & & @V\pi_1VV\\
Fl_{\Bbb H}(V) @>\tau_2>> X @<\pi_2<< LG_nV
\endCD
$$
\smallskip
\noindent
where $Fl_{\Bbb H}(V)$ is the Quaternionic (complete) flag bundle.
Let $x_1,\ldots,x_n$ be the sequence of the Chern roots of the tautological
quotient bundle on $LG_nV$.
By Corollary 5.6(i) we know that if there exists an even $i_p$, then
$(\pi_2\circ\pi_1)_*\bigl(x_1^{i_1}\cdot\ldots\cdot x_n^{i_n}\bigr)=0$.
(Calculating the other way arround, this follows easily from the
projection formula.)
On the other hand, iff all $i_p$ are odd, then (see Proposition 5.5)
$$
s_{\rho_n}(x_1,\ldots,x_n)\cdot
(\pi_2\circ\pi_1)_*\bigl(x_1^{i_1}\cdot\ldots\cdot x_n^{i_n}\bigr)=
s_{I- \rho_{n-1}}(x_1,\ldots,x_n).
$$

\noindent
Putting $i_p=2j_p+1$ and calculating the other way around, we get

$$
\aligned
(\tau_2\tau_1)_*&\Bigl(x_1^{2j_1+1}x_2^{2j_2+1}\ldots x_n^{2j_n+1}\Bigr)=\\
&=(\tau_2)_*\Bigl( (x_1^2)^{j_1}\cdot(x_2^2)^{j_2}\cdot\ldots\cdot
      (x_n^2)^{j_n}\Bigr)\\
&=s_{J-\rho_{n-1}}\bigl(x_1^2,\ldots,x_n^2\bigr).
\endaligned
$$

\noindent
Indeed, recalling the notation from 10.1 we have $y_p=x_p^2, \ p=1,\ldots,n$
(see [B, 31.1]), and we use the fact that $(\tau_2)_*$ is induced by
the Jacobi symmetrizer (recalled in the proof of Corollary 5.6(ii)
and that of Lemma 5.7(ii) )
this time applied to $y_1,\ldots,y_n$.
The latter statement follows from 10.1 by exactly the same reasoning as
that used in the proof of Lemma 2.4 in [P1].
Comparison of the results of both computations, yields the desired
identity.
\bigskip\smallskip

\centerline {\bf Appendix B : Introduction to Schubert polynomials
\`a la polonaise}

\medskip
We provide here a brief sketch of a theory of symplectic
Schubert polynomials which has grown up from the present work.
For details and further developments as well as for the orthogonal
Schubert polynomials, we refer the reader to [L-P-R].

Let $(x_1,x_2,\ldots)$ be a sequence of independent variables.
Let $w_0$ be the longest element in the Weyl group $W_n$ of type $C_n$\,.
Define
$$
\Cal C_{w_0}:=\Cal C_{w_0}(x_1,\ldots,x_n):=
(-1)^{n(n-1)/2}\,x_1^{n-1}\,x_2^{n-2}\,\ldots\, x^1_{n-1}\,x^0_n\,
\widetilde Q_{\rho_n}(x_1,\ldots,x_n)\,,
$$
and for an arbitrary $w\in W_n$\,,
$$
\Cal C_w:=\Cal C_w(x_1,\ldots,x_n):=\partial_{w^{-1}w_0}'(\Cal C_{w_0})\,.
$$
Above, by $\partial'_w$ $(w\in W_n)$ we understand the composition of the
divided difference operators $\partial'_i$ defined by
$$\aligned
\partial'_0(f) & =(f-s_0f)/2x_1\\
\partial'_i(f) & =(f-s_if)/(x_{i+1}-x_i)\qquad i=1,2,\ldots,n-1\,,
\endaligned$$
associated in a usual way with an arbitrary reduced decomposition of $w$
using $s_i$\,, $i=0,1,\ldots,n-1$\,.

These polynomials satisfy the following properties.
\smallskip

{\parindent=20pt
\item{1.} (Stability) Suppose that $m>n$\,.
Let $W_n\hookrightarrow W_m$ be the embedding via the first $n$ components.
Then, for any $w\in W_n$, the following equality holds:
$$
\Cal C_w(x_1,\ldots,x_m)|_{x_{n+1}=\ldots=x_m=0}=\Cal C_w(x_1,\ldots,x_n)\,.
$$
\smallskip
\item{2.} (the Grassmannian case) Let $I=(i_1>\ldots>i_k>0)$ be a strict
partition contained in $\rho_n$\,. Set
$$
w_I=(\overline{i_1},\ldots,\overline{i_k},j_1<j_2<\ldots<j_{n-k})\,,
$$
where $\{i_1,\ldots,i_k,j_1,\ldots,j_{n-k}\}=\{1,2,\ldots,n\}$\,. Then
$$
\Cal C_{w_I}(x_1,\ldots,x_n)=\widetilde Q_I(x_1,\dots,x_n)\,.
$$
\par}
\smallskip

As we know from Section 4, $\widetilde Q_I(x_1,\dots,x_n)$
is a positive sum of monomials.
The polynomial $\Cal C_w$ has not this property.
Also, it is in general neither negative nor positive sum of monomials.

The following is the list of symplectic Schubert polynomials for $n=2$.
$$\aligned
\Cal C_{(\overline{1},\overline{2})} = & -x_1^3x_2-x_1^2x_2^2 \\
\Cal C_{(1,\overline{2})}=-x_1^2x_2\qquad
& ,\qquad \Cal C_{(\overline{2},\overline{1})}=x_1^2x_2+x_1x_2^2 \\
\Cal C_{(\overline{2},1)}=x_1x_2\qquad
& ,\qquad \Cal C_{(2,\overline{1})}=x_2^2 \\
\Cal C_{(2,1)}=x_2 \qquad & ,\qquad \Cal C_{(\overline{1},2)}=x_1+x_2 \\
\Cal C_{(1,2)} & =1\,.
\endaligned$$

\bye